\documentclass{article}

\usepackage{lipsum}
\usepackage{amsfonts}
\usepackage{graphicx}
\usepackage{epstopdf}
\usepackage{algorithmic}
\usepackage{enumerate}
\usepackage{csquotes}
\usepackage{todonotes}
\usepackage{subcaption}
\usepackage{wrapfig}
\ifpdf
\DeclareGraphicsExtensions{.eps,.pdf,.png,.jpg}
\else
\DeclareGraphicsExtensions{.eps}
\fi

\usepackage{soul}

\usepackage{algorithm}
\usepackage{algorithmic}
\usepackage{amsmath}
\usepackage{amsthm}
\usepackage{amssymb}
\usepackage{booktabs}
\usepackage{hyperref}
\usepackage{cleveref}
\usepackage{url}
\usepackage[numbers]{natbib}
\usepackage[margin=1.4in]{geometry}

\usepackage{lipsum}
\usepackage{amsfonts}
\usepackage{amssymb} 
\usepackage{graphicx}
\usepackage{epstopdf}
\usepackage{algorithmic}
\usepackage{bbm}
\usepackage{mathrsfs}
\usepackage{todonotes}
\usepackage{csquotes}
\usepackage{subcaption}
\usepackage{booktabs}
\ifpdf
  \DeclareGraphicsExtensions{.eps,.pdf,.png,.jpg}
\else
  \DeclareGraphicsExtensions{.eps}
\fi

\usepackage[shortlabels]{enumitem}
  \setlist[enumerate]{leftmargin=.5in}
  \setlist[itemize]{leftmargin=.5in}

\title{Regularising orientation estimation in Cryo-EM 3D map refinement through measure-based lifting over Riemannian manifolds\thanks{Our code for the main algorithms is available at \href{https://github.com/wdiepeveen/Cryo-EM}{\texttt{https://github.com/wdiepeveen/Cryo-EM}} .
}
}

\author{Willem Diepeveen\thanks{Faculty of Mathematics, University of Cambridge, Cambridge 
  (wd292@cam.ac.uk, cbs31@cam.ac.uk).}
\and Jan Lellmann\thanks{Institute of Mathematics and Image Computing, University of L\"ubeck, 
  (lellmann@mic.uni-luebeck.de).}
\and Ozan \"Oktem\thanks{Department of Mathematics, KTH–Royal Institute of Technology (ozan@kth.se)}
\and Carola-Bibiane Sch\"onlieb\footnotemark[2]}
\date{}
\usepackage{amsopn}

\newcommand{\revA}[1]{#1}

\newcommand{\Natural}{\mathbb{N}}

\newcommand{\Real}{\mathbb{R}}

\newcommand{\Sphere}{\mathbb{S}}

\newcommand{\DimInd}{d} 
\newcommand{\probSpace}{\mathbb{P}} 
\newcommand{\ballVol}{\omega} 

\newcommand{\eA}{x} 
\newcommand{\eB}{y} 
\newcommand{\eC}{z} 
\newcommand{\ePoint}{x} 
\newcommand{\dom}{\Omega} 
\newcommand{\coordInd}{i} 
\newcommand{\eball}{\mathbb{B}} 

\newcommand{\manifold}{\mathcal{M}} 
\newcommand{\tangent}{\mathcal{T}} 
\newcommand{\mPoint}{p} 
\newcommand{\mPointB}{q} 
\newcommand{\mTVector}{X} 
\newcommand{\mTVectorB}{Y} 
\newcommand{\mPoints}{P} 
\newcommand{\chart}{\phi} 
\newcommand{\differential}{D} 
\newcommand{\ball}{\mathcal{B}} 
\newcommand{\ellipsoid}{\mathcal{E}} 

\newcommand{\geodesic}{\gamma} 
\newcommand{\distance}{d} 

\newcommand{\Grad}{\operatorname{grad}} 
\newcommand{\Hess}{\operatorname{Hess}} 

\newcommand{\dummyFunctionA}{F} 
\newcommand{\dummyFunctiona}{f} 

\newcommand{\dummyBilForm}{Q} 

\newcommand{\rot}{g} 

\newcommand{\optim}[1]{{#1}^{*}} 

\newcommand{\PDens}{\mathsf{P}}
\newcommand{\stochastic}[1]{\mathbbm{#1}}

\newcommand{\measure}{\mu} 
\newcommand{\dirac}{\delta} 
\newcommand{\wasserstein}{W} 

\newcommand{\indicator}{\chi} 

\newcommand{\vol}{u} 
\newcommand{\volSpace}{\mathscr{X}} 
\newcommand{\volLtwoSpace}{L^2(\dom)} 

\newcommand{\waveNumber}{k} 
\newcommand{\defocus}{\Delta z} 
\newcommand{\spAberration}{C_s} 
\newcommand{\ampContrast}{\alpha} 

\newcommand{\noiseLevelParam}{\sigma} 
\newcommand{\volRegParam}{\tau} 

\newcommand{\tim}{t} 

\newcommand{\forward}{\mathcal{A}} 
\newcommand{\projection}{\mathcal{P}} 
\newcommand{\RerCTF}{H} 
\newcommand{\ReCTF}{\hat{H}} 
\newcommand{\fourier}[1]{\widehat{#1}} 
\newcommand{\rampKernel}{\mathcal{K}}

\newcommand{\regA}{\mathcal{R}} 
\newcommand{\regParam}{\gamma} 

\newcommand{\rotsDensCoef}{\beta} 
\newcommand{\simplex}{\Delta} 

\newcommand{\img}{f} 
\newcommand{\stimg}{\stochastic{\img}}
\newcommand{\imgSpace}{\mathscr{Y}} 
\newcommand{\numImgs}{N} 
\newcommand{\imgInd}{i} 
\newcommand{\noise}{e}
\newcommand{\stnoise}{\stochastic{\noise}}
\newcommand{\PDDataNoise}{\PDens_{\text{noise}}}

\newcommand{\rotsSampling}{\mathcal{X}} 
\newcommand{\numRots}{|\rotsSampling|} 

\newcommand{\eigenValue}{\lambda} 

\newcommand{\iterInd}{k} 
\newcommand{\iterIndB}{l} 
\newcommand{\sumIndA}{i} 
\newcommand{\sumIndB}{j} 
\newcommand{\sumTotA}{m} 
\newcommand{\exponentA}{\eta} 
\newcommand{\radiusA}{R} 
\newcommand{\radiusB}{r} 
\newcommand{\radiusC}{\rho} 
\newcommand{\dummyIntegerB}{J} 

\newcommand{\Dim}{\operatorname{dim}}

\newcommand{\argmin}{\operatorname*{argmin}}
\newcommand{\argmax}{\operatorname*{argmax}}

\newtheorem{theorem}{Theorem}[section]

\newtheorem{lemma}[theorem]{Lemma}
\newtheorem{proposition}[theorem]{Proposition}
\newtheorem{definition}[theorem]{Definition}
\newtheorem{remark}[theorem]{Remark}

\begin{document}

\maketitle

\begin{abstract}
	\revA{Motivated by the trade-off between noise-robustness and data-consistency for joint 3D map reconstruction and rotation estimation in single particle cryogenic-electron microscopy (Cryo-EM), w}e propose ellipsoidal support lifting (ESL), a \revA{measure}-based \revA{lifting} scheme for \revA{regularising and} approximating the global minimiser of a smooth function over a Riemannian manifold. 
	Under a uniqueness assumption on the minimiser we show several theoretical results, in particular \revA{well-posedness of the method and }an error bound \revA{due to the induced bias} with respect to the global minimiser. Additionally, we use the developed theory to integrate the \revA{measure}-based \revA{lifting} scheme into an alternating \revA{update} method for joint homogeneous \revA{3D map} reconstruction and rotation estimation, where typically tens of thousands of manifold-valued minimisation problems have to be solved \revA{and where regularisation is necessary because of the high noise levels in the data}. The joint recovery method
	is used to test both the theoretical predictions and algorithmic performance through numerical experiments with Cryo-EM data. \revA{In particular, the induced bias due to the regularising effect of ESL empirically estimates better rotations, i.e., rotations closer to the ground truth, than global optimisation would.}

\end{abstract}



\section{Introduction}
\label{sec:introduction}
\revA{Cryogenic-electron microscopy (Cryo-EM) single particle analysis is a tomographic imaging technique for determining the 3D structure (map) of a biomolecule from multiple 2D transmission electron microscope (TEM) images, each representing the \enquote{projection} of the biomolecule in some unknown rotation (\cite[Figure~1]{singer2018mathematics}). Jointly recovering the volumetric 3D map of the biomolecule along with particle specific rotation from the 2D TEM images is typically approached through an alternating update scheme. These alternating schemes typically intertwine iterates for updating the 3D map given the rotations with iterates for updating the rotations given the 3D map. Arguably, the challenging part in such an \emph{alternating joint refinement} approach is the estimation of rotations given the high noise-levels that data has in Cryo-EM (\cite[Figure~7]{bendory2020single}). 

To see why joint refinement requires particular care, consider
the family of loss functions $\{\dummyFunctionA_\imgInd^{cryo}\}_{\imgInd=1}^\numImgs$ with $\dummyFunctionA_\imgInd^{cryo}:\mathrm{SO}(3)\to \Real$ given by
\begin{equation}
    F_{\imgInd}^{cryo}(\rot) := \frac{1}{2\sigma}\|\forward (\rot . \vol^\iterInd) - \img_\imgInd \|_2^2, \quad \text{for $\imgInd=1,\ldots, \numImgs$.}
    \label{eq:intro-cryo-rotation-loss}
\end{equation}
In the above, $\{\img_\imgInd\}_{\imgInd=1}^\numImgs \subset \imgSpace$ are the noisy particle specific 2D TEM (projection) images that reside in some function space $\imgSpace$, $\vol^\iterInd\subset \volSpace$ is the 3D map at iteration $\iterInd$ of a joint refinement scheme residing in some function space $\volSpace$, $\forward:\volSpace\to\imgSpace$ is the Cryo-EM forward operator, $\rot . \vol(x):= \vol (\rot^{-1} x)$ and $\noiseLevelParam$ is a parameter representing the noise level. 
Say we want to estimate the rotation that fits the data best in the sense that it yields a low loss. One way of going about that is to seek the global minimiser to the functions $\dummyFunctionA_\imgInd^{cryo}$ for all $\imgInd=1,\ldots,\numImgs$. Because of the non-convexity of these loss functions and the large number, $\numImgs$, of optimisation problems, in practice one would take the point with lowest value over a properly constructed sampling set, e.g., through branch-and-bound as done in the popular software package cryoSPARC \cite{punjani2017cryosparc}. 

\begin{figure}[h!]
    \centering
    \includegraphics[width=0.49\linewidth]{"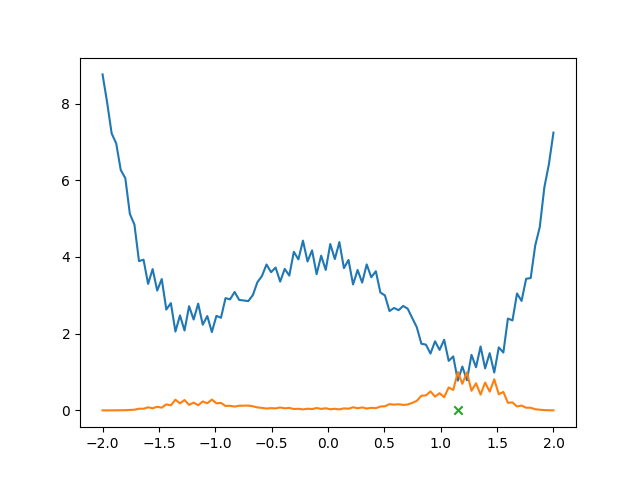"}
    \caption{\revA{A smooth toy loss function $\dummyFunctionA^{toy}:(-2,2)\to\Real$ evaluated on a coarse grid (blue). Because of the irregular loss landscape there is an ambiguity regarding what input variable would be an objectively better choice. Two approaches are shown: soft assignment through computing a probability measure (orange) or hard assignment through picking the point with lowest loss (green).}}
    \label{fig:1D-examples}
\end{figure}
These loss functions change after the 3D map is updated, the latter through fixing the rotations and optimising across all 2D TEM images. The non-convexity, the high noise levels in data, and the increasing resolution of the 3D map results in an increasingly irregular optimisation landscape. Hence, such a strategy potentially has a hard time converging in practice, i.e., there is \emph{lack of robustness to noise}. In addition, in practice, the ground truth rotation hardly ever coincides with the global minimiser, although they tend to be close to each other. This begs the question whether optimisation is a good approach in the first place. A popular -- more noise-robust -- approach is to marginalise over the space of rotations, e.g., by assigning a probability to every rotation in the discretized space based on the loss at that rotation and average over rotations against that probability measure in the 3D map update step. The latter underlies the popular software package RELION \cite{scheres2012relion}. Such a method tends to be more robust to noise as the rotations with low losses empirically tend to stay at least in the same neighbourhood.  \cref{fig:1D-examples} shows these two approaches for a 1D toy example with a irregular and undersampled loss landscape.


Whereas the marginalisation approach has a clear advantage in the rotation update step, the optimisation approach has an advantage in the 3D map update step. For the latter approach, consider the 3D map update through
\begin{equation}
    \inf_{\vol\in\volSpace} \Bigl\{\frac{1}{2\sigma} \sum_{\imgInd=1}^\numImgs \|\forward (\rot_\imgInd^{\iterInd} . \vol) - \img_\imgInd \|_2^2 + \regA(\vol) \Bigr\},
    \label{eq:data-consistent-optimisation}
\end{equation}
whereas, the former would update the 3D map as
\begin{equation}
    \inf_{\vol\in\volSpace} \Bigl\{\frac{1}{2\sigma} \sum_{\imgInd=1}^\numImgs \int_{\mathrm{SO}(3)} \|\forward (\rot . \vol) - \img_\imgInd \|_2^2 \; \mathrm{d}\measure_\imgInd^{\iterInd}(\rot)+ \regA(\vol) \Bigr\},
\end{equation}
where $\regA:\volSpace\to\Real$ is some suitably chosen regularisation functional. The solutions to these 3D map optimisation problems are regularised back-projections \cite{scheres2012relion}. So for the marginalisation approach, back-projecting each image for a multiple rotations (each weighted by the probability measure) \enquote{smears out the data}, causing a \emph{loss of data-consistency} in the final 3D map. The optimisation approach does not suffer from this drawback assuming the rotations are estimated correctly. 

The methods described above requires one to choose between noise-robustness and data-consistency. In this paper we develop a method that combines the best of both worlds by being \enquote{in between optimisation and marginalisation}. To concretise, consider as a first step the following:
\begin{equation}
    \optim{\rot}_\imgInd \in \argmin_{\rot\in \mathrm{SO}(3)} \dummyFunctionA_\imgInd^{\text{cryo}} (\rot) \quad \Longleftrightarrow \quad \dirac_{\optim{\rot}_\imgInd} \in \argmin_{\measure\in \mathbb{P}(\mathrm{SO}(3))} \int_{\mathrm{SO}(3)} \dummyFunctionA_\imgInd^{\text{cryo}} (\rot) \; \mathrm{d}\measure.
\end{equation}
The right hand setting allows us to move from optimisation towards marginalisation, e.g., by adding regularisation on the measure. Doing this can yield more noise-robustness, but one loses data-consistency for reasons mentioned above. However, similarly to the toy example in \cref{fig:1D-examples}, most of the weight of the measure in Cryo-EM refinement has very local support and this behaviour is typically enforced through a local search routine later on in the refinement \cite{scheres2012relion}. Furthermore, integrating a function with measure $\measure$ having local support can be approximated by function evaluation at the barycentre of that measure, i.e., for the $\imgInd$'th image
\begin{equation}
    \int_{\mathrm{SO}(3)} \dummyFunctionA_\imgInd^{\text{cryo}} (\rot) \; \mathrm{d}\measure_\imgInd(\rot) \approx \dummyFunctionA_\imgInd^{\text{cryo}} (\rot_{\measure_\imgInd}), \quad \text{where}\quad \rot_{\measure_\imgInd}\in \argmin_{\rot'\in\mathrm{SO}(3)} \int_{\mathrm{SO}(3)} \distance_{\mathrm{SO}(3)} (\rot', \rot)^2 \; \mathrm{d}\measure_\imgInd (\rot).
    \label{eq:intro-barycentre-idea}
\end{equation}
Now, the barycentre rotation $\rot_{\measure_\imgInd}$ encodes the same information as the measure $\measure_\imgInd$ without introducing any smearing of the data in the 3D map update step (which ensures data-consistency). In other words, noise-robustness is inherited from $\measure_\imgInd$ and we improve upon data-consistency. Finally, if the width of the measure can be tuned for locality, the induced bias in the barycentre can be tuned as well. Indeed as $\measure\to \dirac_{\optim{\rot}_\imgInd}$, we expect that $\rot_{\measure} \to \optim{\rot}_\imgInd$, i.e., the barycentre converges to the energy minimising rotation. Then, regularising the measure and computing the barycentre is indeed \enquote{in between optimisation and marginalisation}. 

Besides having both noise-robustness and data-consistency for the rotation estimation, a key question is whether the estimated rotations are closer to the ground truth rotations than the global minimiser. In other words, does adding bias pay off in practice? Investigating this will be a central theme throughout this work.

As a final note, observe that the problem of estimating particle specific rotations is a special case of a more general problem of recovering particle specific nuisance parameters. In single particle Cryo-EM refinement, solvers typically also account for symmetric molecules \cite{jiang2017atomic} and in-plane shift due to improper centering of the 2D TEM images. The former reduces to estimating rotations on the quotient manifold $\mathrm{SO}(3)/G$, where $G$ is the symmetry group of the molecule. The latter results in an estimation problem over the product manifold $\mathrm{SO}(3) \times \Real^2$. The framework we develop can also be used to handle this setting. We therefore consider a generic setting that amounts to regularise and solve variational problems of the following form:
\begin{equation}
\label{eq:model-problem}
    \inf_{\mPoint\in \manifold} \dummyFunctionA(\mPoint)
    \quad\text{for some smooth energy function $\dummyFunctionA \colon \manifold\to \Real$.}
\end{equation} 
Our approach is based on the framework of measure-based lifting for global optimisation on Riemannian manifolds.

}

\subsection{Related work}
\label{sec:lifting-math-relatedwork}
The paper relates mainly to two fields that have developed separately from each other, rotation estimation in Cryo-EM and solving inverse problems on Riemannian manifolds. For that reason, we will discuss these two separately.

\paragraph{Rotation estimation in Cryo-EM}
\revA{The loss function \cref{eq:intro-cryo-rotation-loss} assumes the forward problem of each $\img_{\imgInd} \in \imgSpace$ being a single sample of the $\imgSpace$–valued random variable $\stimg_{\imgInd}$ where
\begin{equation}\label{eq:CryoEMInvProb}
    \stimg_\imgInd = \forward(\rot_\imgInd . \vol) + \stnoise_\imgInd 
    \quad\text{with $\rot_\imgInd\in \mathrm{SO}(3)$ and $\stnoise_{\imgInd} \sim \PDDataNoise^{\imgInd}$ for $\imgInd=1,\ldots, \numImgs$.}
\end{equation}
The inverse problem is then the task of \emph{jointly} recovering the 3D map $\vol \in \volSpace$ of the biomolecule while estimating (or marginalising out) the unknown particle specific rotations $\rot_\imgInd \in \mathrm{SO}(3)$ for $\imgInd = 1,\ldots,\numImgs$.

Mathematically, the first  natural question that comes to mind when analysing an inverse problem is uniqueness. 
For tomographic imaging problems this typically means determining whether uniqueness holds assuming continuum noise-free data generated by a perfectly known forward operator.
For \cref{eq:CryoEMInvProb}, this means to determine whether uniqueness holds for the joint problem of recovering the 3D map and estimating the particle specific rotations. 
When the forward operator $\forward$ in \cref{eq:CryoEMInvProb} is the parallel beam ray transform, then one can, except for some special cases, prove uniqueness up to a global orthogonal transformation \cite{lamberg2008unique}. 
Interestingly, adopting a more realistic TEM imaging model (based on Born approximation and TEM optics) allows one to prove uniqueness without any global orthogonal transformation \cite{Kurlberg:2021aa}.

The next natural mathematical topic is to characterise the degree of ill-posedness of the inverse problem in \cref{eq:CryoEMInvProb}.
If particle specific rotations are known, then \cref{eq:CryoEMInvProb} reduces to inverting the parallel beam ray transform on complete data (actually data is here the entire 4D manifold of lines, so the problem is over-determined), which is known to be a mildly ill-posed problem.
As to be expected, the ill-posedness worsens when rotations are also unknown.
Even though the authors are unaware of any mathematically rigorous characterisations of the degree of ill-posedness in this setting, numerical experience strongly suggests that the inverse problem in \cref{eq:CryoEMInvProb} is severely ill-posed. 

The above theoretical issues represent a small part of the research within computational sciences inspired by single particle Cryo-EM.
The main focus within the Cryo-EM community has been on the design and implementation of algorithms for solving \cref{eq:CryoEMInvProb}. 
This inverse problem amounts to jointly solving two coupled inverse problems, namely one related to recovering the 3D map and the other on recovering particle specific rotations. 
These two inverse problems, which are coupled in \cref{eq:CryoEMInvProb}, have very different characteristics.
Hence, it was natural to seek methods that \enquote{disentangle} them.

Early approaches relied on the observation that if the forward model is the parallel beam ray transform, then one can use integral geometric techniques to estimate rotations without having access to a 3D map. 
With this disentanglement, one can first estimate particle specific rotations, then use these to recover the 3D map. Such an approach was proposed independently as early as in 1986 \cite{Vainshtein1986} and again in 1987 \cite{van1987angular} for retrieving relative rotations of 2D projections in Fourier space through common lines. Soon after, the method of moments \cite{goncharov1987methods,goncharov1988integral,goncharov1988determination,salzman1990method} was proposed and later further improved and analyzed in \cite{basu2000feasibility,basu2000uniqueness}, which did not need any information on the 3D map either. 
A different statistical approach that has recently become popular in the mathematics community is Kam's auto-correlation analysis, which dates back to 1980 \cite{kam1980reconstruction}. Here, the recovery of the 3D map bypasses the estimation of the rotations. Instead, the recovery of the 3D map is based on using the 2D image patches to compute the autocorrelation and higher order correlation functions of the 3D map in Fourier space. 

One limitation with many of the above methods relates to difficulties in handling the high noise levels present in contemporary low-dose single particle Cryo-EM data. 
The common lines method will in its standard formulation break down under high noise conditions, but several extensions to it has resulted in formulations that are more resilient to noise \cite{bandeira2020non,shkolnisky2012viewing,singer2010detecting,singer2011three,wang2013orientation}. Similarly, the method of moments and Kam's method have also served as an inspiration for several extensions \cite{bhamre2015orthogonal,bandeira2017estimation, levin20183d,sharon2020method}. 

Apart from noise, a second hurdle relates to the forward operator, and in particular to the varying \emph{contrast transfer function} (CTF) parameters (most notably the defocus) between micrographs. For this reason and other factors such as structural heterogeneity of particles, there is typically a refinement stage after the \textit{ab initio} stage as described above. 

Early work on iterative refinement can be traced back as early as 1972 \cite{gilbert1972iterative} and to follow up works such as \cite{harauz1983direct}, where typically a 3D map update step alternated with an rotation alignment step. These type of algorithms were also the cornerstone of the first generation of software packages \cite{frank1996spider,ludtke1999eman,van1996new} and most software packages being released in the years after followed this approach \cite{grigorieff2007frealign,hohn2007sparx,shaikh2008spider,sorzano2004xmipp,tang2007eman2}.

Arguably, starting from the 2010s there has been a shift in rotation estimation in refinement, and there are now two packages dominating the field: RELION \cite{scheres2012relion} and cryoSPARC \cite{punjani2017cryosparc}. Rather than picking one rotation after every alignment step -- which can suffer from high noise levels -- RELION gives probabilities to each rotation and marginalizes over a large set, which was originally proposed in \cite{sigworth1998maximum}. Although losing consistency to data, this method tends to be very robust to noise because of this averaging approach. CryoSPARC on the other hand uses a branch-and-bound algorithm for the rotation estimation, which is more sophisticated than plain alignment and faster than the RELION approach. It also solves the original problem rather than a modification, but arguably pays in noise-robustness.

Even more recently, deep learning-based methods have been introduced to joint refinement. Here too approaches can be categorized by whether they work around the rotations through marginalisation \cite{gupta2021cryogan} or try to actually find them through a multi-scale solver \cite{zhong2021cryodrgn2}.
To summarise, model- or learning-based refinement methods proposed so far are choosing between noise-robustness and data-consistency. As we will argue in this work, there is plenty to be gained here.
}

\paragraph{Manifold-valued inverse problems}
\revA{Popular approaches for solving inverse problems over Riemannian manifolds are optimisation-based. Work on optimisation problems on manifolds} similar to \cref{eq:model-problem} can at least be traced back to 1972 \cite{luenberger1972gradient}, where the paradigm was to rely on Euclidean solvers while accounting for the manifold as a constraint set in $\Real^\DimInd$ on which each iterate had to be projected throughout the optimisation scheme. Such methods are not generally applicable since \revA{it is not always possible to project an arbitrary point in the embedded space onto the} manifold, as is the case with the symmetric positive definite matrices. Also, alternatively resorting to local charts and performing the optimisation using a Euclidean solver has several drawbacks. As an example, symmetries of the manifold are poorly accounted for, distortions of the metric may corrupt convergence rates, and restricting ourselves to local neighbourhoods makes global analysis impossible \cite{grohs2016nonsmooth}. These issues can be resolved by passing to \emph{intrinsic} methods, where optimisation is formalised in terms of several structures and mappings on the manifold that do not depend on the specific embedding or chart. For this reason, the development and study of intrinsic algorithms for optimisation on manifolds has become a flourishing field, especially after seminal works in the mid-1990s \cite{smith1994optimization,udriste1994convex}. \revA{These developments have now been collected in popular works such as \cite{absil2009optimization,boumal2020introduction}.}

\revA{In inverse problems optimisation is only part of the problem. Proper regularisation of the loss functional is key to retrieving good reconstructions and has amounted to a wide variety of smooth and non-smooth models \cite{bacak2016second,bergmann2018priors,bredies2018total,holler2020non,weinmann2014total} and specialized solvers \cite{bergmann2019recent,persch2018optimization}. However, the geometry of the problem determines the degree of difficulty to finding the global minimiser to these variational models. In particular, i}n the Euclidean case, convexity is key for efficiently finding a unique global minimiser of a function and having non-smoothness is no real issue for \revA{convex optimisation} algorithms \cite{clason2020introduction}. For manifolds, only a certain subset allows for existence of convex functions. That is, on \emph{Hadamard manifolds} -- Riemannian manifolds that are complete, simply connected, and have non-positive sectional curvature -- convex \revA{functions} can be generalized and has been crucial in 
\revA{the design of the above-mentioned variational methods \cite{weinmann2014total} and in} the convergence proofs \revA{of solvers}. Indeed, besides the standard gradient-based methods \cite{zhang2016first} a myriad of algorithms exploiting special structure of the \revA{-- possibly non-smooth --} convex function of interest have been generalized from the Euclidean case only very recently \cite{bacak2014computing,banert2014backward,bergmann2015restoration,bergmann2016parallel,bergmann2021fenchel,diepeveen2021inexact}.


\revA{The case of rotation estimation in Cryo-EM joint refinement poses new challenges. One key challenge is regularisation, i.e., in the typical absence of an initial guess for the rotations, distance-based regularisers -- i.e., $\mPoint \mapsto \frac{1}{s}\distance_\manifold(\mPoint, \mPointB)^s$ generalizing $\ell^s$ regularisation for $s\geq 1$ --  are not suitable. Similarly, the problem lacks the structure needed for total variation-like regularisation. Optimisation-wise, $\mathrm{SO}(3)$ being non-Hadamard is another indication for being careful when optimising. 
Whereas there is a large amount of work on global optimisation on compact manifolds -- e.g., simulated annealing \cite{arnaudon2014stochastic,holley1989asymptotics,holley1988simulated}, genetic algorithms \cite{fong2022stochastic}, particle swarm optimisation \cite{fornasier2020consensus,chizat2021sparse}, Nelder-Mead \cite{dreisigmeyer2007direct}, mesh adaptive direct search (MADS) \cite{dreisigmeyer2007direct,dreisigmeyer2018direct} --, these methods are not necessarily regularising (and additionally might have trouble scaling up under the large amounts of optimisation problems to be solved). }

\revA{Measure}-based \revA{lifting} methods on the other hand have shown great potential performance-wise \cite{lellmann2013total,vogt2020lifting,vogt2019measure} \revA{and -- as will be explored in this work -- are very suitable for building in regularisation. That is, even though these methods are initially proposed as global optimisation methods -- because some lifting-based relaxations actually yield the global minimiser \cite{pock2010global} --, the process of projecting a solution back onto the original space cannot
be expected to be a minimiser of the original functional. Moreover, by choosing a suitable relaxation, the introduced error -- or bias -- can potentially be leveraged to find a regularised solution that is still close to the global minimiser. In other words, this type of method can approximate global minimisers while regularising, which is an idea that has -- to best of our knowledge -- not been explored in a similar fashion so far.}




\subsection{Specific contributions}
\begin{figure}[h!]
    \centering
    \begin{subfigure}{0.49\textwidth}
    \includegraphics[width=\linewidth]{"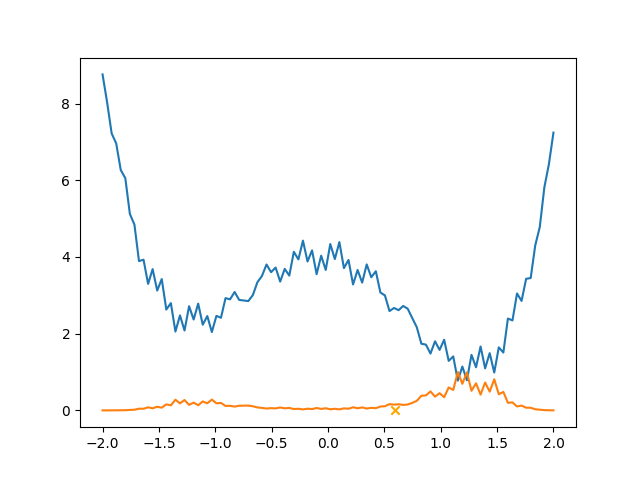"}
    \caption{ }
    \label{fig:toy-relion-weights}
    \end{subfigure}
    \hfill
    \begin{subfigure}{0.49\textwidth}
    \centering
    \includegraphics[width=\linewidth]{"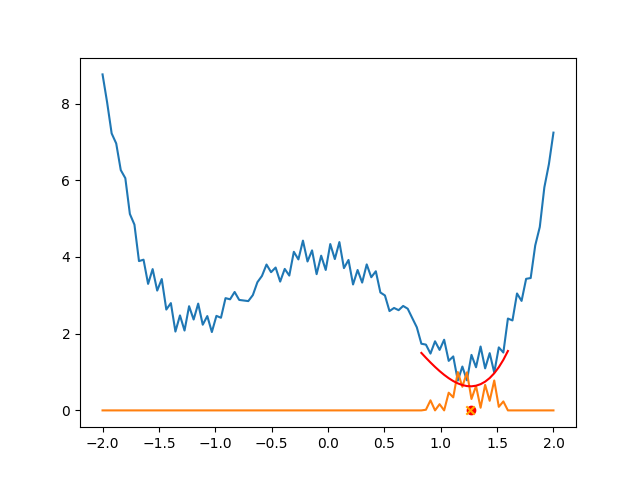"}
    \caption{ }
    \label{fig:toy-ESL-weights}
    \end{subfigure}
    \caption{\revA{The evaluation of the means of the measure generated by the weights as in \cref{eq:relion-like-weights} (a) and of the measure generated by the proposed method (b). Because of outliers on the left-most part of the measure in (a), the mean does not yield a suitable solution. The proposed method resolves this locality issue and the measure mean approximates the minimiser of the smooth lower envelope of the function (red) in this toy example, which is distinct from the global minimiser.}}
\end{figure}

\revA{Using \cref{eq:intro-barycentre-idea} to directly compute the barycentre of the measure as in RELION will run into severe issues. In particular, for $\rotsSampling \subset \mathrm{SO}(3)$ a discretisation of $\mathrm{SO}(3)$, RELION picks $\measure_\imgInd^{\iterInd}:= \sum_{\rot\in\rotsSampling} \rotsDensCoef_{\imgInd \rot}^{\iterInd} \delta_\rot$, where $\rotsDensCoef_{\imgInd \rot}^{\iterInd}$ is computed along similar lines as
\begin{equation}
    \rotsDensCoef_{\imgInd \rot}^{k}= \frac{ \exp \Bigl(-\frac{\|\forward \rot. \vol^k - \img_\imgInd\|_2^2}{2\noiseLevelParam}\Bigr) }{\sum_{\rot\in\rotsSampling} \exp \Bigl(- \frac{\|\forward \rot. \vol^k - \img_\imgInd\|_2^2}{2\noiseLevelParam}\Bigr)}.
    \label{eq:relion-like-weights}
\end{equation}
Now, the first issue can be visualized using the toy example from \cref{fig:1D-examples}. Even though the measure is mostly concentrated around the right-most part of the domain, there will always be non-zero coefficients in other parts -- here the left-most part -- of the domain. The barycentre reduces to the weighted mean in the scalar case, and the mean is not robust to outliers as becomes clear in \cref{fig:toy-relion-weights}. Furthermore, the obtained mean does not at all have a low loss and is not a suitable solution to the problem. The second issue arises from the geometry of the underlying Riemannian manifold $\manifold$. That is, the barycentre problem for a non-local measure on non-Hadamard manifolds, e.g., $\mathrm{SO}(3)$, yields a non-convex optimisation problem \cite[Figure~1]{vogt2020lifting}, which might be as difficult to solve as the original problem. 

Having a local measure with zero weights outside some neighbourhood can have a regularising effect even though it does not suffer from these issues. In particular, for the toy example \cref{fig:toy-ESL-weights} shows that the regularised solution using the approach proposed in this work is close to the minimiser of the lower envelope of the toy function, which can be a more suitable solution. In other words, adding bias arguably pays off for this toy example. Finally, this brings us to the contributions of this work.}

The overall contribution of this paper is in the theory and algorithm development \revA{of a measure-based lifting scheme for regularising and} solving minimisation problems of the form \cref{eq:model-problem} where the Riemanian manifold $\manifold$ is not necessarily Hadamard. Particular emphasis is on \revA{guaranteeing locality of the measure and} deriving error estimates \revA{due to the introduced bias (regularised} global optimisation over Riemannian manifolds).
\revA{As discussed above, s}uch optimisation problems arise when regularising certain challenging large-scale inverse problems with high noise, like the one in Cryo-EM.
This brings us to the second part of our contribution, namely to ensure the optimisation solver is noise-robust while maintaining sufficient consistency against measured data (noise-robust and data-consistent rotation estimation in Cryo-EM).

\paragraph{regularised global optimisation over Riemannian manifolds}
For \revA{regularised} global optimisation over manifolds we propose ellipsoidal support lifting (ESL), a new measure-based lifting scheme which allows for an error analysis and gives ellipsoidal locality of the support of the measure around the global minimiser that can be used to construct a practical algorithm. Theoretical results are tested in numerical experiments.

\paragraph{Noise-robust and data-consistent rotation estimation in Cryo-EM}
The joint rotation estimation and recovery of the 3D map typically involves the trade-off between noise-robustness and data-consistency. Based on the theoretical framework for regularised optimisation over manifolds that is developed, we propose a new joint recovery scheme with a lifting-based rotation estimation part that arguably takes into account both noise-robustness and data-consistency, while still being numerically feasible. \revA{Performance is tested in numerical experiments, suggesting that the induced bias additionally yields better solutions, i.e., rotations closer to the ground truth, than global minimisation can provide.}

\subsection{Outline}
\Cref{sec:prelim} introduces basic notions from differential geometry, Riemannian geometry and optimal transport which then are used in measure-based lifting. In \cref{sec:lifting-based-global-optimisation} we propose a new measure-based lifting scheme \revA{for regularising the model problem \cref{eq:model-problem} and show that this method solves the locality issue of the measure and yields a well-posed regularised solution within error bounds to the global minimiser}. \Cref{sec:so3-lifting-analysis} contains proofs of the main theoretical results. Next, the theory and algorithms from \cref{sec:lifting-based-global-optimisation} are adapted and integrated into a joint reconstruction and rotation estimation method for single particle Cryo-EM data assuming homogeneous particles (\cref{sec:application-to-cryo-EM}). Numerical results are shown in \cref{sec:numerics} and \cref{sec:conclusions} concludes the paper with a summary and discussion of our proposed theory and method for Cryo-EM.

\section{Preliminaries}
\label{sec:prelim}
The aim here is to provide essential notions from Riemannian geometry and optimal transport (\cref{sec:riemannian geometry}).
We then introduce basic ideas underlying lifting methods for optimisation on manifolds (\cref{sec:basics-lifting}). 


\subsection{Riemannian geometry and optimal transport}
\label{sec:riemannian geometry}


This overview is based on \cite{boothby2003introduction,carmo1992riemannian,lee2013smooth,sakai1996riemannian} for differential and Riemannian geometry and on  \cite{villani2009optimal} for optimal transport.

\paragraph{Riemannian Geometry}
Let $\manifold$ be a smooth manifold. Its \emph{tangent space} at $\mPoint \in \manifold$ \revA{-- defined as the space of all \emph{derivations} at $\mPoint$ --} is denoted by $\tangent_\mPoint \manifold$ and for \emph{tangent vectors} we write $\mTVector_\mPoint$ or just $\mTVector$.
Next, $\partial\manifold$ denotes its boundary (assuming it exists).
Furthermore, $\manifold$ is a \emph{Riemannian manifold} if, at each point $\mPoint \in \manifold$, there is a smoothly varying \emph{metric tensor} $(\,\cdot\,, \,\cdot\,)_\mPoint \colon \tangent_\mPoint \manifold \times \tangent_\mPoint \manifold \to \mathbb{R}$. This tensor induces a \emph{(Riemannian) metric} $\distance_{\manifold} \colon \manifold\times\manifold\to\Real_+$. 

The metric tensor can also be used to construct a unique affine connection, the \emph{Levi-Civita connection} (or \emph{covariant derivative}) that is denoted by $\nabla_{(\,\cdot\,)}(\,\cdot\,)$. 
This connection is in turn the cornerstone of a myriad of manifold mappings.
One is the notion of a \emph{geodesic}, which for two points $\mPoint,\mPointB \in \manifold$ is defined as \revA{a} curve $\geodesic_{\mPoint,\mPointB} \colon [0,1] \to \manifold$ with minimal length that connects $\mPoint$ with $\mPointB$. This notion is well-defined if the manifold is \emph{geodesically connected}, i.e., any two points $\mPoint,\mPointB\in\manifold$ can be connected with a curve that is contained in $\manifold$.
Another closely related notion is the curve $\tim \mapsto \geodesic_{\mPoint,\mTVector}(\tim)$  for a geodesic starting from $\mPoint\in\manifold$ with velocity $\dot{\geodesic}_{\mPoint,\mTVector} (0) = \mTVector \in \tangent_\mPoint\manifold$. This can be used to define the \emph{exponential map} $\exp_\mPoint \colon \mathcal{G}_\mPoint \to \manifold$ as 
\[ 
\exp_\mPoint(\mTVector) := \geodesic_{\mPoint,\mTVector}(1)
\quad\text{where $\mathcal{G}_\mPoint \subset \tangent_\mPoint\manifold$ is the set on which $\geodesic_{\mPoint,\mTVector}(1)$ is defined.} 
\]
The manifold $\manifold$ is said to be \emph{complete} whenever $\mathcal{G}_p = \mathcal{T}_{\mPoint}\manifold$.
Furthermore, the \emph{logarithmic map} $\log_\mPoint \colon \exp(\mathcal{G}'_\mPoint ) \to \mathcal{G}'_\mPoint$ is defined as the inverse of $\exp_\mPoint$, so it is well-defined on  $\mathcal{G}'_{\mPoint} \subset \mathcal{G}_{\mPoint}$ where $\exp_\mPoint$ is a diffeomorphism.
Moreover, for a smooth function $\dummyFunctionA \colon \manifold\to \Real$ and a point $\mPoint\in \manifold$, the metric tensor $(\,\cdot\,, \,\cdot\,)_\mPoint \colon \tangent_\mPoint \manifold \times \tangent_\mPoint \manifold \to \mathbb{R}$ allows to construct the \emph{Riemannian gradient} at $\mPoint\in \manifold$ as the unique tangent vector $\Grad_\mPoint \dummyFunctionA \in \tangent_\mPoint\manifold$ such that 
\[
    \bigl(\,
      \Grad_\mPoint \dummyFunctionA,
      \,\mTVector \,\bigr)_\mPoint := \mTVector \dummyFunctionA
    \quad\text{holds for any $\mTVector \in \tangent_\mPoint\manifold$.}
\]
When combined with the Levi-Civita connection, we in addition can define the \emph{Riemannian Hessian}  
\[ \Hess_\mPoint \dummyFunctionA \colon \tangent_\mPoint\manifold\times \tangent_\mPoint\manifold\to\Real
\quad\text{at $\mPoint\in \manifold$,}
\]
as the unique bilinear form that satisfies 
\[
  \Hess_\mPoint \dummyFunctionA(\mTVector,\mTVectorB)
  := (\nabla_\mTVector \Grad_\mPoint \dummyFunctionA, \mTVectorB)_\mPoint
  \quad\text{for any $\mTVector,\mTVectorB\in \tangent_\mPoint\manifold$}
\]

The final two notions derived from the Riemannian metric are the eigenvalues and a determinant of a symmetric bi-linear form. These are less frequently used in general theory of Riemannian geometry, but they will be important for us.

\begin{definition}[Eigenvalue of a bi-linear form]
Let $\manifold$ be a Riemannian manifold and the mapping  $\dummyBilForm \colon \tangent_\mPoint\manifold \times \tangent_\mPoint\manifold\to \Real$ is a bi-linear form at some given $\mPoint\in \manifold$. We say that $\eigenValue \in \Real$ is an \emph{eigenvalue} of $\dummyBilForm$ if there is a non-trivial tangent vector $\mTVector \in \tangent_\mPoint\manifold$ at $\mPoint$ such that
\begin{equation}
               \dummyBilForm(\mTVector, \mTVectorB)=
               \eigenValue(\mTVector, \mTVectorB)_{\mPoint} 
               \quad\text {for all $\mTVectorB \in \tangent_\mPoint\manifold$.}
\end{equation}
$\mTVector$ is called the \emph{corresponding eigenvector} of $\dummyBilForm$ on $\tangent_\mPoint\manifold$.
\end{definition}

\begin{definition}[Determinant of a bi-linear form]
\label{def:det-bilform}
Let $\manifold$ be a Riemannian manifold and let $\dummyBilForm \colon \tangent_\mPoint\manifold \times \tangent_\mPoint\manifold\to \Real$ be a symmetric bi-linear form at a given $\mPoint\in \manifold$ with real eigenvalues $\eigenValue_1, \ldots,\eigenValue_\DimInd$, where $\DimInd:= \dim (\manifold)$. We say that the \emph{determinant} $\det(\dummyBilForm)$ of the bi-linear form $\dummyBilForm$ is the product of the eigenvalues of $\dummyBilForm$, i.e.,
\begin{equation}
    \det(\dummyBilForm) := \prod_{\sumIndA=1}^\DimInd \eigenValue_{\sumIndA}.
\end{equation}
\end{definition}



\paragraph{Optimal transport}
Lifting is a key strategy for solving non-convex minimisation over some Riemannian manifold $\manifold$ (\cref{sec:basics-lifting}).
The idea is to replace the minimisation over $\manifold$ with one over probability measures on $\manifold$. It therefore becomes necessary to work with probability measures on $\manifold$ and to have a metric over the space of such measures. We consider measures with finite $s$-th moment and use the Wasserstein distance as a metric on such measures.
\begin{definition}[Wasserstein space and distance, {\cite[Def.6.4]{villani2009optimal}}]
Let $\manifold$ be a connected
Riemannian manifold with metric $\distance_\manifold$ and $\probSpace(\manifold)$ is the space of Borel probability measures 
over $\manifold$. The \emph{Wasserstein space of order $s \geq 1$ on $\manifold$} is defined as
\begin{equation}\label{eq:WassersteinSpace}
    \probSpace_{s}(\manifold):=\biggl\{\measure \in \probSpace(\manifold) \;\mid\; \int_{\manifold} \distance_\manifold (\mPoint_{0}, \mPoint)^{s} \mathrm{d}\mu(\mPoint) <+\infty, \; \text{for all $\mPoint_0 \in \manifold$} \biggr\},
\end{equation}
and the \emph{Wasserstein-$s$ distance} is a mapping $\wasserstein_{s} \colon \probSpace_{s}(\manifold) \times \probSpace_{s}(\manifold) \to \Real_+$ defined as
\begin{equation}
    \wasserstein_{s}(\mu, \nu):=\Bigl(\,
      \inf _{\pi \in \Gamma(\mu, \nu)} \int_{\manifold\times\manifold} \distance_\manifold(\mPoint, \mPointB)^{s} \; \mathrm{d} \pi(\mPoint, \mPointB)
    \Bigr)^{1 / s}.
\end{equation}
Here, $\Gamma(\mu, \nu) \subset \probSpace(\manifold \times \manifold)$ is the set of all \emph{couplings}, i.e., probability measures on $\manifold \times \manifold$ with marginals $\mu$ and $\nu$.
\end{definition}
The Wasserstein-$s$ distance defines a (finite) distance on $\probSpace_{s}(\manifold)$ that defines the topology in which we consider convergence of sequences of measures. 

\subsection{Lifting methods for global optimisation on manifolds}
\label{sec:basics-lifting}
The minimisation in \cref{eq:model-problem} is typically a non-convex optimisation problem, so solving it with standard techniques from convex optimisation run into the risk of getting stuck in local minima. 
One approach to address such non-convex manifold optimisation problems is through \emph{lifting}, which provides -- for low-dimensional manifolds -- a computationally feasible convex relaxation.
In case of measure-based lifting, the method consists of roughly four main steps: convexification, projection, discretisation and relaxation.
The first three steps are fairly standard when lifting is phrased within the space of probability measures, but the choice of relaxation may differ in-between approaches (a topic we deal with in \cref{sec:general-lifting} through the lense of regularisation).



\paragraph{Convexification}
For convexification, one \emph{lifts} the manifold $\manifold$ into $\probSpace_2(\manifold)$, which is the Wasserstein space of order 2 that is defined in \cref{eq:WassersteinSpace}. 
More precisely, the idea is replace \cref{eq:model-problem} with the convex optimisation problem
\begin{equation}
    \inf_{\measure \in\probSpace_2(\manifold)} \Bigl\{ 
    \int_{\manifold} \dummyFunctionA(\mPoint) \mathrm{d}\measure(\mPoint)  \Bigr\}.
    \label{eq:lifted-function}
\end{equation}

\paragraph{Projection}
If one can guarantee existence of solutions to \cref{eq:lifted-function}, then such a solution 
$\optim{\measure}\in \probSpace_2(\manifold)$ can be projected back onto the manifold $\manifold$ through a (possibly set-valued) \emph{(Riemannian) barycentre mapping} 
\begin{equation}
    P \colon \probSpace_2(\manifold)\rightrightarrows\manifold
    \quad\text{where}\quad
    \mPoints(\measure) 
    := \argmin_{\mPointB\in \manifold}\Bigl\{ \int_{\manifold} \distance_{\manifold}(\mPoint, \mPointB)^2 \mathrm{d}\measure(\mPoint) \Bigr\}
    \text{ for $\measure \in \probSpace_2(\manifold)$.}
    \label{eq:lifted-projection}
\end{equation}
The above is by \cite[Proposition~1]{le2017existence} well-defined for all measures $\measure\in \probSpace_2(\manifold)$ whenever $\manifold$ is geodesically connected. For connected manifolds existence holds only locally, i.e., for measures with support on a geodesically connected subset.

Furthermore, if $\optim{\measure}$ is a unique solution to \cref{eq:lifted-function} that is given as a Dirac measure, i.e.,  $\optim{\measure}=\dirac_q$ for some $q \in \manifold$, then $q \in \manifold$ also solves \cref{eq:model-problem} and it is recovered by the barycenter mapping. However, $\optim{\measure}$ may differ from a Dirac measure, e.g., it may even be diffuse.

\paragraph{Discretisation}
A numerical implementation for a method that solves \cref{eq:lifted-function} needs to include a discretisation.
A common approach is to restrict the control variable in the optimisation to an appropriate finite dimensional subspace of $\probSpace_2(\manifold)$. This procedure will give rise to a \emph{discretisation error}.

For the (geodesically) connected Riemannian manifold $\manifold$, there is a natural discretisation of $\probSpace_2(\manifold)$ obtained by considering the sub-space of measures that are linear combinations of Dirac distributions located at finitely many sample points in $\manifold$.
To formalise this, we note first that \emph{any} measure $\measure\in \probSpace_2(\manifold)$ can be approximated by a sum of Dirac measures \cite[Lemma 6.2]{kim2017wasserstein}.
Hence, it is natural to consider probability measures on $\manifold$ that are given as linear combinations of Dirac measures, i.e., measures of the form
\begin{equation}
    \sum_{\mPoint\in \rotsSampling} \rotsDensCoef_{\mPoint} \dirac_{\mPoint}
    \quad\text{for some $\rotsDensCoef_{\mPoint} \in \Real$ with fixed finite set $\rotsSampling\subset \manifold$ (sampling points).}
\label{eq:DicsretisedProbMeas}
\end{equation} 
Note that the linear combination in \cref{eq:DicsretisedProbMeas} defines a probability measure on $\manifold$ whenever the coefficient array $\rotsDensCoef=( \rotsDensCoef_{\mPoint} )_{\mPoint\in \rotsSampling} \in \Real^{\numRots}$ in \cref{eq:DicsretisedProbMeas} is in the unit $\numRots$-dimensional simplex:
\begin{equation} 
\rotsDensCoef\in \simplex^{\numRots}:= 
\bigl\{
  \ePoint\in\Real^{\numRots} \mid 
  \ePoint\geq 0,\mathbf{1} \cdot \ePoint =1 
\bigr\}.
\label{eq:UnitSimplex}
\end{equation}
To summarise, the set of measures on $\manifold$ of the form in \cref{eq:DicsretisedProbMeas} for some fixed $\rotsSampling\subset \manifold$ is a discretisation of $\probSpace_2(\manifold)$, i.e., a finite dimensional subspace of $\probSpace_2(\manifold)$.
Note also that if $\measure \in \probSpace_2(\manifold)$ is a measure of the form \cref{eq:DicsretisedProbMeas}, then integrating $\dummyFunctionA \colon \manifold \to \Real$ against $\measure$ reduces to taking the inner product of  $\dummyFunctionA_\rotsSampling \in \Real^{\numRots}$ (the array obtained by sampling $\dummyFunctionA$ on $\rotsSampling$) and $\rotsDensCoef \in \simplex^{\numRots} \subset \Real^{\numRots}$:
\begin{equation}
    \int_{\manifold} \dummyFunctionA(\mPoint) \mathrm{d}\measure(\mPoint)
= \revA{\sum_{\mPoint\in \rotsSampling} \rotsDensCoef_\mPoint \dummyFunctionA(\mPoint)}
=
\dummyFunctionA_\rotsSampling \!\cdot \rotsDensCoef 
\quad\text{where $\dummyFunctionA_\rotsSampling:= (\dummyFunctionA(\mPoint))_{\mPoint\in\rotsSampling}\in \Real^{\numRots}$.}
\label{eq:filled-in-measure-disc}
\end{equation}

\paragraph{Relaxation}
Discretising the lifted problem \cref{eq:lifted-function} by \cref{eq:filled-in-measure-disc} amounts to minimising 
\begin{equation}
\inf_{\rotsDensCoef \in\Real^{\numRots}} 
\Bigl\{ 
  \dummyFunctionA_\rotsSampling \!\cdot \rotsDensCoef 
  + \indicator_{\simplex^{\numRots}}(\rotsDensCoef)  
\Bigr\},
\label{eq:lifted-function-discr-no-reg}
\end{equation}
where $\indicator_{\simplex^{\numRots}}\colon \Real^{\numRots} \to \{0,\infty\}$ is defined as  
\begin{equation}
\indicator_{\simplex^{\numRots}}(\rotsDensCoef):= 
\begin{cases}
  0, & \text{if $\rotsDensCoef \in \simplex^{\numRots}$},
  \\ 
  \infty, & \text{otherwise.}
\end{cases}
\label{eq:UnitSimplexIndicator}
\end{equation}
The solution to optimisation problem \cref{eq:lifted-function-discr-no-reg} will collapse into the specific sampling point $\mPoint\in\rotsSampling$ where $\dummyFunctionA$ has its  lowest value \revA{or will be a convex combination of \enquote{one-hot} vectors if there are multiple points with the same lowest value}. 
\revA{In the former case} the corresponding barycentre will then equal that point, so without further adaptations to the optimisation problem, one will need a very fine mesh -- and unfeasibly many sample points -- to get an accurate approximation of the global minimiser to $\dummyFunctionA$.
A typical approach to overcome this, is to make adaptations to the discretized energy \cref{eq:lifted-function-discr-no-reg} such that $\rotsDensCoef$ has multiple non-zero entries. The barycentre mapping then allows one to find values between the sampling points so that one is no longer directly restricted by the mesh density.

Approaches for relaxation vary depending on the underlying optimisation problem \cite{lellmann2013total,vogt2020lifting,vogt2019measure}.
Typically, the optimisation problem after relaxation is computationally feasible. However, it is difficult to perform a mathematical analysis of the errors due to convexification, projection, discretisation and relaxation. This still remains an open problem in the general setting.

\section{Ellipsoidal support lifting for Riemannian manifolds}
\label{sec:lifting-based-global-optimisation}
\revA{Unlike the classical goal of lifting, i.e., approximating the global minimiser of a function, we propose a new measure-based lifting scheme for regularising the model problem \cref{eq:model-problem} and show that this method solves the locality issue of the measure and yields a well-posed regularised solution within error bounds. 

We begin this section by briefly considering the type of functions a measure-based lifting method should be able to handle in a realistic Cryo-EM rotation estimation setting (\cref{sec:lifting-for-cryo-rotations}). Subsequently, w}e propose a new relaxation method for measure-based lifting on $\manifold$ \revA{and discuss several observations regarding computational feasibility and the method's potential for sparsity of the measure without collapsing into a single point (\cref{sec:general-lifting}). These observations will be made more concrete in our main theoretical results} (\cref{sec:main-results}). Additionally, it yields a computationally feasible algorithm (\cref{sec:lift-algo}). The main theoretical results stated in \cref{sec:main-results} are proved in \cref{sec:so3-lifting-analysis}. 
Throughout this section, $\ballVol_\DimInd$ denotes the volume of a $\DimInd$-dimensional unit ball.

\subsection{Measure-based lifting for Cryo-EM rotation estimation}
\label{sec:lifting-for-cryo-rotations}
\revA{The design choices of the proposed method and the provided theory are highly motivated by Cryo-EM rotation loss functions. In particular, Cryo-EM rotation losses are \emph{smooth}, but have very large Lipschitz gradients near steep drops into \emph{several narrow valleys of very low loss}. Additionally, the \emph{global minimiser -- living in one of these valleys -- is hardly ever the ground truth rotation}. However, if the noise level is not too high, empirically \emph{the ground truth rotation tends to be in the same valley of very low loss}. The proposed method and the provided analysis caters to these types of functions. In particular, in the absence of a good model for the ground truth rotation, the proposed model uses the key assumption that a good solution should have a low loss in a small neighbourhood. 

Then, a good regularisation method allows to pick the size of this neighbourhood and if the regularisation goes to zero, we would recover the global minimiser. Because lifting methods cannot be expected to give the global minimiser for a finite sampling set, the regularisation parameter should be a function of the sampling set size. Now, as we increase the amount of points in the sampling set, the error we make with respect to the minimiser should go down, but it should go down in a regular fashion, i.e., the resulting measure should not collapse into some point when increasing the size. As we will see in the following, we can invoke smoothness assumptions at the global minimiser to alleviate this issue and with that construct a proper regularisation method that is suitable for Cryo-EM rotation loss functions.
}

\subsection{A new relaxation strategy}
\label{sec:general-lifting}
Adding the Tikhonov-type regulariser $\rotsDensCoef \mapsto \regParam \|\rotsDensCoef\|_2^2/2\numRots^\exponentA$ as a relaxation on the coefficients $\rotsDensCoef$ in \cref{eq:lifted-function-discr-no-reg} leads to the following variational problem:
\begin{equation}
\inf_{\rotsDensCoef \in\Real^{\numRots}} 
\Bigl\{ 
  \dummyFunctionA_\rotsSampling \!\cdot \rotsDensCoef 
  + \indicator_{\simplex^{\numRots}}(\rotsDensCoef)  
  + \frac{\regParam}{2\numRots^\exponentA}\|\rotsDensCoef\|_2^2 
\Bigr\}
\quad\text{for regularisation parameters $\regParam>0$ and $\exponentA \in \Real$.}
\label{eq:lifted-function-discr}
\end{equation}
The projection map in \cref{eq:lifted-projection} applied to the measure $\optim{\measure}:= \sum_{\mPoint\in \rotsSampling}\optim{\rotsDensCoef}_{\mPoint}\dirac_{\mPoint}$ that solves \cref{eq:lifted-function-discr} yields the barycentre, which is the point in $\manifold$ given as
\begin{equation}
P\Bigl(
  \sum_{\mPoint\in \rotsSampling}
    \optim{\rotsDensCoef}_{\mPoint}\dirac_{\mPoint}
  \Bigr)
  = \argmin_{\mPointB\in \manifold}\Bigl\{
      \sum_{\mPoint\in \rotsSampling} \optim{\rotsDensCoef}_{\mPoint} \, \distance_{\manifold}(\mPoint,\mPointB)^2  
  \Bigr\}.
\label{eq:lifted-projection-discr}
\end{equation}
\revA{We make three observations.}

\paragraph{Computational feasibility}
Choosing a Tikhonov type regulariser re\revA{sults in} computational feasibility. In fact, the (finite-dimensional) minimisation in \cref{eq:lifted-function-discr} is not only strongly convex, its unique minimiser can be expressed in closed form. To see this, note first that the unique minimiser $\optim{\rotsDensCoef} \in \simplex^{\numRots}$ to \cref{eq:lifted-function-discr} can be expressed more conveniently by rewriting the objective functional
\begin{equation}
\begin{split}
\optim{\rotsDensCoef} 
  &= \argmin_{\rotsDensCoef \in\Real^{\numRots}}\Bigl\{ 
       \dummyFunctionA_\rotsSampling \cdot \!\rotsDensCoef 
       + \indicator_{\simplex^{\numRots}}(\rotsDensCoef)  
       + \frac{\regParam}{2\numRots^\exponentA}\|\rotsDensCoef\|_2^2 
     \Bigr\} 
\\
  &= \argmin_{\rotsDensCoef \in\Real^{\numRots}}\biggl\{
       \frac{\regParam}{2\numRots^\exponentA}\Bigl\|
           \rotsDensCoef + \frac{\numRots^\exponentA}{\regParam}
           \dummyFunctionA_\rotsSampling 
         \Bigr\|_2^2 
         + \indicator_{\simplex^{\numRots}}(\rotsDensCoef) 
     \biggr\} 
     = \Pi_{\simplex^{\numRots}}\left(-\frac{\numRots^\exponentA}{\regParam}\dummyFunctionA_\rotsSampling\right),
\end{split}
\label{eq:lifted-function-discr-argmin}
\end{equation}
where $\Pi_{\simplex^{\numRots}}\colon \Real^{\numRots} \to \simplex^{\numRots}$ is the projection onto the unit $\numRots$-simplex $\simplex^{\numRots}$ in \cref{eq:UnitSimplex}. \revA{Note that the last equality holds by definition of projection onto a convex set.}
As we show next, the final term in \cref{eq:lifted-function-discr-argmin} has a closed-form expression that can be computed efficiently.
\begin{lemma}[{\cite[Lemma 1\&2]{duchi2008efficient}}]
\label{thm:simplex-projection}
For any $\ePoint \in \Real^{\DimInd}$ with $\ePoint_{(1)}\geq \ldots\geq \ePoint_{(\DimInd)}$ are the components of $\ePoint$ sorted in descending order, there is a unique \emph{cut-off integer} $\dummyIntegerB=\dummyIntegerB(\ePoint)$, so that
\begin{equation}
    \ePoint_{(\dummyIntegerB)} \geq \frac{1}{\dummyIntegerB} \Bigl(\sum_{j=1}^{\dummyIntegerB} \ePoint_{(j)}\Bigr) - \frac{1}{\dummyIntegerB} \geq \ePoint_{(\dummyIntegerB+1)} 
    \quad\text{or}\quad
    \ePoint_{(\DimInd)}\geq \frac{1}{\dummyIntegerB} \Bigl(\sum_{j=1}^{\dummyIntegerB} \ePoint_{(j)}\Bigr) -  \frac{1}{\dummyIntegerB}. 
    \label{eq:cutoff-inequalities}
\end{equation}
holds.

Then, the 
projection $\Pi_{\simplex^{\!\DimInd}} \colon \Real^{\DimInd} \to \simplex^{\!\DimInd}$ of 
$\Real^{\DimInd}$ onto the $\DimInd$-simplex $\simplex^{\!\DimInd}$ can be computed as
    \begin{equation}
        \Pi_{\simplex^{\!\DimInd}}(\ePoint)=\biggl(\ePoint-\frac{1}{\dummyIntegerB} \Bigl(\sum_{j=1}^{\dummyIntegerB} \ePoint_{(j)}\Bigr) + \frac{1}{\dummyIntegerB}\biggr)_{+}.
        \label{eq:simplex-projection-form}
    \end{equation}
In particular, note that $\dummyIntegerB = \bigl\|\Pi_{\simplex^{\!\DimInd}}(\ePoint) \bigr\|_0$ holds.
\end{lemma}

\begin{remark}
\label{rem:monotonicity-of-projection-correction}
In the analysis of the relaxation method -- in particular in the proof of \cref{thm:new-num-non-zero-coeffs} -- it will be important to consider the uniqueness of the cut-off integer $\dummyIntegerB$ closer. That is, observe that $n\mapsto \frac{1}{n} \Bigl(\sum_{j=1}^{n} \ePoint_{(j)}\Bigr) -  \frac{1}{n}$ is monotonically increasing. Then, the uniqueness of the cut-off integer from \cref{thm:simplex-projection} tells us that $\|(\ePoint-\frac{1}{\dummyIntegerB'} \Bigl(\sum_{j=1}^{\dummyIntegerB'} \ePoint_{(j)}\Bigr) + \frac{1}{\dummyIntegerB'})_{+}\|_0 < \dummyIntegerB$ is equivalent to having $\dummyIntegerB'>\dummyIntegerB$, i.e., there will be \emph{too few} non-zero coefficients iff $\dummyIntegerB'$ \enquote{were to be used as cut-off integer}. Conversely, choosing an $\dummyIntegerB'<\dummyIntegerB$ is equivalent to having \emph{too many} non-zero coefficients. 
\end{remark}


\paragraph{Sparsity and locality of the measure}
The second \revA{observation} relates to having a local and sparse measure concentrated around the global minimiser of \cref{eq:model-problem}.
Note that \revA{locally and up to higher-order terms} the coefficients generated by the Tikhonov-type lifting scheme in \cref{eq:lifted-function-discr} \revA{can be understood through} the Riemannian Hessian of the objective at the global minimiser $\optim{\mPoint}$.
This can be seen through Taylor expansion 
\begin{equation}
\begin{split}
    \optim{\rotsDensCoef}_\mPoint = \Pi_{\simplex^{\numRots}}\Bigl(- \frac{\numRots^{\exponentA}}{\regParam}\dummyFunctionA_{\rotsSampling} \Bigr)_\mPoint 
    &= \biggl(- \frac{\numRots^{\exponentA}}{\regParam}\dummyFunctionA(\mPoint)-\frac{1}{\dummyIntegerB} \Bigl(\sum_{j=1}^{\dummyIntegerB} - \frac{\numRots^{\exponentA}}{\regParam}\dummyFunctionA(\mPoint_{(j)})-1\Bigr)\biggr)_{+} 
    \\
    &= \biggl( \frac{1}{\dummyIntegerB} +\frac{\numRots^{\exponentA}}{\regParam} \Bigl( \frac{1}{\dummyIntegerB} \Bigl(\sum_{j=1}^{\dummyIntegerB}  \dummyFunctionA(\mPoint_{(j)})\Bigr) - \dummyFunctionA(\mPoint) \Bigr) \biggr)_{+} 
    \\
    &\approx \biggl( C^{\exponentA,\regParam}_\rotsSampling - \frac{\numRots^{\exponentA}}{2\regParam\dummyIntegerB} \Hess_{\optim{\mPoint}} \dummyFunctionA\bigl(\log_{\optim{\mPoint}}(\mPoint), \log_{\optim{\mPoint}}(\mPoint)\bigr) \biggr)_{+}
\end{split}   
\label{eq:riem-hess-shape-coeffs}
\end{equation}
for $\mPoint \in \rotsSampling$ and where 
\[ 
C^{\exponentA,\regParam}_\rotsSampling:= \frac{1}{\dummyIntegerB} +\frac{\numRots^{\exponentA}}{2\regParam} \frac{1}{\dummyIntegerB} \Bigl(\sum_{j=1}^{\dummyIntegerB}  \Hess_{\optim{\mPoint}} \dummyFunctionA \bigl(\log_{\optim{\mPoint}}(\mPoint_{(j)}), \log_{\optim{\mPoint}}(\mPoint_{(j)})\bigr)\Bigr).
\]
From \cref{eq:riem-hess-shape-coeffs} we also see that the coefficients of $\rotsDensCoef$ become zero outside some level set of the $\mPoint\mapsto \Hess_{\optim{\mPoint}} \dummyFunctionA\bigl(\log_{\optim{\mPoint}}(\mPoint), \log_{\optim{\mPoint}}(\mPoint)\bigr)$ mapping.

For the approximation to become \revA{useful for regularisation purposes} -- while ensuring a sparse and local measure as suggested above -- we will show in the following that we need to set
\begin{equation}
 \label{eq:thm-new-num-non-zero-coeffs-gamma}
     \regParam := \frac{1}{2(\DimInd+2)} \dummyIntegerB_0^{\frac{\DimInd+2}{\DimInd}} \sqrt[\DimInd]{\det (\Hess_{\optim{\mPoint}} \dummyFunctionA)  \Bigl(\frac{\operatorname{vol}(\manifold)}{\ballVol_\DimInd}\Bigr)^2 }.
 \end{equation}
where $\DimInd:= \dim (\manifold)$ and $\dummyIntegerB_0$ is roughly the amount of non-zero coefficients one chooses to have, and
the parameter $\eta \in  \bigl(1/(\DimInd+1),2/\DimInd\bigr)$.

\paragraph{Resolving collapsing}
The third \revA{observation} relates to avoiding collapsing into a single point. Consider the optimisation problem \cref{eq:lifted-function-discr} for setting the coefficients $\rotsDensCoef$.
As sparse (and local) solutions have been suggested above to be key for the lifting approach, one might wonder why a Tikhonnov-type  $\ell^2$-regulariser has been chosen rather than $\ell^1$. 
The reason is that every $\ell^1$-regulariser will collapse into a single point in practice. Indeed
\begin{equation}
\inf_{\rotsDensCoef \in\Real^{\numRots}} 
\Bigl\{ 
  \dummyFunctionA_\rotsSampling \!\cdot \rotsDensCoef 
  + \indicator_{\simplex^{\numRots}}(\rotsDensCoef)  + \|\rotsDensCoef\|_1 
\Bigr\} = \inf_{\rotsDensCoef \in\simplex^{\numRots}} 
\Bigl\{ 
  \dummyFunctionA_\rotsSampling \!\cdot \rotsDensCoef 
  + \mathbf{1} \cdot \rotsDensCoef 
\Bigr\}
= \inf_{\rotsDensCoef \in\simplex^{\numRots}} 
\Bigl\{ 
  (\dummyFunctionA_\rotsSampling +\mathbf{1})  \!\cdot \rotsDensCoef 
\Bigr\},
\end{equation}
which will collapse into the specific sampling point $\mPoint\in\rotsSampling$ where $\dummyFunctionA$ has its lowest value. Hence, the choice for $\ell^1$ is not suitable,
whereas we will argue that a properly scaled Tikhonov-type  regulariser can overcome the collapsing issue while preserving sparsity (\cref{thm:new-num-non-zero-coeffs}).

\begin{remark}
\label{rem:lifting-tikh-improved}
Finally, one can wonder whether any properly scaled strongly convex regulariser could yield similar or better results to using a Tikhonov-type regulariser as a relaxation on the coefficients $\rotsDensCoef$ -- that is, on top of computing the coefficients in a computationally feasible manner and overcoming the collapsing issue while preserving sparsity and locality \revA{-- and furthermore, whether it could yield a better regularised solution to the original problem}. 
\end{remark}
At this point \cref{rem:lifting-tikh-improved} is an open question as the lack of closed form solutions for the resulting optimisation problem \revA{for the coefficients $\rotsDensCoef$} makes analysis difficult. Indeed, even for simple $\ell^s$-regularisation for $s\neq 1,2$ these solutions are not available to the best of our knowledge. \revA{More importantly for the case of Cryo-EM rotation estimation, lack of closed form solutions results in computational infeasibility, because the optimisation over the weights will have to be performed for every image and throughout several refinement steps, i.e., in practice this problem will have to be solved millions of times.}



\subsection{Main theoretical results}
\label{sec:main-results}
As suggested above with \emph{locality and sparsity}, we will only be interested in the \emph{local} behaviour of functions around its minimiser. In particular, we will be interested in certain ellipsoids generated by the Riemannian Hessian of the function $\dummyFunctionA$ at its minimiser. More generally, one can define ellipsoids generated by symmetric positive definite bi-linear forms.

\begin{definition} 
\label{def:ellipsoid}
Let $\manifold$ be a Riemannian manifold with dimension  $\DimInd:= \Dim(\manifold)$, $\mPoint\in \manifold$, and let $\dummyBilForm \colon \tangent_\mPoint\manifold\times \tangent_\mPoint\manifold\to \Real$ be a symmetric positive definite bi-linear form. 
We define the \emph{open ellipsoid $\ellipsoid_\radiusB^\dummyBilForm(\mPoint) \subset \manifold$ with radius $\radiusB>0$ around $\mPoint \in \manifold$ generated by $\dummyBilForm$} as
\begin{equation}\label{eq:Ellipsoid}
  \ellipsoid_\radiusB^\dummyBilForm(\mPoint) 
  := \bigl\{ \mPointB \in \manifold \mid  \dummyBilForm(\log_{\mPoint} \mPointB, \log_{\mPoint} \mPointB) 
  <
  \sqrt[\DimInd]{\operatorname{det}(\dummyBilForm)}
  \, \radiusB^2
  \bigr\}.
\end{equation}
\end{definition}
For the following results to hold it is key that the sampling set $\rotsSampling$ is \enquote{\emph{distributed uniformly}} over the manifold. What it means for sampling sets $\rotsSampling$ to cover the manifold well will be formalized by the \emph{sampling set's ability to approximate local integration}. For completeness, the definition used throughout this work is provided below, but for the reader's convenience it will be \revA{deconstructed} in \cref{sec:local-low-discr} as it is rather technical.



\begin{definition}[local low-discrepancy sequence]
\label{def:eta-local-low-discrepancy-sequence}
Let $\manifold$ be a finite-volume Riemannian manifold with dimension $\DimInd:= \dim (\manifold)$ and $S(\manifold):= \bigl\{
    \rotsSampling \subset \manifold 
    \mid 
    |\rotsSampling|<\infty 
\bigr\}$ denotes the collection of all finite sampling sets in $\manifold$.
A sequence $(\rotsSampling_\sumTotA)_{\sumTotA=1}^\infty \subset S(\manifold)$ is a \emph{local low-discrepancy sequence on $\manifold$} for \revA{parameters $0 < \exponentA < 2/\DimInd$ and $\dummyIntegerB>0$}
if $|\rotsSampling_\sumTotA|<|\rotsSampling_{\sumTotA+1}|$ for all $\sumTotA\in \Natural$ and 
for any symmetric positive-definite bi-linear form $\dummyBilForm\colon \tangent_\mPoint\manifold\times\tangent_\mPoint\manifold\to\Real$ at any point $\mPoint\in\manifold$ and for
\[ \radiusA_\sumTotA:= \sqrt[\DimInd]{\dummyIntegerB\dfrac{\operatorname{vol}(\manifold)}{\ballVol_\DimInd}} |\rotsSampling_\sumTotA|^{-\frac{1+\exponentA}{\DimInd+2}}
\]
the local integration approximation has the following asymptotic behaviour as $|\rotsSampling_\sumTotA|\to\infty$:
\begin{enumerate}
\item Local integration approximation for constant functions: 
\begin{equation}
\displayindent0pt
\displaywidth\textwidth
\biggl| \frac{|\ellipsoid_{\radiusA_\sumTotA}^\dummyBilForm(\mPoint) \cap \rotsSampling_\sumTotA|}{|\rotsSampling_\sumTotA|} - \frac{\operatorname{vol}\bigl(\ellipsoid_{\radiusA_\sumTotA}^\dummyBilForm(\mPoint)\bigr)}{\operatorname{vol}(\manifold)} \biggr|
    \in o\Bigl(|\rotsSampling_\sumTotA|^{-\frac{\DimInd(1+\exponentA)}{\DimInd+2}}\Bigr) 
\label{eq:eta-local-low-discrepancy-sequence-1}
\end{equation}
\item Local integration approximation for bi-linear forms: 
\begin{equation}
\displayindent0pt
\displaywidth\textwidth
 \biggl|\sum_{\mPointB\in \ellipsoid_{\radiusA_\sumTotA}^\dummyBilForm(\mPoint) \cap \rotsSampling_\sumTotA}\!\!\!\!
    \frac{\dummyBilForm\bigl(\log_\mPoint(\mPointB),\log_\mPoint(\mPointB)\bigr)}{|\rotsSampling_\sumTotA|} - \frac{\int_{\ellipsoid_{\radiusA_\sumTotA}^\dummyBilForm(\mPoint)}\dummyBilForm\bigl(\log_\mPoint(\mPointB),\log_\mPoint(\mPointB)\bigr)\; \mathrm{d}\mPointB}{\operatorname{vol}(\manifold)} \biggr|
    \in o\Bigl(|\rotsSampling_\sumTotA|^{-(1+\exponentA)}\Bigr).
    \label{eq:eta-local-low-discrepancy-sequence-2}
\end{equation}
\end{enumerate}
\end{definition}

With the above notion of local low-discrepancy sequence of sampling sets on $\manifold$, we can formulate the main results, \cref{thm:new-num-non-zero-coeffs,thm:distance-bound-post-lifting}. The proofs are provided in \cref{sec:so3-lifting-analysis}.

\Cref{thm:new-num-non-zero-coeffs} justifies the choice of Tikhonov-type regularisation term in \cref{eq:lifted-function-discr} and can be summarized by two important insights that hold when the corresponding regularisation parameter $\regParam$ is appropriately chosen.  
The first insight is that non-zero coefficients associated with the probability measure will be concentrated on an ellipsoid around the global minimiser, the second is that there will be very few non-zero coefficients, i.e., $\optim{\rotsDensCoef}$ is sparse:
\begin{theorem}[Sparsity \& locality of the optimal measure]
\label{thm:new-num-non-zero-coeffs}
Consider the minimisation in \cref{eq:lifted-function-discr} where $\dummyFunctionA \colon \manifold\to \Real$ is a smooth function defined on a finite-volume Riemannian manifold $\manifold$ with dimension $\DimInd:= \dim (\manifold)$. 
Assume also that $\dummyFunctionA$ has a unique global minimiser $\optim{\mPoint}\in \manifold$ and its Riemannian Hessian $\Hess_{\optim{\mPoint}} \dummyFunctionA$ is positive definite at $\optim{\mPoint}$.
Furthermore, \revA{let $\epsilon_1, \epsilon_2, \epsilon_3, \epsilon_4 \in (\frac{\DimInd}{2\DimInd + 2},1)$}, let $1/(\DimInd+1) < \exponentA < 2/\DimInd$ and set $\regParam$ as in \cref{eq:thm-new-num-non-zero-coeffs-gamma} for some $\dummyIntegerB_0\geq 3^{-\frac{\DimInd}{\DimInd+2}}$. Finally, assume $(\rotsSampling_\sumTotA)_{\sumTotA=1}^\infty$ is a local low-discrepancy sequence \revA{for all pairs of parameters $(\eta, \dummyIntegerB)$ where $\dummyIntegerB \in \{\dummyIntegerB_0,\bigl( \frac{\DimInd + 1 + \epsilon_1}{\DimInd + 2}\bigr)^{\frac{\DimInd}{2}}\dummyIntegerB_0, \bigl(\frac{\DimInd + 3 - \epsilon_2}{\DimInd + 2}\bigr)^{\frac{\DimInd}{2}}\dummyIntegerB_0, 3^{-\frac{\DimInd}{\DimInd+2}}\dummyIntegerB_0,\bigl(\frac{\DimInd + 1 + \epsilon_3}{\DimInd + 2}\bigr)^{\frac{\DimInd}{2}}3^{-\frac{\DimInd}{\DimInd+2}}\dummyIntegerB_0, \bigl(\frac{\DimInd + 3 - \epsilon_4}{\DimInd + 2}\bigr)^{\frac{\DimInd}{2}}3^{-\frac{\DimInd}{\DimInd+2}}\dummyIntegerB_0 \}$}.


Then, there exists $\sumTotA'\in \Natural$ such that for all $\sumTotA\geq \sumTotA'$ the points $\mPoint\in \rotsSampling_\sumTotA$ with non-zero coefficients must lie in an ellipsoid  $\ellipsoid_{\radiusA_\sumTotA}^{\dummyBilForm}\bigl(\optim{\mPoint}\bigr)$ of the type in \cref{eq:Ellipsoid} where $\dummyBilForm\colon \tangent_\mPoint\manifold\times\tangent_\mPoint\manifold\to\Real$
is $\dummyBilForm:=\Hess_{\optim{\mPoint}} \dummyFunctionA$ and $\radiusA_\sumTotA > 0$ given as
\begin{equation}
\label{eq:EllipsRadii}
\radiusA_\sumTotA := \sqrt[\DimInd]{\dummyIntegerB_0 \frac{\operatorname{vol}(\manifold)}{\ballVol_\DimInd}} |\rotsSampling_\sumTotA|^{-\frac{1+\exponentA}{\DimInd+2}}.    
\end{equation} 
Furthermore, if $\optim{\rotsDensCoef}_\sumTotA \in \simplex^{|\rotsSampling_\sumTotA|}$ is the unique minimiser to \cref{eq:lifted-function-discr} for sampling set $\rotsSampling_\sumTotA$, then 
\begin{equation}
    \label{eq:thm-J-bounds}
        3^{-\frac{\DimInd}{\DimInd+2}} \dummyIntegerB_0|\rotsSampling_\sumTotA|^{\frac{2 - \DimInd\exponentA}{\DimInd+2}} 
    \leq 
    \| \optim{\rotsDensCoef}_\sumTotA \|_0
    \leq \dummyIntegerB_0|\rotsSampling_\sumTotA|^{\frac{2 - \DimInd\exponentA}{\DimInd+2}}.
\end{equation}
\end{theorem}
The next result, \cref{thm:distance-bound-post-lifting}, establishes existence and uniqueness of the barycentre and provides a bound on the error made by convexification, projection, discretisation, and relaxation. \revA{In terms of regularising the original problem, the introduced bias will vanish asymptotically as the size of the sampling set increases. However, for a finite sampling set tuning the support of the measure -- through choosing a large enough value for $\dummyIntegerB_0$ -- prevents us from estimating the global mimimiser and ensures a regularising effect through averaging over a suitable amount of points:}

\begin{theorem}[Well-posedness \& error bounds]
\label{thm:distance-bound-post-lifting}
Consider the minimisation in \cref{eq:lifted-function-discr} with the same assumptions as in \cref{thm:new-num-non-zero-coeffs}.

Then, there exists $\sumTotA'\in \Natural$ such that for all $\sumTotA\geq \sumTotA'$
the Riemannian barycentre w.r.t.\@ the probability measure given by the unique minimiser $\optim{\rotsDensCoef}_\sumTotA \in \simplex^{|\rotsSampling_\sumTotA|}$ to \cref{eq:lifted-function-discr} for sampling set $\rotsSampling_\sumTotA$  exists and is single valued, i.e., 
\begin{equation}
    \optim{\mPoint}_{\rotsSampling_\sumTotA}:= \argmin_{\mPointB\in \manifold}
     \biggl\{\sum_{\mPoint\in \rotsSampling_\sumTotA} 
    (\optim{\rotsDensCoef}_\sumTotA)_{\mPoint}
    \,
    \distance_{\manifold}(\mPoint, \mPointB)^2  
    \biggr\}
    \quad\text{is well-defined point in $\manifold$.}
    \label{eq:thm-barycentre}
\end{equation}
Furthermore, the distance between this unique barycentre $\optim{\mPoint}_{\rotsSampling_\sumTotA}$ and the unique minimiser $\optim{\mPoint}$ of $\dummyFunctionA$ is bounded as
\begin{equation}
    \distance_{\manifold}(\optim{\mPoint}_{\rotsSampling_\sumTotA},\optim{\mPoint}) 
    \leq 
    2\sqrt{\frac{\sqrt[\DimInd]{\det(\Hess_{\optim{\mPoint}} \dummyFunctionA)}}{\eigenValue_{\min}(\Hess_{\optim{\mPoint}} \dummyFunctionA)}} 
    \,\radiusA_{\sumTotA}
    \quad\text{with $\radiusA_{\sumTotA} > 0$ given as in \cref{eq:EllipsRadii}.}
\label{eq:cor-bound-distance}
\end{equation}
In the above, $\eigenValue_{\min}(\Hess_{\optim{\mPoint}} \dummyFunctionA)$ is the smallest eigenvalue of the symmetric positive-definite bilinear form $\Hess_{\optim{\mPoint}} \dummyFunctionA \colon \tangent_\mPoint\manifold\times\tangent_\mPoint\manifold\to\Real$.
\end{theorem}


\begin{remark}
    \revA{\cref{thm:new-num-non-zero-coeffs,thm:distance-bound-post-lifting} should be read and interpreted with care. That is, if the cost function is $\epsilon$-close to the minimum function value somewhere else, one might need an unrealistically large sampling set before the behaviour in the Theorems kick in. This will indeed be an issue if this happens far away from the global minimiser. If this happens closer by, e.g., as in the toy example \cref{fig:toy-ESL-weights}, the points with non-zero coefficients will be the ones in the lower envelope of the function. Although this case is not covered by the theory, we expect similar behaviour in the sense that the method picks up on the Hessian of the lower envelope. In that case, the regularised solution is still expected to be a good estimation.}
\end{remark}

\subsection{Algorithm}
\label{sec:lift-algo}
From Karcher \cite[Theorem~1.2]{karcher1977riemannian} we know that the problem \cref{eq:thm-barycentre} restricted to a small enough geodesic ball around the global minimiser $\optim{\mPoint}$ -- rather than the full manifold $\manifold$ -- is actually convex. In other words, the theoretical results suggest that for $1/(\Dim(\manifold)+1) < \exponentA <2/\Dim(\manifold)$ and a proper choice of regularisation parameter $\regParam>0$ and sampling set $\rotsSampling\subset\manifold$ one can use a Riemannian gradient descent scheme -- initialised with the point $\mPoint\in \rotsSampling$ with the highest value for~$\optim{\rotsDensCoef}_{\mPoint}$ -- to solve the optimisation in \cref{eq:lifted-projection-discr} and thereby compute the Riemannian barycentre for approximating the solution to the non-convex minimisation problem in \cref{eq:model-problem}.
If the gradient step size is $1/2$, then this scheme reads as 
\begin{equation}
\mPoint^{\iterInd+1} 
  :=\exp _{\mPoint^{\iterInd}}\!\Bigl(
      \sum_{\mPoint \in \mathcal{X}} \optim{\rotsDensCoef}_{\mPoint} \log_{\mPoint^{\iterInd}}(\mPoint)
    \Bigr).
\end{equation}
To see this, simply note that the gradient of $\distance_{\manifold}(\mPoint, \,\cdot\,)^2$ at some point $\mPointB \in \manifold$ is given as
\[ \Grad\bigl( 
     \distance_{\manifold}(\mPoint, \,\cdot\,)^2
     \bigr)\big|_{\mPointB} 
= -2 \log_\mPointB(\mPoint)
\quad\text{for $\mPointB\in \manifold$.}
\]
Combined with the computation of the coefficients, the full optimisation scheme 
is outlined in \cref{alg:l2-scaled-lifting}. We will refer to this approach as \emph{ellipsoidal support lifting} (ESL).

\begin{algorithm}[h!]
\caption{The ellipsoidal support lifting (ESL) scheme}
\label{alg:l2-scaled-lifting}
\begin{algorithmic}
\STATE{\textit{Initialisation}: $\rotsSampling\subset\manifold$, $\DimInd:= \dim(\manifold)$, $\regParam>0$, $\exponentA\in(\frac{1}{\DimInd+1},\frac{2}{\DimInd})$, $\iterInd:=0$}
\STATE{$\optim{\rotsDensCoef}=\Pi_{\simplex^{\numRots}}\left(-\frac{\numRots^\exponentA}{\regParam}\dummyFunctionA_\rotsSampling\right)$}
\STATE{$\mPoint^0:=\displaystyle{\argmax_{\mPoint\in \rotsSampling}} \,\optim{\rotsDensCoef}_{\mPoint}$}
\WHILE{not converged}
\STATE{$\mPoint^{\iterInd+1}:=\exp _{\mPoint^{\iterInd}}
\!\left(\sum_{\mPoint \in \mathcal{X}} 
\optim{\rotsDensCoef}_{\mPoint} \log_{\mPoint^{\iterInd}}(\mPoint)\right)$}
\STATE{$\iterInd:=\iterInd+1$}
\ENDWHILE
\end{algorithmic}
\end{algorithm}

Note that it remains to have a practically realisable scheme for choosing/constructing $\rotsSampling$ and $\regParam$. Without going too much into detail, the sampling set $\rotsSampling$ -- or rather the sequence of sampling sets  $(\rotsSampling_\sumTotA)_{\sumTotA=1}^\infty$ -- can be constructed through mesh refinement, i.e., one starts from some small uniformly distributed set of points $\rotsSampling_0$, constructs a triangulation, and applies some refinement rule. Choosing such a $\rotsSampling_0$ and refinement rule is highly dependent on the manifold of interest and the application. Both will be addressed below in the Cryo-EM case study in \cref{sec:application-to-cryo-EM}. Then, the most challenging aspect in constructing $\regParam$ will be determining the (determinant of the) Riemannian Hessian. In exceptional cases the Hessian is known, but in general one will need to resort to approximations. One final consequence of our theoretical analysis is the following result, which provides such an approximation for the $\regParam$ in \cref{eq:thm-new-num-non-zero-coeffs-gamma} required by the main results.
\begin{proposition}
\label{prop:choice-regulariser-lifting}
\revA{Consider the minimisation in \cref{eq:lifted-function-discr} with the same assumptions as in \cref{thm:new-num-non-zero-coeffs}, except that $(\rotsSampling_\sumTotA)_{\sumTotA=1}^\infty$ is a local low-discrepancy sequence \revA{for all pairs of parameters $(\eta, \dummyIntegerB)$ where $\dummyIntegerB \in \{\dummyIntegerB_0,\bigl( \frac{\DimInd + 1 + \epsilon_1}{\DimInd + 2}\bigr)^{\frac{\DimInd}{2}}\dummyIntegerB_0, \bigl(\frac{\DimInd + 3 - \epsilon_2}{\DimInd + 2}\bigr)^{\frac{\DimInd}{2}}\dummyIntegerB_0\}$}. Furthermore, define for each $\rotsSampling_\sumTotA$ the labeled points $\mPoint^\sumTotA_{(1)}, \ldots,\mPoint^\sumTotA_{(|\rotsSampling_\sumTotA|)}$ sorted such that  $\dummyFunctionA(\mPoint^\sumTotA_{(1)})\leq \ldots\leq \dummyFunctionA(\mPoint^\sumTotA_{(|\rotsSampling_\sumTotA|)})$, and let 
\begin{equation}
    \epsilon^-\in \bigl(0,\revA{\frac{(\DimInd + 2)(\DimInd + \epsilon_1) - \DimInd(\DimInd + 3 - \epsilon_2)}{\DimInd + 2}}\bigr) \quad \text{and} \quad  \epsilon^{+} \in \bigl(0, \frac{\DimInd(\DimInd + 1+ \epsilon_1) - (\DimInd +2 )(\DimInd - \epsilon_2)}{\DimInd + 2}\bigr).
\end{equation}
}

Then, \revA{there exists $\sumTotA'\in \Natural$ such that for all $\sumTotA\geq \sumTotA'$}
\begin{equation}
     \revA{\frac{1 + \epsilon^-}{2}\regParam \leq \frac{1}{2}\dummyIntegerB_0 |\rotsSampling_\sumTotA|^{\frac{2 + 2\exponentA}{\DimInd + 2}} \biggl(\dummyFunctionA(\mPoint^\sumTotA_{(\dummyIntegerB_\sumTotA+1)})-\frac{1}{\dummyIntegerB_\sumTotA} \sum_{\sumIndB=1}^{\dummyIntegerB_\sumTotA} \dummyFunctionA(\mPoint^\sumTotA_{(\sumIndB)})\biggr) 
\leq \frac{3 - \epsilon^+}{2}\regParam,}
\end{equation}
where $\dummyIntegerB_\sumTotA:= \lfloor \dummyIntegerB_0|\rotsSampling_\sumTotA|^{\frac{2 - \DimInd\exponentA}{\DimInd+2}} \rfloor$.
\end{proposition}

\begin{remark}
    \revA{Note that for $\epsilon_1,\epsilon_2$ close to $1$ we can choose $\epsilon^-$ and $\epsilon^+$ to be close to $1$, i.e., the approximation becomes tight.}
\end{remark}


\section{Proofs of the main results}
\label{sec:so3-lifting-analysis}
This part focuses on the proofs of \cref{thm:new-num-non-zero-coeffs,thm:distance-bound-post-lifting,prop:choice-regulariser-lifting}. 
The starting point is to elaborate on the notion of an local low-discrepancy sequence (\cref{sec:local-low-discr}),
which is followed by the proof of \cref{thm:new-num-non-zero-coeffs} (\cref{sec:proof-sparsity}), i.e., showing that \cref{eq:lifted-function-discr-argmin} generates a sparse measure concentrated on an ellipsoid. 
The latter also yields techniques for proving  \cref{thm:distance-bound-post-lifting} (\cref{sec:lifting-wellposed}), which gives well-posedness of the barycentre and error bounds on \cref{alg:l2-scaled-lifting}, and for proving \cref{prop:choice-regulariser-lifting} (\cref{sec:approx-regu}).



\revA{Before that, we would like to point out that all of the above relies on the} following lemma, \revA{which} allows us to approximate ellipsoid volumes and provides a second integration approximation.

\begin{lemma}
\label{lem:ellipsoid-integration-rules}
Let $\manifold$ be a Riemannian manifold with dimension $\DimInd:= \dim (\manifold)$, $\mPoint\in \manifold$, and let $\dummyBilForm \colon \tangent_\mPoint\manifold\times \tangent_\mPoint\manifold\to \Real$ be a symmetric positive definite bi-linear form.
Then, for radii $\radiusB\to 0$ we have 
\begin{equation}
\operatorname{vol}\bigl(\ellipsoid_\radiusB^\dummyBilForm(\mPoint)\bigr) 
  = \ballVol_{\DimInd}\radiusB^{\DimInd} + \mathcal{O}(\radiusB^{\DimInd+1})
\label{eq:ellips-integral-1}
\end{equation}
and
\begin{equation}
\int_{\ellipsoid_\radiusB^\dummyBilForm(\mPoint)}\!\!
  \dummyBilForm\bigl(
    \log_\mPoint(\mPointB),\log_\mPoint(\mPointB)
  \bigr) \mathrm{d} \mPointB 
= \sqrt[\DimInd]{\det(\dummyBilForm)}\,
  \frac{\ballVol_{\DimInd}\DimInd}{\DimInd+2}\radiusB^{\DimInd+2}  
  + \mathcal{O}(\radiusB^{\DimInd+3}).
\label{eq:ellips-integral-2}
\end{equation}
\end{lemma}
\begin{proof}
Let $\{\mTVector_{\coordInd}\}_{\coordInd=1}^\DimInd\subset \tangent_\mPoint\manifold$ be an orthonormal basis of eigenvectors of $\dummyBilForm$ with corresponding eigenvalues $\eigenValue_\coordInd(\dummyBilForm)$. Furthermore, let $\chart_{\ellipsoid} \colon \ellipsoid_\radiusB^\dummyBilForm(\mPoint) \to \Real^\DimInd$ be the \emph{re-scaled} normal coordinate chart generated by $\{\mTVector_{\coordInd}\}_{\coordInd=1}^\DimInd$, with the $\coordInd$:th component given by
\begin{equation}
\chart_\ellipsoid(\mPointB)_\coordInd
:= \frac{\sqrt{\eigenValue_{\coordInd}(\dummyBilForm)}}
     {\sqrt[2\DimInd]{\operatorname{det}(\dummyBilForm)}}\bigl(\mTVector_\coordInd,
       \log_{\mPoint}(\mPointB)
     \bigr)_{\mPoint}.
\end{equation}
The inverse $\chart^{-1}_{\ellipsoid} \colon \chart_\ellipsoid(\ellipsoid^\dummyBilForm_\radiusB(\mPoint))\to \ellipsoid_\radiusB(\mPoint)$ is then
\begin{equation}
  \chart^{-1}_\ellipsoid(\ePoint)= 
  \exp_{\mPoint}\Bigl( 
    \sum_{\coordInd=1}^\DimInd 
      \frac{\sqrt[2\DimInd]{\operatorname{det}(\dummyBilForm)}}
      {\sqrt{\eigenValue_{\coordInd}(\dummyBilForm)}}
      \, \ePoint^\coordInd \mTVector_\coordInd
  \Bigr).
\end{equation}
Note in particular that $\chart_\ellipsoid(\ellipsoid^\dummyBilForm_\radiusB(\mPoint))=\mathbb{B}_\radiusB := \bigl\{\ePoint\in \Real^\DimInd \mid \|\ePoint\|_2 < \radiusB \bigr\}$.
        
Next, we claim that $\det(\differential_{0}\chart_\ellipsoid^{-1})=1$.
To see this, note first that 
\begin{align}
\differential_{0} \chart_\ellipsoid^{-1} \Bigl[\frac{\partial}{\partial \ePoint^\coordInd}\Bigr]
\dummyFunctiona 
= \frac{\mathrm{d}}{\mathrm{d}\tim} \dummyFunctiona\biggl( 
\exp_{\mPoint}\Bigl(
  \frac{\sqrt[2\DimInd]{\operatorname{det}(\dummyBilForm)}}{\sqrt{\eigenValue_{\coordInd}(\dummyBilForm)}} 
  \tim \mTVector_\coordInd
\Bigr) \biggr) \bigg|_{\tim=0}\!\!
= \frac{\sqrt[2\DimInd]{\operatorname{det}(\dummyBilForm)}}{\sqrt{\eigenValue_{\coordInd}(\dummyBilForm)}} \mTVector_\coordInd \dummyFunctiona
\quad\text{for $\dummyFunctiona\in C^{\infty}(\manifold)$.}
\end{align}
Hence, $\differential_{0} \chart_\ellipsoid^{-1} \bigl[\frac{\partial}{\partial \ePoint^\coordInd}\bigr] = \frac{\sqrt[2\DimInd]{\operatorname{det}(\dummyBilForm)}}{\sqrt{\eigenValue_{\coordInd}(\dummyBilForm)}} \mTVector_\coordInd$ and we have a diagonal matrix in the $\{\mTVector_{\coordInd}\}_{\coordInd=1}^\DimInd\subset \tangent_\mPoint\manifold$ basis. 
This gives 
\begin{equation}
\det(\differential_{0} \chart_\ellipsoid^{-1}) 
= \prod_{\coordInd=1}^\DimInd \frac{\sqrt[2\DimInd]{\operatorname{det}(\dummyBilForm)}}{\sqrt{\eigenValue_{\coordInd}(\dummyBilForm)}} 
= \frac{\sqrt[2]{\det(\dummyBilForm)}}{\prod_{\coordInd=1}^\DimInd\sqrt{\eigenValue_{\coordInd}(\dummyBilForm)}} 
= 1
\quad\text{as $\det(\dummyBilForm)=\prod_{\coordInd=1}^\DimInd\eigenValue_{\coordInd}(\dummyBilForm)$,}
\end{equation}
which proves that $\det(\differential_{0}\chart_\ellipsoid^{-1})=1$.

Using this result, the proof of the claim in \cref{eq:ellips-integral-1} follows directly from the following calculation:
\begin{align*}
\operatorname{vol}(\ellipsoid_\radiusB^\dummyBilForm(\mPoint)) 
&= \int_{\ellipsoid_\radiusB^\dummyBilForm(\mPoint)}
   \!\! \mathrm{d} \mPointB 
= \int_{\mathbb{B}_\radiusB}    
   \bigl|\det(\differential_{\ePoint} \chart_\ellipsoid^{-1}) \bigr| 
   \mathrm{d}\ePoint
= \int_{\mathbb{B}_\radiusB}  \det(\differential_{\ePoint} \chart_\ellipsoid^{-1}) \mathrm{d}\ePoint
\\
& = \int_{\mathbb{B}_\radiusB}
    \det(\differential_{0} \chart_\ellipsoid^{-1}) 
    + \mathcal{O}(\radiusB) \mathrm{d}\ePoint
= \operatorname{vol}(\mathbb{B}_\radiusB) 
    + \mathcal{O}(\radiusB^{\DimInd+1})
= \ballVol_\DimInd \radiusB^\DimInd 
    + \mathcal{O}(\radiusB^{\DimInd+1}),
\end{align*}
where we used in the third equality the fact that for small enough radius $\radiusB$ the determinant must have a constant (positive) sign. The fourth equality stems from a Taylor expansion around 0 -- this is justified, as the differential of the local chart varies smoothly and so will its eigenvalues, and hence also its determinant.

Similarly, the proof of equality \cref{eq:ellips-integral-2} reads as follows:
\begin{multline*}
\int_{\ellipsoid_\radiusB^\dummyBilForm(\mPoint)} \dummyBilForm\bigl(\log_\mPoint(\mPointB),\log_\mPoint(\mPointB)\bigr) \mathrm{d} \mPointB 
= \int_{\eball_\radiusB} \dummyBilForm\Bigl(\log_\mPoint(\chart_\ellipsoid^{-1} \bigl(\ePoint)\bigr),\log_\mPoint\bigl(\chart_\ellipsoid^{-1} \bigl(\ePoint)\bigr)\Bigr) \det(\differential_\ePoint\chart_\ellipsoid^{-1}) \mathrm{d} \ePoint
\\
= \sqrt[\DimInd]{\det(\dummyBilForm)} \int_{\eball_\radiusB} \|\ePoint\|_2^2 \det(\differential_\ePoint\chart_\ellipsoid^{-1}) \mathrm{d} \ePoint
= \sqrt[\DimInd]{\det(\dummyBilForm)} \int_{\eball_\radiusB} \Bigl[ \|\ePoint\|_2^2 \det(\differential_0\chart_\ellipsoid^{-1}) + \mathcal{O}(\radiusB^{3}) \Bigr] \mathrm{d} \ePoint
\\
= \sqrt[\DimInd]{\det(\dummyBilForm)} \int_{\eball_\radiusB}  \|\ePoint\|_2^2 \mathrm{d} \ePoint + \mathcal{O}(\radiusB^{\DimInd+3})
= \sqrt[\DimInd]{\det(\dummyBilForm)} \frac{\ballVol_{\DimInd}\DimInd}{\DimInd+2}\radiusB^{\DimInd+2} + \mathcal{O}(\radiusB^{\DimInd+3}).
\end{multline*}
The 2nd equality above follows from  $\dummyBilForm\Bigl(\log_\mPoint(\chart_\ellipsoid^{-1}(\ePoint)),\log_\mPoint\bigl(\chart_\ellipsoid^{-1}(\ePoint)\bigr)\Bigr) = \sqrt[\DimInd]{\det(\dummyBilForm)} \|\ePoint\|_2^2$.
\end{proof}

\subsection{A note on local low-discrepancy sequences}
\label{sec:local-low-discr}

    As \cref{def:eta-local-low-discrepancy-sequence} is rather technical, we will briefly \revA{deconstruct} it. There are four main choices to justify: the choice of the sets we want to integrate over, the choice of functions we integrate, the choice of the decay rates on the discrepancies we ask for, and the choice of the bounds for the parameter $\exponentA$ in \cref{def:eta-local-low-discrepancy-sequence}. Existence of such sequences will be discussed afterwards.

\paragraph{Choice of sets}
First consider the sets $\ellipsoid_{\radiusA_\sumTotA}^\dummyBilForm(\mPoint)$.
As indicated before, the solution to \cref{eq:lifted-function-discr} yields a sparse measure that is supported around the global minimiser of the function $\dummyFunctionA$. 
The shape of the region containing the support is determined by the Riemannian Hessian, if it is a symmetric positive-definite bi-linear form, and its size scales with radius 
\[ \radiusA_\sumTotA= \sqrt[\DimInd]{\dummyIntegerB\frac{\operatorname{vol}(\manifold)}{\ballVol_\DimInd}} |\rotsSampling_\sumTotA|^{-\frac{1+\exponentA}{\DimInd+2}}. \]
In other words, we are only concerned with local integration.

\paragraph{Choice of functions}
This naturally motivates the choice of functions.
Since the sampling sets only have to integrate accurately in a local neighbourhood, they only need to be accurate enough on local approximations of functions. In particular, 
the integration scheme only needs to be able to integrate a Taylor expansion of a function around its minimiser up to second-order term. In other words, we only need to be able to integrate over constants \cref{eq:eta-local-low-discrepancy-sequence-1} and Riemannian Hessians \cref{eq:eta-local-low-discrepancy-sequence-2}, which come in the form of a symmetric positive-definite bi-linear forms, as the gradient is zero at the global minimiser.


\paragraph{Choice of decay rates and $\exponentA$ bounds}
The rates in \cref{eq:eta-local-low-discrepancy-sequence-1,eq:eta-local-low-discrepancy-sequence-2} and the $\exponentA$ bounds are chosen to ensure that, as the sample set size $|\rotsSampling_\sumTotA|$ increases, the number of points in the ellipsoids $\ellipsoid_{\radiusA_\sumTotA}^\dummyBilForm(\mPoint)$ approaches infinity. This is despite the fact that the size of the ellipsoids vanishes asymptotically. More precisely, in \cref{lem:ellipsoid-integration-rules} it \revA{was} shown that
\begin{equation}
   \operatorname{vol}(\ellipsoid_{\radiusA_\sumTotA}^\dummyBilForm(\mPoint))= \ballVol_\DimInd \radiusA_\sumTotA^\DimInd  + \mathcal{O}(\radiusA_\sumTotA^{\DimInd+1}) = \ballVol_\DimInd \radiusA_\sumTotA^\DimInd  + \mathcal{O}(|\rotsSampling_\sumTotA|^{-\frac{(\DimInd+1)(1+\exponentA)}{\DimInd+2}}),
   \label{eq:def-integration-motivation-vol}
\end{equation}
and
\begin{multline}
    \int_{\ellipsoid_{\radiusA_\sumTotA}^\dummyBilForm(\mPoint)}\!\!
  \dummyBilForm\bigl(
    \log_\mPoint(\mPointB),\log_\mPoint(\mPointB)
  \bigr) \mathrm{d} \mPointB 
= \sqrt[\DimInd]{\det(\dummyBilForm)}\,
  \frac{\ballVol_{\DimInd}\DimInd}{\DimInd+2}\radiusA_\sumTotA^{\DimInd+2}  
  + \mathcal{O}(\radiusA_\sumTotA^{\DimInd+3}) 
  \\
  = \sqrt[\DimInd]{\det(\dummyBilForm)}\,
  \frac{\ballVol_{\DimInd}\DimInd}{\DimInd+2}\radiusA_\sumTotA^{\DimInd+2}  
  + \mathcal{O}(|\rotsSampling_\sumTotA|^{-\frac{(\DimInd+3)(1+\exponentA)}{\DimInd+2}}).
  \label{eq:def-integration-motivation-Q}
\end{multline}
Note that for $1/(\DimInd+1) < \exponentA < 2/\DimInd$ the number of points that would fit in an ellipsoid of radius $\radiusA_\sumTotA= \sqrt[\DimInd]{\revA{\dummyIntegerB}\frac{\operatorname{vol}(\manifold)}{\ballVol_\DimInd}} |\rotsSampling_\sumTotA|^{-\frac{1+\exponentA}{\DimInd+2}}$ grows indefinitely. Indeed,
\begin{equation}
    \frac{\operatorname{vol}(\ellipsoid_{\radiusA_\sumTotA}^\dummyBilForm(\mPoint))}{\operatorname{vol}(\manifold)}|\rotsSampling_\sumTotA| =  \revA{\dummyIntegerB}|\rotsSampling_\sumTotA|^{\frac{2 - \DimInd\exponentA}{\DimInd+2}} + \mathcal{O}(|\rotsSampling_\sumTotA|^{\frac{1 - (\DimInd+1)\exponentA}{\DimInd+2}}),
\end{equation}
and the first term on the right hand side will grow infinitely large, whereas the second term vanishes as $\sumTotA\to\infty$, i.e., the number of points in this family of ellipsoids grows arbitrarily large even though the sets themselves contract.

\begin{remark}
    Thus, we simply ask in \cref{def:eta-local-low-discrepancy-sequence} that for an increasing amount of points in a contracting local volume, the approximations in \cref{eq:eta-local-low-discrepancy-sequence-1,eq:eta-local-low-discrepancy-sequence-2} should be as exact as the leading term in \cref{eq:def-integration-motivation-vol} respectively \cref{eq:def-integration-motivation-Q}.
\end{remark}

Although \cref{def:eta-local-low-discrepancy-sequence} -- accompanied with this motivation -- sounds reasonable, existence of such sequences \revA{on general manifolds} is still open, as analysis of error rates of sampling sets is typically developed for global integration \cite{brandolini2010quadrature}. To take a first step, \revA{a proof of existence of a local low-discrepancy sequence on $(0,1)\subset \Real$ is provided in \cref{sec:low-discr-R}. Furthermore,} the numerical experiments provided in \cref{sec:numerics} -- that show the error rates -- will give some additional empirical evidence for the local low-discrepancy property being well-defined \revA{on general manifolds}.

\subsection{Sparse and concentrated measures through relaxation}
\label{sec:proof-sparsity}
In order to prove \cref{thm:new-num-non-zero-coeffs}, we will argue that the points with non-zero coefficients must be contained in an ellipsoid and that the scaling behaviour of the ellipsoid's volume yields a bound for the scaling of number of non-zero coefficients.

The following lemma is a consequence of \Cref{lem:ellipsoid-integration-rules} and essential for the proof of \cref{thm:new-num-non-zero-coeffs}.
\begin{lemma}
\label{lem:numpoints-in-ellipsoid}
Let $\manifold$ be a finite-volume Riemannian manifold with dimension $\DimInd:= \dim (\manifold)$. Consider any point $\mPoint\in\manifold$, symmetric positive definite bi-linear form $\dummyBilForm\colon \tangent_\mPoint\manifold\times\tangent_\mPoint\manifold\to\Real$ and $\revA{\epsilon  \in (0,1)}$. 
Furthermore, assume that $(\rotsSampling_\sumTotA)_{\sumTotA=1}^\infty$ is a local low-discrepancy sequence with $1/(\DimInd+1) < \exponentA < 2/\DimInd$ \revA{and $\dummyIntegerB_0 >0$ and define the sequences of vanishing radii  $(\radiusA_\sumTotA)_{\sumTotA=1}^\infty\subset \Real_+$ by $\radiusA_\sumTotA:= \sqrt[\DimInd]{\dummyIntegerB_0\frac{\operatorname{vol}(\manifold)}{\ballVol_\DimInd}} |\rotsSampling_\sumTotA|^{-\frac{1+\exponentA}{\DimInd+2}}$}.
\revA{
\begin{enumerate}
    \item Define $(\dummyIntegerB^-_\sumTotA)_{\sumTotA=1}^\infty\subset\Natural$ by $\dummyIntegerB^-_\sumTotA := \lfloor \dummyIntegerB_0|\rotsSampling_\sumTotA|^{\frac{2 - \DimInd\exponentA}{\DimInd+2}} \rfloor$, and define the sequences of vanishing radii $(\radiusA^-_\sumTotA)_{\sumTotA=1}^\infty,(\radiusB_\sumTotA)_{\sumTotA=1}^\infty \subset \Real_+$ by $\radiusA^-_\sumTotA := \sqrt[\DimInd]{\frac{\operatorname{vol}(\manifold)}{\ballVol_\DimInd} \frac{\dummyIntegerB^-_\sumTotA}{|\rotsSampling_\sumTotA|}}$ and $\radiusB_\sumTotA := \sqrt{\frac{\DimInd+1 + \epsilon}{\DimInd+2}} \radiusA_\sumTotA$. If $(\rotsSampling_\sumTotA)_{\sumTotA=1}^\infty$ is also a local low discrepancy sequence with $\exponentA$ and $\bigl(\frac{\DimInd+1 + \epsilon}{\DimInd+2}\bigr)^{\frac{\DimInd}{2}}\dummyIntegerB_0$, then there exists $\sumTotA'\in \Natural$ such that for all $\sumTotA\geq \sumTotA'$ we have
\begin{equation}
 \dummyIntegerB_\sumTotA^- > |\ellipsoid_{\radiusB_\sumTotA}^{\dummyBilForm}(\mPoint) \cap \rotsSampling_\sumTotA|.
\label{eq:lem-J-as-points-in-ellipsoid->}
\end{equation}
\item Define $(\dummyIntegerB^+_\sumTotA)_{\sumTotA=1}^\infty\subset\Natural$ by $\dummyIntegerB^+_\sumTotA := \lceil \dummyIntegerB_0|\rotsSampling_\sumTotA|^{\frac{2 - \DimInd\exponentA}{\DimInd+2}} \rceil$, and define the sequences of vanishing radii $(\radiusA^+_\sumTotA)_{\sumTotA=1}^\infty ,(\radiusC_\sumTotA)_{\sumTotA=1}^\infty \subset \Real_+$ by $\radiusA^+_\sumTotA := \sqrt[\DimInd]{\frac{\operatorname{vol}(\manifold)}{\ballVol_\DimInd} \frac{\dummyIntegerB^+_\sumTotA}{|\rotsSampling_\sumTotA|}}$ and $\radiusC_\sumTotA := \sqrt{\frac{\DimInd+3 - \epsilon}{\DimInd+2}} \radiusA_\sumTotA$. If $(\rotsSampling_\sumTotA)_{\sumTotA=1}^\infty$ is also a local low discrepancy sequence with $\exponentA$ and $\bigl(\frac{\DimInd+3 - \epsilon}{\DimInd+2}\bigr)^{\frac{\DimInd}{2}}\dummyIntegerB_0$, then there exists $\sumTotA'\in \Natural$ such that for all $\sumTotA\geq \sumTotA'$ we have 
\begin{equation}
\dummyIntegerB_\sumTotA^+
  < |\ellipsoid_{\radiusC_\sumTotA}^{\dummyBilForm}(\mPoint) \cap \rotsSampling_\sumTotA|.
\label{eq:lem-J-as-points-in-ellipsoid-<}
\end{equation}
\end{enumerate}
}

\end{lemma}

\begin{proof}
\revA{For the given $\epsilon\in (0,1)$, $1/(\DimInd+1) < \exponentA < 2/\DimInd$ and $\dummyIntegerB_0 >0$ f}ix a local low-discrepancy sequence $(\rotsSampling_\sumTotA)_{\sumTotA=1}^\infty \subset S(\manifold)$ \revA{fulfilling the conditions in the statement}.



To see \cref{eq:lem-J-as-points-in-ellipsoid->,eq:lem-J-as-points-in-ellipsoid-<},
first we \revA{note} that by \cref{eq:ellips-integral-1} in \cref{lem:ellipsoid-integration-rules}
\begin{multline}
    \frac{\operatorname{vol}(\ellipsoid_{\radiusA_\sumTotA^{\revA{\pm}}}^{\dummyBilForm}\bigl(\optim{\mPoint}\bigr))}{\operatorname{vol}(\manifold)}|\rotsSampling_\sumTotA| =
    \frac{\ballVol_{\DimInd} {(\radiusA_\sumTotA^{\revA{\pm}})}^{\DimInd} + \mathcal{O}({(\radiusA_\sumTotA^{\revA{\pm}})}^{\DimInd+1})}{\operatorname{vol}(\manifold)}|\rotsSampling_\sumTotA|
    \\
    = \dummyIntegerB_\sumTotA^{\revA{\pm}} + \mathcal{O}(\dummyIntegerB_\sumTotA^{\revA{\pm}} \radiusA_\sumTotA^{\revA{\pm}}) = \dummyIntegerB_\sumTotA^{\revA{\pm}} + \mathcal{O}( |\rotsSampling_\sumTotA|^{\frac{1 - (\DimInd+1)\exponentA}{\DimInd+2}})
    \label{eq:vol-ellipsoid-J-relation}
\end{multline}
and the right-hand side approaches $\dummyIntegerB_\sumTotA^{\revA{\pm}}$ arbitrarily close as the remainder term vanishes because $\exponentA> \frac{1}{\DimInd+1}$. This tells us that the ellipsoid can hold as many points as $\dummyIntegerB_\sumTotA^{\revA{\pm}}$. \revA{
Next, \cref{eq:vol-ellipsoid-J-relation} combined with $\exponentA < \frac{2}{\DimInd}$ also gives that $\dummyIntegerB_\sumTotA^{\revA{\pm}} \to \infty$ as $\sumTotA\to\infty$ and
\begin{equation}
    \lim_{\sumTotA\to\infty} \frac{\operatorname{vol}\bigl(\ellipsoid_{\radiusA_\sumTotA^{\revA{\pm}}}^\dummyBilForm(\optim{\mPoint})\bigr)}{\operatorname{vol}(\ellipsoid_{\radiusA_\sumTotA}^\dummyBilForm(\optim{\mPoint}))} =  \lim_{\sumTotA\to\infty}\frac{\dummyIntegerB_\sumTotA^{\revA{\pm}}}{\dummyIntegerB_0|\rotsSampling_\sumTotA|^{\frac{2 - \DimInd\exponentA}{\DimInd+2}}} = 1.
    \label{eq:-lem-upper-lower-J-ellipsoids-lim1}
\end{equation}
Combining \cref{eq:-lem-upper-lower-J-ellipsoids-lim1} with the local low-discrepancy property \cref{eq:eta-local-low-discrepancy-sequence-1} of $(\rotsSampling_\sumTotA)_{\sumTotA=1}^\infty$ gives that $\dummyIntegerB_\sumTotA^{\revA{\pm}}$ is close to $|\ellipsoid_{\radiusA_\sumTotA}^\dummyBilForm(\optim{\mPoint})\cap \rotsSampling_\sumTotA|$ for large $\sumTotA$, because $|\ellipsoid_{\radiusA_\sumTotA^-}^\dummyBilForm(\optim{\mPoint})\cap \rotsSampling_\sumTotA| \leq |\ellipsoid_{\radiusA_\sumTotA}^\dummyBilForm(\optim{\mPoint})\cap \rotsSampling_\sumTotA| \leq |\ellipsoid_{\radiusA_\sumTotA^+}^\dummyBilForm(\optim{\mPoint})\cap \rotsSampling_\sumTotA|$. Now, for showing \cref{eq:lem-J-as-points-in-ellipsoid->,eq:lem-J-as-points-in-ellipsoid-<} is it sufficient to show that
\begin{equation}
    \lim_{\sumTotA\to\infty} \frac{\operatorname{vol}\bigl(\ellipsoid_{\radiusB_\sumTotA}^\dummyBilForm(\optim{\mPoint})\bigr)}{\operatorname{vol}(\ellipsoid_{\radiusA_\sumTotA}^\dummyBilForm(\optim{\mPoint}))}  < 1, \quad \text{and} \quad \lim_{\sumTotA\to\infty} \frac{\operatorname{vol}\bigl(\ellipsoid_{\radiusC_\sumTotA}^\dummyBilForm(\optim{\mPoint})\bigr)}{\operatorname{vol}(\ellipsoid_{\radiusA_\sumTotA}^\dummyBilForm(\optim{\mPoint}))} > 1,
    \label{eq:-lem-upper-lower-J-ellipsoids-lim<>1}
\end{equation}
Indeed, by the additional local low discrepancy assumptions on $(\rotsSampling_\sumTotA)_{\sumTotA=1}^\infty$ in (1) and (2) in combination with \cref{eq:eta-local-low-discrepancy-sequence-1} and by the fact that the amount of points in both ellipsoids grows infinitely large -- even though the ellipsoid contracts as $\sumTotA\to \infty$ -- because $\exponentA<\frac{2}{\DimInd}$, this tells us that the amount of points in the smaller respectively larger ellipsoid will grow apart from $|\ellipsoid_{\radiusA_\sumTotA}^{\dummyBilForm}(\mPoint) \cap \rotsSampling_\sumTotA|$ and hence from $\dummyIntegerB_\sumTotA^{\revA{\pm}}$ for large $\sumTotA$. This then yields \cref{eq:lem-J-as-points-in-ellipsoid->,eq:lem-J-as-points-in-ellipsoid-<}.

By invoking \cref{eq:ellips-integral-1} in \cref{lem:ellipsoid-integration-rules} once more, the left hand inequality in \cref{eq:-lem-upper-lower-J-ellipsoids-lim<>1} follows from
\begin{equation}
 \lim_{\sumTotA\to\infty} \frac{\operatorname{vol}\bigl(\ellipsoid_{\radiusB_\sumTotA}^\dummyBilForm(\optim{\mPoint})\bigr)}{\operatorname{vol}(\ellipsoid_{\radiusA_\sumTotA}^\dummyBilForm(\optim{\mPoint}))} = \Bigl ( \frac{\DimInd+1 + \epsilon}{\DimInd+2}\Bigr)^{\frac{\DimInd}{2}} < 1.
\end{equation}
which proves \cref{eq:lem-J-as-points-in-ellipsoid->} in (1) and the right hand inequality in \cref{eq:-lem-upper-lower-J-ellipsoids-lim<>1} follows from
\begin{equation}
\lim_{\sumTotA\to\infty} \frac{\operatorname{vol}\bigl(\ellipsoid_{\radiusC_\sumTotA}^\dummyBilForm(\optim{\mPoint})\bigr)}{\operatorname{vol}(\ellipsoid_{\radiusA_\sumTotA}^\dummyBilForm(\optim{\mPoint}))} = \Bigl ( \frac{\DimInd+3 - \epsilon}{\DimInd+2}\Bigr)^{\frac{\DimInd}{2}} > 1,
\end{equation}
which proves \cref{eq:lem-J-as-points-in-ellipsoid-<} in (2). 
}
\end{proof}



We are now ready to prove \cref{thm:new-num-non-zero-coeffs}.


\begin{proof}[Proof of \cref{thm:new-num-non-zero-coeffs}]
\revA{For the given $\epsilon_1, \epsilon_2, \epsilon_3, \epsilon_4 \in (\frac{\DimInd}{2\DimInd + 2},1)$, $1/(\DimInd+1) < \exponentA < 2/\DimInd$ and $\dummyIntegerB_0 \geq 3^{-\frac{\DimInd}{\DimInd +2}}$ f}ix a local low-discrepancy sequence $(\rotsSampling_\sumTotA)_{\sumTotA=1}^\infty \subset S(\manifold)$ \revA{fulfilling the conditions in the statement}. Define $(\optim{\rotsDensCoef}_\sumTotA)_{\sumTotA=1}^\infty \subset \ell_0$ as the the zero-padded solution operator of \cref{eq:lifted-function-discr}, i.e.,
\begin{equation}
    \optim{\rotsDensCoef}_\sumTotA:= \biggl( \Pi_{\simplex^{|\rotsSampling_\sumTotA|}}\Bigl(- \frac{|\rotsSampling_\sumTotA|^{\exponentA}}{\regParam}\dummyFunctionA_{\rotsSampling_\sumTotA} \Bigr), 0,\ldots\biggr).
\end{equation}
Additionally, for any $\sumTotA \in \Natural$ define the points $\mPoint^\sumTotA_{(1)}, \ldots,\mPoint^\sumTotA_{(|\rotsSampling_\sumTotA|)}$ by labeling the elements of $\rotsSampling_\sumTotA$ such that  $\dummyFunctionA(\mPoint^\sumTotA_{(1)})\leq \ldots\leq \dummyFunctionA(\mPoint^\sumTotA_{(|\rotsSampling_\sumTotA|)})$. 

Note that with these definitions, the non-zero coefficients in $\optim{\rotsDensCoef}_\sumTotA$ will be at the indices corresponding to the $\|\optim{\rotsDensCoef}_\sumTotA\|_0$ smallest components of $\dummyFunctionA_{\rotsSampling_\sumTotA}$.

\vspace{0.5cm}




(i) For the first step we will (i-a) prove the upper bound in \cref{eq:thm-J-bounds} and (i-b) prove that the points corresponding to the non-zero coefficients lie in the ellipsoid $\ellipsoid_{\radiusA_\sumTotA}^{\dummyBilForm}\bigl(\optim{\mPoint}\bigr)$ with $\dummyBilForm:= \Hess_{\optim{\mPoint}} \dummyFunctionA$ with $\radiusA_\sumTotA$ as in \cref{eq:EllipsRadii}.

\vspace{0.5cm}

(i-a) For the upper bound in \cref{eq:thm-J-bounds} we will show that, for large enough $\sumTotA$, the number of non-zero coefficients $\|\optim{\rotsDensCoef}_\sumTotA\|_0$ is strictly bounded by $\|\optim{\rotsDensCoef}_\sumTotA\|_0 < \dummyIntegerB_\sumTotA^{\revA{+}}$, where 
\begin{equation}
\label{eq:thm-dummyinteger_sumtotA}
    \dummyIntegerB_\sumTotA^{\revA{+}} := \lceil \dummyIntegerB_0|\rotsSampling_\sumTotA|^{\frac{2 - \DimInd\exponentA}{\DimInd+2}} \rceil.
\end{equation}
In particular, we will assume that $\|\optim{\rotsDensCoef}_\sumTotA\|_0 \geq \dummyIntegerB_\sumTotA^{\revA{+}}$. If we can show that, for large $\sumTotA$, 
\begin{equation}
    \label{eq:thm-num-nonzero-contradiction}
    \biggl\|
    \biggl(- \frac{|\rotsSampling_\sumTotA|^{\exponentA}}{\regParam}\dummyFunctionA_{\rotsSampling_\sumTotA} -\frac{1}{\dummyIntegerB_\sumTotA^{\revA{+}}} \Bigl(\sum_{j=1}^{\dummyIntegerB_\sumTotA^{\revA{+}}} - \frac{|\rotsSampling_\sumTotA|^{\exponentA}}{\regParam} \dummyFunctionA(\mPoint^\sumTotA_{(\sumIndB)})\Bigr) + \frac{1}{\dummyIntegerB_\sumTotA^{\revA{+}}}\biggr)_{+}
    \biggr\|_0 < \dummyIntegerB_\sumTotA^{\revA{+}}
\end{equation}
holds, then also
\begin{equation}
    \biggl\|
    \biggl(- \frac{|\rotsSampling_\sumTotA|^{\exponentA}}{\regParam}\dummyFunctionA_{\rotsSampling_\sumTotA} -\frac{1}{\dummyIntegerB_\sumTotA^{\revA{+}}} \Bigl(\sum_{j=1}^{\dummyIntegerB_\sumTotA^{\revA{+}}} - \frac{|\rotsSampling_\sumTotA|^{\exponentA}}{\regParam} \dummyFunctionA(\mPoint^\sumTotA_{(\sumIndB)})\Bigr) + \frac{1}{\dummyIntegerB_\sumTotA^{\revA{+}}}\biggr)_{+}
    \biggr\|_0 
    < \|\optim{\rotsDensCoef}_\sumTotA\|_0
\end{equation}
holds by assumption. However, the latter is equivalent to $\|\optim{\rotsDensCoef}_\sumTotA\|_0 < \dummyIntegerB_\sumTotA^{\revA{+}}$ by \cref{rem:monotonicity-of-projection-correction}, which yields our contradiction.

\vspace{0.5cm}
It remains to show inequality \cref{eq:thm-num-nonzero-contradiction}. For any point $\mPoint \in \rotsSampling_\sumTotA$,
\begin{equation}
    \label{eq:pluscoeffzero}
    \biggl(
    - \frac{|\rotsSampling_\sumTotA|^{\exponentA}}{\regParam}\dummyFunctionA(\mPoint) -\frac{1}{\dummyIntegerB_\sumTotA^{\revA{+}}} \Bigl(\sum_{j=1}^{\dummyIntegerB_\sumTotA^{\revA{+}}} - \frac{|\rotsSampling_\sumTotA|^{\exponentA}}{\regParam} \dummyFunctionA(\mPoint^\sumTotA_{(\sumIndB)})\Bigr) + \frac{1}{\dummyIntegerB_\sumTotA^{\revA{+}}}
    \biggr)_{+} = 0
\end{equation}
is equivalent to
\begin{equation}
   \Bigl(- \frac{|\rotsSampling_\sumTotA|^{\exponentA}}{\regParam}\dummyFunctionA(\mPoint)\Bigr) - \frac{1}{\dummyIntegerB_\sumTotA^{\revA{+}}} \Bigl(\sum_{\sumIndB=1}^{\dummyIntegerB_\sumTotA^{\revA{+}}} - \frac{|\rotsSampling_\sumTotA|^{\exponentA}}{\regParam} \dummyFunctionA(\mPoint^\sumTotA_{(\sumIndB)}) \Bigr) + \frac{1}{\dummyIntegerB_\sumTotA^{\revA{+}} } \leq 0 ,
\end{equation}
or equally
\begin{equation}
  \dummyFunctionA(\mPoint) - \frac{1}{\dummyIntegerB_\sumTotA^{\revA{+}}} \sum_{\sumIndB=1}^{\dummyIntegerB_\sumTotA^{\revA{+}}} \dummyFunctionA(\mPoint^\sumTotA_{(\sumIndB)})\geq \frac{\regParam}{\dummyIntegerB_\sumTotA^{\revA{+}} |\rotsSampling_\sumTotA|^{\exponentA}}.
  \label{eq:Jm'-ineq}
\end{equation}



We will show that, for large enough $\sumTotA$, too many points $\mPoint$ -- that is more than $|\rotsSampling_\sumTotA|- \dummyIntegerB_\sumTotA^{\revA{+}}$ -- satisfy the last inequality \cref{eq:Jm'-ineq}; hence there are less than $\dummyIntegerB_\sumTotA^{\revA{+}}$ points $\mPoint \in \rotsSampling_\sumTotA$ for which the left side of \cref{eq:pluscoeffzero} can be non-zero, which yields \cref{eq:thm-num-nonzero-contradiction} and therefore the claim.

It remains to show that more than $|\rotsSampling_\sumTotA|- \dummyIntegerB_\sumTotA^{\revA{+}}$ points $\mPoint$ satisfy \cref{eq:Jm'-ineq} for large enough $\sumTotA$. To see this, first define $(\radiusA_\sumTotA^{\revA{+}})_{\sumTotA=1}^\infty \subset \Real_+$ with
\begin{equation}
\label{eq:thm-non-zero-coeffs-radius-def}
\revA{\radiusA_\sumTotA^{+} :=  \sqrt[\DimInd]{\dummyIntegerB_0\frac{\operatorname{vol}(\manifold)}{\ballVol_\DimInd}}|\rotsSampling_\sumTotA|^{-\frac{1+\exponentA}{\DimInd+2}} }, 
\end{equation}
as in \cref{lem:numpoints-in-ellipsoid} and note that
{\allowdisplaybreaks
\begin{multline}
\label{eq:JR-ineq}
\frac{\regParam}{\dummyIntegerB_\sumTotA^{\revA{+}} |\rotsSampling_\sumTotA|^{\exponentA}} 
\overset{\cref{eq:thm-dummyinteger_sumtotA}}{\leq} \frac{\regParam}{\dummyIntegerB_0|\rotsSampling_\sumTotA|^{\frac{2 - \DimInd\exponentA}{\DimInd+2}} |\rotsSampling_\sumTotA|^\exponentA} 
= \frac{\regParam}{\dummyIntegerB_0|\rotsSampling_\sumTotA|^{\frac{2 + 2\exponentA}{\DimInd+2}}}
\\
= \frac{\regParam}{\dummyIntegerB_0}\Bigl(\frac{\ballVol_\DimInd}{\dummyIntegerB_0\operatorname{vol}(\manifold)}\Bigr)^{\frac{2}{\DimInd}}\Bigl(\frac{\dummyIntegerB_0\operatorname{vol}(\manifold)}{\ballVol_\DimInd}\Bigr)^{\frac{2}{\DimInd}} |\rotsSampling_\sumTotA|^{-\frac{2 + 2\exponentA}{\DimInd+2}}
\\
\overset{\cref{eq:thm-non-zero-coeffs-radius-def}}{\revA{=}} \frac{\regParam}{\dummyIntegerB_0}\Bigl(\frac{\ballVol_\DimInd}{\dummyIntegerB_0\operatorname{vol}(\manifold)}\Bigr)^{\frac{2}{\DimInd}}(\radiusA_\sumTotA^{\revA{+}})^{2}
\overset{\cref{eq:thm-new-num-non-zero-coeffs-gamma}}{=} \frac{\sqrt[\DimInd]{\det(\Hess_{\optim{\mPoint}} \dummyFunctionA)}}{2(\DimInd+2)}(\radiusA_\sumTotA^{\revA{+}})^{2}.
\end{multline}}
Combining \cref{eq:Jm'-ineq} and \cref{eq:JR-ineq}, we see that any point $\mPoint\in \rotsSampling_\sumTotA$ that satisfies
\begin{equation}
  \dummyFunctionA(\mPoint) - \frac{1}{\dummyIntegerB_\sumTotA^{\revA{+}}} \sum_{\sumIndB=1}^{\dummyIntegerB_\sumTotA^{\revA{+}}} \dummyFunctionA(\mPoint^\sumTotA_{(\sumIndB)})\geq \frac{\sqrt[\DimInd]{\det(\Hess_{\optim{\mPoint}} \dummyFunctionA)}}{2(\DimInd+2)}(\radiusA_\sumTotA^{\revA{+}})^{2} .
\label{eq:thm-J-bound-step-2-opt-ineq}
\end{equation}
also satisfies \cref{eq:Jm'-ineq}.
We now \revA{note that $0< \frac{(\DimInd + 2)(\DimInd + \epsilon_1) - \DimInd(\DimInd + 3 - \epsilon_2)}{\DimInd + 2} < 1$ for $\epsilon_1, \epsilon_2 \in (\frac{1}{2\DimInd + 2}, 1)$,} fix $0< \epsilon^{\revA{+}}< \revA{\frac{(\DimInd + 2)(\DimInd + \epsilon_1) - \DimInd(\DimInd + 3 - \epsilon_2)}{\DimInd + 2}} < 1$ and consider the inequality
\begin{equation}
  \revA{\liminf_{\sumTotA\to \infty}} \frac{\dummyFunctionA(\mPoint) - \frac{1}{\dummyIntegerB_\sumTotA^{\revA{+}}} \sum_{\sumIndB=1}^{\dummyIntegerB_\sumTotA^{\revA{+}}} \dummyFunctionA(\mPoint^\sumTotA_{(\sumIndB)})}{(\radiusA_\sumTotA^{\revA{+}})^{2}}\geq (1 +\epsilon^{\revA{+}}) \frac{\sqrt[\DimInd]{\det(\Hess_{\optim{\mPoint}} \dummyFunctionA)}}{2(\DimInd+2)}.
\label{eq:thm-J-bound-step-2-opt-ineq-eps}
\end{equation}
Then there exists an $\sumTotA'$ such that for all $\sumTotA \geq \sumTotA'$
too many points $\mPoint\in \rotsSampling_\sumTotA$ satisfy \cref{eq:thm-J-bound-step-2-opt-ineq} -- and through \cref{eq:Jm'-ineq} also \cref{eq:thm-num-nonzero-contradiction} --, which yields our claim. The remainder of step (i-a) is devoted to showing that asymptotically too many $\mPoint\in \rotsSampling_\sumTotA$ satisfy \cref{eq:thm-J-bound-step-2-opt-ineq-eps}.



\vspace{0.5cm}

Consider the smaller ellipsoid $\ellipsoid_{\radiusB_\sumTotA^{\revA{+}}}^\dummyBilForm(\optim{\mPoint})$, again with $\dummyBilForm:= \Hess_{\optim{\mPoint}} \dummyFunctionA$, but with smaller radius $\radiusB_\sumTotA^{\revA{+}}:= \sqrt{\frac{\DimInd+1 + \epsilon_{\revA{1}}}{\DimInd+2}}\radiusA_\sumTotA^{\revA{+}}$. We will show that all points $\mPoint\in\rotsSampling_\sumTotA$ outside of $\ellipsoid_{\radiusB_\sumTotA'}^\dummyBilForm(\optim{\mPoint})$ asymptotically satisfy \cref{eq:thm-J-bound-step-2-opt-ineq-eps}. Then,
\begin{equation}
    \dummyIntegerB_\sumTotA^{\revA{+}} \revA{\geq \lfloor \dummyIntegerB_0|\rotsSampling_\sumTotA|^{\frac{2 - \DimInd\exponentA}{\DimInd+2}} \rfloor} > |\ellipsoid_{\radiusB_\sumTotA^{\revA{+}}}^{\dummyBilForm}(\optim{\mPoint}) \cap \rotsSampling_\sumTotA|
    \label{eq:Jm'discrepency-ER}
\end{equation}
holds by \cref{eq:lem-J-as-points-in-ellipsoid->} in \cref{lem:numpoints-in-ellipsoid}, which yields the claim.

For showing that for large $\sumTotA$ all points outside of the ellipsoid $\ellipsoid_{\radiusB_\sumTotA^{\revA{+}}}^\dummyBilForm(\optim{\mPoint})$ satisfy \cref{eq:thm-J-bound-step-2-opt-ineq-eps}, it suffices to show to show that the inequality \cref{eq:thm-J-bound-step-2-opt-ineq-eps} holds on the boundary $\partial \ellipsoid_{\radiusB_\sumTotA^{\revA{+}}}^\dummyBilForm(\optim{\mPoint})$. Indeed, since the function $\dummyFunctionA$ has a unique minimiser, the family of ellipsoids $\ellipsoid_{\radiusB}^\dummyBilForm(\optim{\mPoint})$ approximates a level sets of $\dummyFunctionA$ for $\radiusB$ small enough. Then, because $\radiusB_\sumTotA^{\revA{+}}\to 0$ for large enough $\sumTotA$ the points in the sampling set $\rotsSampling_\sumTotA$ outside $\ellipsoid_{\radiusB_\sumTotA^{\revA{+}}}^\dummyBilForm(\optim{\mPoint})$ must have a larger value of $\dummyFunctionA$ than any element of the boundary $\partial \ellipsoid_{\radiusB_\sumTotA^{\revA{+}}}^\dummyBilForm(\optim{\mPoint})$, which implies that these points must have a zero-valued coefficient.

 So, let $\mPoint\in \partial \ellipsoid_{\radiusB_\sumTotA^{\revA{+}}}^\dummyBilForm(\optim{\mPoint})$, i.e., we have 
\begin{equation}
\label{eq:ellipsoid-boundary-hessF}
\Hess_{\optim{\mPoint}} \dummyFunctionA \bigl(\log_{\optim{\mPoint}}(\mPoint), \log_{\optim{\mPoint}}(\mPoint)\bigr) 
  \overset{\cref{eq:Ellipsoid}}{=} \sqrt[\DimInd]{\det(\Hess_{\optim{\mPoint}} \dummyFunctionA)}\,(\radiusB_\sumTotA^{\revA{+}})^2
  = \sqrt[\DimInd]{\det(\Hess_{\optim{\mPoint}} \dummyFunctionA)}\frac{\DimInd+1 +\epsilon_{\revA{1}}}{\DimInd+2} (\radiusA_\sumTotA^{\revA{+}})^2.
\end{equation}
\revA{Now also consider the larger ellipsoid $\ellipsoid_{\radiusC_\sumTotA^{\revA{+}}}^\dummyBilForm(\optim{\mPoint})$, again with $\dummyBilForm:= \Hess_{\optim{\mPoint}} \dummyFunctionA$, but with larger radius $\radiusC_\sumTotA^{+}:= \sqrt{\frac{\DimInd+3 - \epsilon_2}{\DimInd+2}}\radiusA_\sumTotA^{\revA{+}}$.} Then we have
\begin{multline}
\revA{\liminf_{\sumTotA\to \infty}} \frac{\dummyFunctionA(\mPoint)-\frac{1}{\dummyIntegerB_\sumTotA^{\revA{+}}} \sum_{\sumIndB=1}^{\dummyIntegerB_\sumTotA^{\revA{+}}} \dummyFunctionA(\mPoint^\sumTotA_{(\sumIndB)})}{(\radiusA_\sumTotA^{\revA{+}})^{2}} 
\revA{\geq}  \revA{\liminf_{\sumTotA\to \infty}} \frac{\dummyFunctionA(\mPoint) - \revA{\frac{1}{|\ellipsoid_{\radiusC_\sumTotA^{+}}^{\dummyBilForm}(\optim{\mPoint})\cap \rotsSampling_\sumTotA|}} \sum_{\mPointB \in\ellipsoid_{\revA{\radiusC_\sumTotA^{+}}}^{\dummyBilForm}(\optim{\mPoint})\cap \rotsSampling_\sumTotA } \dummyFunctionA(\mPointB)}{(\radiusA_\sumTotA^{\revA{+}})^{2}}
\\
=  \revA{\liminf_{\sumTotA\to \infty}} \frac{\dummyFunctionA(\mPoint) - \dummyFunctionA(\optim{\mPoint}) - \revA{\frac{1}{|\ellipsoid_{\radiusC_\sumTotA^{+}}^{\dummyBilForm}(\optim{\mPoint})\cap \rotsSampling_\sumTotA|}} \Bigl(\sum_{\mPointB \in\ellipsoid_{\revA{\radiusC_\sumTotA^{+}}}^{\dummyBilForm}(\optim{\mPoint})\cap \rotsSampling_\sumTotA } \dummyFunctionA(\mPointB)- \dummyFunctionA(\optim{\mPoint}) \Bigr)}{(\radiusA_\sumTotA^{\revA{+}})^{2}},
\end{multline}
where we use the \revA{upper bound} in \cref{eq:lem-J-as-points-in-ellipsoid-<} \revA{and a monotonicity argument similar to \cref{rem:monotonicity-of-projection-correction}} to see that
\[
\revA{\frac{1}{\dummyIntegerB_\sumTotA^{\revA{+}} }}\sum_{\sumIndB=1}^{\dummyIntegerB_\sumTotA^{\revA{+}} } \dummyFunctionA(\mPoint^\sumTotA_{(\sumIndB)}) 
\revA{\leq}
\,
\revA{\frac{1}{|\ellipsoid_{\radiusC_\sumTotA^{+}}^{\dummyBilForm}(\optim{\mPoint})\cap \rotsSampling_\sumTotA|}}\sum_{\mPointB \in\ellipsoid_{\revA{\radiusC_\sumTotA^{+}}}^{\dummyBilForm}(\optim{\mPoint})\cap \rotsSampling_\sumTotA}
\!\!\!
\dummyFunctionA(\mPointB)
\]
holds since $\mPoint^\sumTotA_{(\sumIndB)}$ are points with minimal $\dummyFunctionA$-value.
Then we can use Taylor expansion \cite[(5.28)]{boumal2020introduction} 
and the fact that $\mathrm{d}\dummyFunctionA (\mTVector)=0$ for all $\mTVector\in \tangent_{\optim{\mPoint}}\manifold$  by first-order optimality conditions to approximate the terms in the numerator:




\begin{multline}
\dummyFunctionA(\mPoint) - \dummyFunctionA(\optim{\mPoint}) 
= \frac{1}{2}\Hess_{\optim{\mPoint}} \dummyFunctionA (\log_{\optim{\mPoint}}(\mPoint), \log_{\optim{\mPoint}}(\mPoint)) + \mathcal{O}((\radiusB_\sumTotA^{\revA{+}})^3) 
\\
\overset{\cref{eq:Ellipsoid}}{=} \frac{1}{2} \sqrt[\DimInd]{\det(\Hess_{\optim{\mPoint}} \dummyFunctionA)}\frac{\DimInd+1 +
\epsilon_{\revA{1}}}{\DimInd+2}(\radiusA_\sumTotA^{\revA{+}})^2 + \mathcal{O}((\radiusA_\sumTotA^{\revA{+}})^3),
\end{multline}
and for large $\sumTotA$
{\allowdisplaybreaks
\begin{multline}
\label{thm:eq-step-second-term-ellipsoid-integration}
\revA{\frac{1}{|\ellipsoid_{\radiusC_\sumTotA^{+}}^{\dummyBilForm}(\optim{\mPoint})\cap \rotsSampling_\sumTotA|}}
\!\!
\sum_{\mPointB \in\ellipsoid_{\revA{\radiusC_\sumTotA^{+}}}^{\dummyBilForm}(\optim{\mPoint})\cap \rotsSampling_\sumTotA}
\!\!\! 
\bigl(\dummyFunctionA(\mPointB)- \dummyFunctionA(\optim{\mPoint})\bigr) 
\\
= \frac{1}{|\ellipsoid_{\revA{\radiusC_\sumTotA^{+}}}^{\dummyBilForm}(\optim{\mPoint})\cap \rotsSampling_\sumTotA|}
\!\!
\sum_{\mPointB \in\ellipsoid_{\revA{\radiusC_\sumTotA^{+}}}^{\dummyBilForm}(\optim{\mPoint})\cap \rotsSampling_\sumTotA}
\!\!\!
\Bigl(\frac{1}{2}\Hess_{\optim{\mPoint}} \dummyFunctionA   
  \bigl(\log_{\optim{\mPoint}}(\mPointB), 
        \log_{\optim{\mPoint}}(\mPointB)
  \bigr) 
  + \mathcal{O}((\revA{\radiusC_\sumTotA^{+}})^3)
\Bigr) 
\\
=  \frac{|\rotsSampling_\sumTotA|}{|\ellipsoid_{\revA{\radiusC_\sumTotA^{+}}}^{\dummyBilForm}(\optim{\mPoint})\cap \rotsSampling_\sumTotA|}
\!\!
\sum_{\mPointB \in\ellipsoid_{\revA{\radiusC_\sumTotA^{+}}}^{\dummyBilForm} (\optim{\mPoint})\cap \rotsSampling_\sumTotA}
\!\!\!
\frac{\Bigl(\frac{1}{2}\Hess_{\optim{\mPoint}} \dummyFunctionA   
  \bigl(\log_{\optim{\mPoint}}(\mPointB), 
        \log_{\optim{\mPoint}}(\mPointB)
  \bigr) 
  + \mathcal{O}((\revA{\radiusC_\sumTotA^{+}})^3)
\Bigr) }{|\rotsSampling_\sumTotA|}
\\
\overset{\cref{eq:eta-local-low-discrepancy-sequence-1}}{=} \frac{\operatorname{vol}(\manifold)}{
    \operatorname{vol}\bigl(
      \ellipsoid_{\revA{\radiusC_\sumTotA^{+}}}^{\dummyBilForm}(\optim{\mPoint})
    \bigr) + \revA{o}((\revA{\radiusC_\sumTotA^{+}})^{\DimInd})}
\!\!
\sum_{\mPointB \in\ellipsoid_{\revA{\radiusC_\sumTotA^{+}}}^{\dummyBilForm} (\optim{\mPoint})\cap \rotsSampling_\sumTotA}
\!\!\!
\frac{\Bigl(\frac{1}{2}\Hess_{\optim{\mPoint}} \dummyFunctionA   
  \bigl(\log_{\optim{\mPoint}}(\mPointB), 
        \log_{\optim{\mPoint}}(\mPointB)
  \bigr) 
  + \mathcal{O}((\revA{\radiusC_\sumTotA^{+}})^3)
\Bigr) }{|\rotsSampling_\sumTotA|}
\\[0.5em]
\overset{\cref{eq:eta-local-low-discrepancy-sequence-2}}{=} \frac{\int_{\ellipsoid_{\revA{\radiusC_\sumTotA^{+}}}^{\dummyBilForm}(\optim{\mPoint})}
\Bigl[ 
\frac{1}{2}\Hess_{\optim{\mPoint}} 
  \dummyFunctionA\bigl(
    \log_{\optim{\mPoint}}(\mPointB),
    \log_{\optim{\mPoint}}(\mPointB)
  \bigr) 
  + \mathcal{O}((\revA{\radiusC_\sumTotA^{+}})^3)
\Bigr] \mathrm{d}\mPointB + \revA{o((\revA{\radiusC_\sumTotA^{+}})^{\DimInd+2})}}{
    \operatorname{vol}\bigl(
      \ellipsoid_{\revA{\radiusC_\sumTotA^{+}}}^{\dummyBilForm}(\optim{\mPoint})
    \bigr)
    + \revA{o}((\revA{\radiusC_\sumTotA^{+}})^{\revA{\DimInd}})}
\\[0.5em]
\overset{\cref{eq:ellips-integral-2}}{=} \frac{1}{2}
\dfrac{\sqrt[\DimInd]{\det(\Hess_{\optim{\mPoint}} \dummyFunctionA)} \frac{\ballVol_{\DimInd}\DimInd}{\DimInd+2}(\revA{\radiusC_\sumTotA^{+}})^{\DimInd+2}  + \mathcal{O}((\revA{\radiusC_\sumTotA^{+}})^{\DimInd+3}) + \revA{o((\revA{\radiusC_\sumTotA^{+}})^{\DimInd+2})}}{\ballVol_{\DimInd}(\revA{\radiusC_\sumTotA^{+}})^{\DimInd} + \revA{o}((\revA{\radiusC_\sumTotA^{+}})^{\revA{\DimInd}})}
\\
= \frac{1}{2}\sqrt[\DimInd]{\det(\Hess_{\optim{\mPoint}} \dummyFunctionA)} \frac{\DimInd}{\DimInd+2}\revA{\frac{\DimInd+3 - \epsilon_2}{\DimInd+2}}(\radiusA_\sumTotA^{\revA{+}})^{2}  + \revA{o}(\revA{(\radiusA_\sumTotA^{\revA{+}})^{2}}).
\end{multline}}
   
Bringing together the above results yields the validity of \cref{eq:thm-J-bound-step-2-opt-ineq-eps}
\begin{multline}
\revA{\liminf_{\sumTotA\to \infty}} \frac{\dummyFunctionA(\mPoint)-\frac{1}{\dummyIntegerB_\sumTotA^{\revA{+}}} \sum_{\sumIndB=1}^{\dummyIntegerB_\sumTotA^{\revA{+}}} \dummyFunctionA(\mPoint^\sumTotA_{(\sumIndB)})}{(\radiusA_\sumTotA^{\revA{+}})^{2}} 
\\
\revA{\geq} \lim_{\sumTotA\to \infty} \frac{1}{2}\frac{ \sqrt[\DimInd]{\det(\Hess_{\optim{\mPoint}} \dummyFunctionA)}\frac{\DimInd+1+\epsilon_{\revA{1}}}{\DimInd+2}(\radiusA_\sumTotA^{\revA{+}})^2 - \sqrt[\DimInd]{\det(\Hess_{\optim{\mPoint}} \dummyFunctionA)} \frac{\DimInd}{\DimInd+2}\revA{\frac{\DimInd+3 - \epsilon_2}{\DimInd+2}}(\radiusA_\sumTotA^{\revA{+}})^{2} + \revA{o}((\radiusA_\sumTotA^{\revA{+}})^{\revA{2}})}{(\radiusA_\sumTotA^{\revA{+}})^{2}} 
\\
\revA{\overset{\epsilon_1, \epsilon_2 \in (\frac{1}{2\DimInd + 2}, 1)}{\geq}} (1+\epsilon^{\revA{+}})\frac{\sqrt[\DimInd]{\det(\Hess_{\optim{\mPoint}} \dummyFunctionA)}}{2(\DimInd+2)},
\label{eq:thm-numpoints-combined-results-lb}
\end{multline}
and for large, but finite $\sumTotA$ still too many points must satisfy \cref{eq:thm-J-bound-step-2-opt-ineq} and hence \cref{eq:thm-num-nonzero-contradiction} follows, which gives our contradiction. We conclude that $\dummyIntegerB_\sumTotA < \dummyIntegerB_\sumTotA^{\revA{+}}$ for large enough $\sumTotA$. Subsequently, we can drop the ceil operation, which gives the bound in \cref{eq:thm-J-bounds}.

(i-b) 
Our claim that the points with non-zero coefficients lie in the ellipsoid -- as stated in the theorem -- also follows by realising that the upper bound from (i-a) gives us that the points with non-zero coefficients fit in $\ellipsoid_{\radiusA_\sumTotA}^{\dummyBilForm}(\optim{\mPoint})$ and by a similar level-set argument as above it must be that these points indeed live in this ellipsoid for large enough $\sumTotA$. 

\vspace{0.5cm}

(ii) For the lower bound in \cref{eq:thm-J-bounds} we will show that $\|\optim{\rotsDensCoef}_\sumTotA\|_0 > \dummyIntegerB_\sumTotA^{\revA{-}}$, where 
\begin{equation}
    \dummyIntegerB_\sumTotA^{\revA{-}} := \lfloor 3^{-\frac{\DimInd}{\DimInd+2}}\dummyIntegerB_0 |\rotsSampling_\sumTotA|^{\frac{2 - \DimInd\exponentA}{\DimInd+2}} \rfloor.
    \label{eq:thm-non-zero-Nm''}
\end{equation}

In particular -- and similarly to step (i-a) --, we will assume that $\|\optim{\rotsDensCoef}_\sumTotA\|_0 \leq \dummyIntegerB_\sumTotA^{\revA{-}}$ but show that, for large $\sumTotA$,
\begin{equation}
    \label{eq:thm-num-nonzero-contradiction-lb}
    \biggl\|\biggl(- \frac{|\rotsSampling_\sumTotA|^{\exponentA}}{\regParam}\dummyFunctionA_{\rotsSampling_\sumTotA} -\frac{1}{\dummyIntegerB_\sumTotA^{\revA{-}}} \Bigl(\sum_{j=1}^{\dummyIntegerB_\sumTotA^{\revA{-}}} - \frac{|\rotsSampling_\sumTotA|^{\exponentA}}{\regParam} \dummyFunctionA(\mPoint^\sumTotA_{(\sumIndB)})\Bigr) + \frac{1}{\dummyIntegerB_\sumTotA^{\revA{-}}}\biggr)_{+}\biggr\|_0 > \dummyIntegerB_\sumTotA^{\revA{-}}
\end{equation}
holds. Then, also
\begin{equation}
    \biggl\|\biggl(- \frac{|\rotsSampling_\sumTotA|^{\exponentA}}{\regParam}\dummyFunctionA_{\rotsSampling_\sumTotA} -\frac{1}{\dummyIntegerB_\sumTotA^{\revA{-}}} \Bigl(\sum_{j=1}^{\dummyIntegerB_\sumTotA^{\revA{-}}} - \frac{|\rotsSampling_\sumTotA|^{\exponentA}}{\regParam} \dummyFunctionA(\mPoint^\sumTotA_{(\sumIndB)})\Bigr) + \frac{1}{\dummyIntegerB_\sumTotA^{\revA{-}}}\biggr)_{+}\biggr\|_0 > \|\optim{\rotsDensCoef}_\sumTotA\|_0
\end{equation}
holds by assumption. However, the latter is equivalent to $\|\optim{\rotsDensCoef}_\sumTotA\|_0 > \dummyIntegerB_\sumTotA^{\revA{-}}$ by \cref{rem:monotonicity-of-projection-correction}, which yields our contradiction. The remainder of (ii) is showing \cref{eq:thm-num-nonzero-contradiction-lb}.

\vspace{0.5cm}

If for any point $\mPoint \in \rotsSampling_\sumTotA$
\begin{equation}
    \biggl(- \frac{|\rotsSampling_\sumTotA|^{\exponentA}}{\regParam}\dummyFunctionA(\mPoint) -\frac{1}{\dummyIntegerB_\sumTotA^{\revA{-}}} \Bigl(\sum_{j=1}^{\dummyIntegerB_\sumTotA^{\revA{-}}} - \frac{|\rotsSampling_\sumTotA|^{\exponentA}}{\regParam} \dummyFunctionA(\mPoint^\sumTotA_{(\sumIndB)})\Bigr) + \frac{1}{\dummyIntegerB_\sumTotA^{\revA{-}}}\biggr)_{+} > 0
\end{equation}
holds, then
\begin{equation}
   \Bigl(- \frac{|\rotsSampling_\sumTotA|^{\exponentA}}{\regParam}\dummyFunctionA(\mPoint)\Bigr) - \frac{1}{\dummyIntegerB_\sumTotA^{\revA{-}}} \Bigl(\sum_{\sumIndB=1}^{\dummyIntegerB_\sumTotA^{\revA{-}}} - \frac{|\rotsSampling_\sumTotA|^{\exponentA}}{\regParam} \dummyFunctionA(\mPoint^\sumTotA_{(\sumIndB)}) \Bigr) + \frac{1}{\dummyIntegerB_\sumTotA^{\revA{-}} } > 0 
\end{equation}
or equivalently
\begin{equation}
  \dummyFunctionA(\mPoint) - \frac{1}{\dummyIntegerB_\sumTotA^{\revA{-}}} \sum_{\sumIndB=1}^{\dummyIntegerB_\sumTotA^{\revA{-}}} \dummyFunctionA(\mPoint^\sumTotA_{(\sumIndB)})< \frac{\regParam}{\dummyIntegerB_\sumTotA^{\revA{-}} |\rotsSampling_\sumTotA|^{\exponentA}}.
  \label{eq:Jm''-ineq}
\end{equation}

Similarly to step (i), we will show that too many points -- that is more than $\dummyIntegerB_\sumTotA^{\revA{-}}$ -- satisfy an even stronger inequality. For that, first define $(\radiusA^{\revA{-}}_\sumTotA)_{\sumTotA=1}^\infty \subset \Real_+$ with
\begin{equation}
\label{eq:thm-non-zero-coeffs-radiuspp-def}
\radiusA_\sumTotA^{\revA{-}} := \sqrt[\DimInd]{\revA{3^{-\frac{\DimInd}{\DimInd+2}}\dummyIntegerB_0}\frac{\operatorname{vol}(\manifold)}{\ballVol_\DimInd} } \revA{|\rotsSampling_\sumTotA|^{-\frac{1+\exponentA}{\DimInd+2}}},
\end{equation}
as in \cref{lem:numpoints-in-ellipsoid} and -- analogously to \cref{eq:JR-ineq} -- note that




\begin{multline}
\frac{\regParam}{\dummyIntegerB_\sumTotA^{\revA{-}} |\rotsSampling_\sumTotA|^{\exponentA}} 
\overset{\cref{eq:thm-non-zero-Nm''}}{\geq} \frac{\regParam}{3^{-\frac{\DimInd}{\DimInd+2}}\dummyIntegerB_0|\rotsSampling_\sumTotA|^{\frac{2 - \DimInd\exponentA}{\DimInd+2}} |\rotsSampling_\sumTotA|^\exponentA} 
= \frac{\regParam}{3^{-\frac{\DimInd}{\DimInd+2}}\dummyIntegerB_0|\rotsSampling_\sumTotA|^{\frac{2 + 2\exponentA}{\DimInd+2}}}
\\
= \frac{\regParam}{3^{-\frac{\DimInd}{\DimInd+2}}\dummyIntegerB_0}\Bigl(\frac{\ballVol_\DimInd}{3^{-\frac{\DimInd}{\DimInd+2}}\dummyIntegerB_0\operatorname{vol}(\manifold)}\Bigr)^{\frac{2}{\DimInd}}\Bigl(\frac{3^{-\frac{\DimInd}{\DimInd+2}}\dummyIntegerB_0\operatorname{vol}(\manifold)}{\ballVol_\DimInd}\Bigr)^{\frac{2}{\DimInd}} |\rotsSampling_\sumTotA|^{-\frac{2 + 2\exponentA}{\DimInd+2}}
\\
\overset{\cref{eq:thm-non-zero-coeffs-radiuspp-def}}{\revA{=}} \frac{\regParam}{3^{-\frac{\DimInd}{\DimInd+2}}\dummyIntegerB_0}\Bigl(\frac{\ballVol_\DimInd}{3^{-\frac{\DimInd}{\DimInd+2}}\dummyIntegerB_0\operatorname{vol}(\manifold)}\Bigr)^{\frac{2}{\DimInd}}(\radiusA_\sumTotA^{\revA{-}})^{2}
\overset{\cref{eq:thm-new-num-non-zero-coeffs-gamma}}{=} 3\frac{\sqrt[\DimInd]{\det(\Hess_{\optim{\mPoint}} \dummyFunctionA)}}{2(\DimInd+2)}(\radiusA_\sumTotA^{\revA{-}})^{2}.
\end{multline}
Hence, for $\mPoint\in \rotsSampling_\sumTotA$ \cref{eq:thm-num-nonzero-contradiction-lb} holds if
\begin{equation}
  \dummyFunctionA(\mPoint) - \frac{1}{\dummyIntegerB_\sumTotA^{\revA{-}}} \sum_{\sumIndB=1}^{\dummyIntegerB_\sumTotA^{\revA{-}}} \dummyFunctionA(\mPoint^\sumTotA_{(\sumIndB)})< 3 \frac{\sqrt[\DimInd]{\det(\Hess_{\optim{\mPoint}} \dummyFunctionA)}}{2(\DimInd+2)}(\radiusA_\sumTotA^{\revA{-}})^{2} .
\label{eq:thm-J''-bound-step-2-opt-ineq}
\end{equation}
Next, \revA{note that $0< \frac{\DimInd(\DimInd + 1+ \epsilon_4) - (\DimInd +2 )(\DimInd - \epsilon_3)}{\DimInd + 2} < 1$ for $\epsilon_3, \epsilon_4 \in (\frac{1}{2\DimInd + 2}, 1)$ and} fix \revA{the parameter} $0< \epsilon^{\revA{-}}< \revA{\frac{\DimInd(\DimInd + 1+ \epsilon_4) - (\DimInd +2 )(\DimInd - \epsilon_3)}{\DimInd + 2}} <1$. Then, by a truncation argument it suffices to consider the asymptotic validity of the yet stronger inequality
\begin{equation}
  \revA{\limsup_{\sumTotA\to \infty}} \frac{\dummyFunctionA(\mPoint) - \frac{1}{\dummyIntegerB_\sumTotA^{\revA{-}}} \sum_{\sumIndB=1}^{\dummyIntegerB_\sumTotA^{\revA{-}}} \dummyFunctionA(\mPoint^\sumTotA_{(\sumIndB)})}{(\radiusA_\sumTotA^{\revA{-}})^{2}}\leq (3-\epsilon^{\revA{-}}) \frac{\sqrt[\DimInd]{\det(\Hess_{\optim{\mPoint}} \dummyFunctionA)}}{2(\DimInd+2)},
\label{eq:thm-J''-bound-step-2-opt-ineq-eps}
\end{equation}
Then for large $\sumTotA$ it still holds that too many points satisfy \cref{eq:thm-J''-bound-step-2-opt-ineq}  -- and through \cref{eq:Jm''-ineq} also \cref{eq:thm-num-nonzero-contradiction-lb} --, which yields our claim. The remainder of step (ii) is devoted to showing that asymptotically too many $\mPoint\in \rotsSampling_\sumTotA$ satisfy \cref{eq:thm-J''-bound-step-2-opt-ineq-eps}.

\vspace{0.5cm}

Consider the ellipsoid $\ellipsoid_{\revA{\radiusC}_\sumTotA^{\revA{-}}}^\dummyBilForm(\optim{\mPoint})$, again with $\dummyBilForm:= \Hess_{\optim{\mPoint}} \dummyFunctionA$, but with $\revA{\radiusC}_\sumTotA^{\revA{-}}:= \sqrt{\frac{\DimInd+3 - \epsilon_{\revA{3}}}{\DimInd+2}}\radiusA_\sumTotA^{\revA{-}}$. We will show that all points $\mPoint\in\rotsSampling_\sumTotA$ outside of $\ellipsoid_{\revA{\radiusC}_\sumTotA^{\revA{-}}}^\dummyBilForm(\optim{\mPoint})$ asymptotically satisfy \cref{eq:thm-J''-bound-step-2-opt-ineq-eps}. Then,
\begin{equation}
    \dummyIntegerB_\sumTotA^{\revA{-}} \revA{\leq  \lceil 3^{-\frac{\DimInd}{\DimInd+2}}\dummyIntegerB_0 |\rotsSampling_\sumTotA|^{\frac{2 - \DimInd\exponentA}{\DimInd+2}} \rceil } < |\ellipsoid_{\revA{\radiusC}_\sumTotA^{\revA{-}}}^{\dummyBilForm}(\optim{\mPoint}) \cap \rotsSampling_\sumTotA|
    \label{eq:Jm''discrepency-ER}
\end{equation}
holds by \cref{eq:lem-J-as-points-in-ellipsoid-<} in \cref{lem:numpoints-in-ellipsoid}, which yields the claim.

For showing that for large $\sumTotA$ all points outside of the ellipsoid $\ellipsoid_{\revA{\radiusC}_\sumTotA^{\revA{-}}}^\dummyBilForm(\optim{\mPoint})$ satisfy \cref{eq:thm-J''-bound-step-2-opt-ineq-eps}, it suffices to show to show that the inequality \cref{eq:thm-J''-bound-step-2-opt-ineq-eps} holds for all points on the boundary $\partial \ellipsoid_{\revA{\radiusC}_\sumTotA^{\revA{-}}}^\dummyBilForm(\optim{\mPoint})$. Indeed, by a similar line of reasoning as in step (i-a) the points with the lowest $\dummyFunctionA$-value will live in this ellipsoid.



So, let $\mPoint\in \partial \ellipsoid_{\revA{\radiusC}_\sumTotA^{\revA{-}}}^\dummyFunctionA(\optim{\mPoint})$, i.e., we have
\begin{equation}
\Hess_{\optim{\mPoint}} \dummyFunctionA \bigl(\log_{\optim{\mPoint}}(\mPoint), \log_{\optim{\mPoint}}(\mPoint)\bigr) 
  \overset{\cref{eq:Ellipsoid}}{=} \sqrt[\DimInd]{\det(\Hess_{\optim{\mPoint}} \dummyFunctionA)}\,(\revA{\radiusC}_\sumTotA^{\revA{-}})^2
  = \sqrt[\DimInd]{\det(\Hess_{\optim{\mPoint}} \dummyFunctionA)}\frac{\DimInd+3 - \epsilon_{\revA{3}}}{\DimInd+2} (\radiusA_\sumTotA^{\revA{-}})^2.
\end{equation}
Then, we have
\begin{multline}
\revA{\limsup_{\sumTotA\to \infty}} \frac{\dummyFunctionA(\mPoint)-\frac{1}{\dummyIntegerB_\sumTotA^{\revA{-}}} \sum_{\sumIndB=1}^{\dummyIntegerB_\sumTotA^{\revA{-}}} \dummyFunctionA(\mPoint^\sumTotA_{(\sumIndB)})}{(\radiusA_\sumTotA^{\revA{-}})^{2}} 
\\
\revA{\leq} \lim_{\sumTotA\to \infty} \frac{1}{2}\frac{ \sqrt[\DimInd]{\det(\Hess_{\optim{\mPoint}} \dummyFunctionA)}\frac{\DimInd+3-\epsilon_\revA{3}}{\DimInd+2}(\radiusA_\sumTotA^{\revA{-}})^2 - \sqrt[\DimInd]{\det(\Hess_{\optim{\mPoint}} \dummyFunctionA)} \frac{\DimInd}{\DimInd+2}\revA{\frac{\DimInd+1+\epsilon_{\revA{4}}}{\DimInd+2}}(\radiusA_\sumTotA^{\revA{-}})^{2} + \revA{o}((\radiusA_\sumTotA^{\revA{-}})^{\revA{2}})}{(\radiusA_\sumTotA^{\revA{-}})^{2}} 
\\
\revA{\overset{\epsilon_3, \epsilon_4 \in (\frac{1}{2\DimInd + 2}, 1)}{\leq}} (3-\epsilon^{\revA{-}})\frac{\sqrt[\DimInd]{\det(\Hess_{\optim{\mPoint}} \dummyFunctionA)}}{2(\DimInd+2)},
\label{eq:thm-numpoints-combined-results-ub}
\end{multline}
where the second term in the second line can be obtained in a similar fashion as done in \cref{thm:eq-step-second-term-ellipsoid-integration}\revA{, but now by passing to a smaller ellipsoid with radius $\radiusB_\sumTotA^{\revA{-}}:= \sqrt{\frac{\DimInd+1 + \epsilon_{\revA{4}}}{\DimInd+2}}\radiusA_\sumTotA^{\revA{-}}$}.
Finally, for large, but finite $\sumTotA$ still too many points must satisfy \cref{eq:thm-J''-bound-step-2-opt-ineq} and hence \cref{eq:thm-num-nonzero-contradiction-lb} follows, which gives our contradiction. We conclude that $\dummyIntegerB_\sumTotA > \dummyIntegerB_\sumTotA^{\revA{-}}$ for large enough $\sumTotA$. Subsequently, we can drop the floor operation, which gives the bound in \cref{eq:thm-J-bounds}.


\end{proof}

\subsection{Well-posedness and consistency of the Riemannian barycentre}
\label{sec:lifting-wellposed}

For proving \cref{thm:distance-bound-post-lifting} we will argue that the error is fully determined by the choice of measure and moreover, that is sufficient to bound the error made by the measure.
To make this more concrete, the lemma below expresses the distance between any point $\optim{\mPoint}$, e.g., a minimiser of \revA{a} function on $\manifold$, and a solution $\optim{\mPoint}_\measure$ to the lifted problem \cref{eq:lifted-projection} in terms of the Wasserstein distance between a point mass at $\optim{\mPoint}$ and the corresponding measure solving \cref{eq:lifted-function}.


\begin{lemma}\label{Lemma:WasserSteinEst} 
Let $\manifold$ be a connected Riemannian manifold, let $\optim{\mPoint}\in \manifold$ and let $\measure \in \probSpace_2(\manifold)$ be a measure with support on a geodesically connected subset of $\manifold$. Then,
\begin{equation}
\distance_{\manifold}(\optim{\mPoint}_\measure,\optim{\mPoint}) 
\leq 
2\wasserstein_2(\measure,\dirac_{\optim{\mPoint}})
\quad\text{for all $\optim{\mPoint}_\measure \in \mPoints(\measure) \subset \manifold$}
\end{equation}
where the (possibly set-valued) solution operator mapping $P \colon \probSpace_2(\manifold)\rightrightarrows\manifold$ is defined as in~\cref{eq:lifted-projection}.
\end{lemma}
\begin{proof}
First, note that for all $\mPoint\in \manifold$
\begin{equation}
    \wasserstein_2(\measure, \dirac_{\mPoint})^2 = \inf_{\pi \in \Gamma(\measure, \dirac_{\mPoint})} \int_{\manifold \times \manifold} \distance_{\manifold}(\mPointB, \mPointB')^{2} \mathrm{d}\pi(\mPointB, \mPointB') = \int_{\manifold} \distance_{\manifold}(\mPointB, \mPoint)^2 \mathrm{d}\measure(\mPointB),
\end{equation}
as every coupling $\pi$ with marginals $\measure$ and $\dirac_{\mPoint}$ must satisfy $\pi=\measure \otimes \dirac_{\mPoint}$.

Next, $\optim{\mPoint}_\measure$ exists since $\measure$ has support on a geodesically connected subset and we have $\wasserstein_2(\measure, \dirac_{\optim{\mPoint}_\measure}) \leq \wasserstein_2(\measure, \dirac_{\optim{\mPoint}})$, since 
\begin{align*}
\wasserstein_2(\measure, \dirac_{\optim{\mPoint}_\measure})^2 
&=  \int_{\manifold} 
\distance_{\manifold}(\mPointB, \optim{\mPoint}_\measure)^2 \mathrm{d}\measure( \mPointB)  
= \inf_{\mPoint \in \manifold}  
\int_{\manifold} 
\distance_{\manifold}(\mPointB, \mPoint)^2 \mathrm{d}\measure(\mPointB) 
\\
& \leq \int_{\manifold} \distance_{\manifold}(\mPointB, \optim{\mPoint})^2 \mathrm{d}\measure(\mPointB)
= \wasserstein_2(\measure, \dirac_{\optim{\mPoint}})^2.
\end{align*}
The claim now follows from fact that $W_2$ is satisfies the triangle inequality:
\begin{equation}
\distance_{\manifold}(\optim{\mPoint}_\measure,\optim{\mPoint}) 
= \wasserstein_2(
\dirac_{\optim{\mPoint}_\measure}, \dirac_{\optim{\mPoint}}
) 
\leq \wasserstein_2(\measure, \dirac_{\optim{\mPoint}_\measure}) + \wasserstein_2(\measure, \dirac_{\optim{\mPoint}}) \leq 2 \wasserstein_2(\measure, \dirac_{\optim{\mPoint}}).
\end{equation}
\end{proof}
\Cref{Lemma:WasserSteinEst} shows that an analysis of the error rate that covers not only the convexification and projection errors, but also discretisation and relaxation errors can be based on bounds for the Wasserstein distance. Indeed, if the non-zero coefficients of $\rotsDensCoef$ concentrate on a geodesically connected subset of $\manifold$ and $\wasserstein_2(\sum_{\mPoint\in \rotsSampling}
    \optim{\rotsDensCoef}_{\mPoint}\dirac_{\mPoint},\dirac_{\optim{\mPoint}})$ goes to $0$ for large enough sampling sets, we obtain convergence of the result of the lifting approach $\optim{\mPoint}_\rotsSampling$ to the solution of the original problem, say $\optim{\mPoint}$, as well.

        \begin{proof}[Proof of \cref{thm:distance-bound-post-lifting}]
\revA{For the given $\epsilon_1, \epsilon_2, \epsilon_3, \epsilon_4 \in (\frac{\DimInd}{2\DimInd + 2},1)$, $1/(\DimInd+1) < \exponentA < 2/\DimInd$ and $\dummyIntegerB_0 \geq 3^{-\frac{\DimInd}{\DimInd +2}}$ f}ix a local low-discrepancy sequence $(\rotsSampling_\sumTotA)_{\sumTotA=1}^\infty \subset S(\manifold)$ \revA{fulfilling the conditions in the statement}. Define $(\optim{\rotsDensCoef}_\sumTotA)_{\sumTotA=1}^\infty \subset \ell_0$ as the the zero-padded solution operator of \cref{eq:lifted-function-discr}, so it is given as
\begin{equation}
    \optim{\rotsDensCoef}_\sumTotA:= \biggl( \Pi_{\simplex^{|\rotsSampling_\sumTotA|}}\Bigl(- \frac{|\rotsSampling_\sumTotA|^{\exponentA}}{\regParam}\dummyFunctionA_{\rotsSampling_\sumTotA} \Bigr), 0,\ldots\biggr).
\end{equation}

First, from the locality statement of \cref{thm:new-num-non-zero-coeffs} -- that is the points with non-zero coefficients being concentrated on an ellipsoid -- existence of the barycentre generated by $\optim{\rotsDensCoef}_\sumTotA$ for $\sumTotA\geq \sumTotA'$ holds on the submanifold $\ball_{\radiusB_{\sumTotA'}}(\optim{\mPoint}) \subset \manifold$ that is the (geodesically convex) geodesic ball with radius $\radiusB_{\sumTotA'}$ such that it is containing all the ellipsoids containing the support of $\optim{\rotsDensCoef}_\sumTotA$ for $\sumTotA\geq \sumTotA'$ by \cite[Proposition~1]{le2017existence}. With the measures concentrated in $\ball_{\radiusB_{\sumTotA'}}(\optim{\mPoint})$ the uniqueness claim of the barycentre holds by the famous result by Kendall \cite[Theorem~7.3]{kendall1990probability}. This shows the first claim \cref{eq:thm-barycentre}.



Consequently, for the bound in the second claim \cref{eq:cor-bound-distance}, $\probSpace_2(\ball_{\radiusB_{\sumTotA'}}(\optim{\mPoint}))$ is well-defined and from \cref{Lemma:WasserSteinEst} we obtain that it is sufficient to show
\begin{equation}
\label{eq:wasserstein-bounddd}
    \wasserstein_2\Bigl(
  \sum_{\mPoint \in \rotsSampling_\sumTotA}(\optim{\rotsDensCoef}_\sumTotA)_{\mPoint} \delta_{\mPoint},\dirac_{\optim{\mPoint}}
\Bigr) \leq \sqrt{\frac{\sqrt[\DimInd]{\det(\Hess_{\optim{\mPoint}} \dummyFunctionA)}}{\eigenValue_{\min}(\Hess_{\optim{\mPoint}} \dummyFunctionA)}} \sqrt[\DimInd]{
\dummyIntegerB_0 \frac{\operatorname{vol}(\manifold)}{\ballVol_\DimInd}} |\rotsSampling_\sumTotA|^{-\frac{1+\exponentA}{\DimInd+2}}.
\end{equation}
Indeed, if the above is established, the additional factor of 2 from the lemma gives the desired result.

\vspace{0.5cm}

So it remains to show \cref{eq:wasserstein-bounddd}. For that, note that for large $\sumTotA$, there must be less than $\dummyIntegerB_0 |\rotsSampling_\sumTotA|^{\frac{2 - \DimInd\exponentA}{\DimInd+2}}$ points with non-zero coefficients by \cref{thm:new-num-non-zero-coeffs}. Additionally, these points lie in the ellipsoid $\ellipsoid_{\radiusA_\sumTotA}^{\dummyBilForm}(\optim{\mPoint})$ with $\dummyBilForm:= \Hess_{\optim{\mPoint}} \dummyFunctionA$, whose size is defined by the radius $\radiusA_\sumTotA\in \Real_+$, defined in \cref{eq:EllipsRadii}, i.e.,
\begin{equation}
\radiusA_\sumTotA := \sqrt[\DimInd]{\dummyIntegerB_0\frac{\operatorname{vol}(\manifold)}{\ballVol_\DimInd}}|\rotsSampling_\sumTotA|^{-\frac{1+\exponentA}{\DimInd+2}}.
\end{equation}


Then, taking into account that
\begin{equation}
\label{eq:ellipsoid-size-eigenvalues}
    \distance_\manifold(\mPoint, \optim{\mPoint})^2 \leq \frac{\sqrt[\DimInd]{\det(\Hess_{\optim{\mPoint}} \dummyFunctionA)}}{\eigenValue_{\min}(\Hess_{\optim{\mPoint}} \dummyFunctionA)}(\radiusA_\sumTotA)^2\quad \text{for $\mPoint\in \ellipsoid_{\radiusA_\sumTotA}^{\dummyBilForm}(\optim{\mPoint})$,}
\end{equation}
the bound on the Wasserstein distance follows directly, as for large $\sumTotA$
\begin{multline}
\wasserstein_2\Bigl(
  \sum_{\mPoint \in \rotsSampling_\sumTotA}(\optim{\rotsDensCoef}_\sumTotA)_{\mPoint} \delta_{\mPoint},\dirac_{\optim{\mPoint}}
\Bigr)^2 
= \sum_{\mPoint \in \rotsSampling_\sumTotA} (\optim{\rotsDensCoef}_\sumTotA)_{\mPoint}\; \distance_{\manifold}(\mPoint,\optim{\mPoint})^2
\overset{\cref{thm:new-num-non-zero-coeffs}}{\leq} \sum_{\mPoint \in \ellipsoid_{\radiusA_\sumTotA}^\dummyBilForm(\optim{\mPoint}) \cap \rotsSampling_\sumTotA} (\optim{\rotsDensCoef}_\sumTotA)_{\mPoint}\; \distance_{\manifold}(\mPoint,\optim{\mPoint})^2
\\
 \overset{\cref{eq:ellipsoid-size-eigenvalues}}{\leq} \sum_{\mPoint \in \ellipsoid_{\radiusA_\sumTotA}^\dummyBilForm(\optim{\mPoint}) \cap \rotsSampling} (\optim{\rotsDensCoef}_\sumTotA)_{\mPoint} \frac{\sqrt[\DimInd]{\det(\Hess_{\optim{\mPoint}} \dummyFunctionA)}}{\eigenValue_{\min}(\Hess_{\optim{\mPoint}} \dummyFunctionA)}(\radiusA_\sumTotA)^2 
 \leq \frac{\sqrt[\DimInd]{\det(\Hess_{\optim{\mPoint}} \dummyFunctionA)}}{\eigenValue_{\min}(\Hess_{\optim{\mPoint}} \dummyFunctionA)}(\radiusA_\sumTotA)^2,
\end{multline}
where the last equality follows from the fact that $\sum_{\mPoint\in \rotsSampling_\sumTotA} (\optim{\rotsDensCoef}_\sumTotA)_\mPoint = 1$.
\end{proof}

        \subsection{Approximating the regularisation parameter}
        \label{sec:approx-regu}
        Finally, it remains to show \cref{prop:choice-regulariser-lifting}. The proof relies heavily on the ideas developed in the proof of \cref{thm:new-num-non-zero-coeffs}.

        \begin{proof}[Proof of \cref{prop:choice-regulariser-lifting}]


\revA{For the given $\epsilon_1, \epsilon_2 \in (\frac{\DimInd}{2\DimInd + 2},1)$, $1/(\DimInd+1) < \exponentA < 2/\DimInd$ and $\dummyIntegerB_0 \geq 3^{-\frac{\DimInd}{\DimInd +2}}$ f}ix a local low-discrepancy sequence $(\rotsSampling_\sumTotA)_{\sumTotA=1}^\infty \subset S(\manifold)$ \revA{fulfilling the conditions in the statement} and define $(\dummyIntegerB_\sumTotA)_{\sumTotA=1}^\infty\subset\Natural$ by $\dummyIntegerB_\sumTotA := \lfloor \dummyIntegerB_0 |\rotsSampling_\sumTotA|^{\frac{2 - \DimInd\exponentA}{\DimInd+2}} \rfloor$ and $(\radiusA_\sumTotA)_{\sumTotA=1}^\infty\subset\Natural$ by $\radiusA_\sumTotA := \revA{\sqrt[\DimInd]{\dummyIntegerB_0\frac{\operatorname{vol}(\manifold)}{\ballVol_\DimInd}}|\rotsSampling_\sumTotA|^{-\frac{1+\exponentA}{\DimInd+2}}}$. Additionally, for any $\sumTotA \in \Natural$ define the points $\mPoint^\sumTotA_{(1)}, \ldots,\mPoint^\sumTotA_{(|\rotsSampling_\sumTotA|)}$ by labeling the elements of $\rotsSampling_\sumTotA$ such that  $\dummyFunctionA(\mPoint^\sumTotA_{(1)})\leq \ldots\leq \dummyFunctionA(\mPoint^\sumTotA_{(|\rotsSampling_\sumTotA|)})$.

Let \revA{$0< \epsilon^{\revA{-}}< \revA{\frac{(\DimInd + 2)(\DimInd + \epsilon_1) - \DimInd(\DimInd + 3 - \epsilon_2)}{\DimInd + 2}}< 1$, $0< \epsilon^{+} < \frac{\DimInd(\DimInd + 1+ \epsilon_1) - (\DimInd +2 )(\DimInd - \epsilon_2)}{\DimInd + 2}<1$} and let $\mPoint\in \partial \ellipsoid_{\revA{\radiusA}_\sumTotA}^\dummyBilForm(\optim{\mPoint})$ with $\dummyBilForm:= \Hess_{\optim{\mPoint}}\dummyFunctionA$. Then, analogously to the argument in \cref{thm:new-num-non-zero-coeffs}, i.e., showing \cref{eq:thm-numpoints-combined-results-lb,eq:thm-numpoints-combined-results-ub} , there is an $\sumTotA'$ such that for all $\sumTotA\geq \sumTotA$'
    \begin{equation}
 \revA{\frac{1 + \epsilon^-}{2}\frac{\sqrt[\DimInd]{\det(\Hess_{\optim{\mPoint}} \dummyFunctionA)}}{2(\DimInd+2)} \leq }\frac{1}{2}\frac{\dummyFunctionA(\mPoint)-\frac{1}{\dummyIntegerB_\sumTotA} \sum_{\sumIndB=1}^{\dummyIntegerB_\sumTotA} \dummyFunctionA(\mPoint^\sumTotA_{(\sumIndB)})}{\radiusA_\sumTotA^{2}} 
\revA{\leq} \frac{3 - \epsilon^+}{2}\frac{\sqrt[\DimInd]{\det(\Hess_{\optim{\mPoint}} \dummyFunctionA)}}{2(\DimInd+2)},
\end{equation}
\revA{where $\epsilon^-$ and $\epsilon^+$}

\revA{Furthermore,} we must have that $\mPoint_{(\dummyIntegerB_\sumTotA+1)}$ gets arbitrarily close to the boundary $\partial\ellipsoid_{\radiusA_\sumTotA}^\dummyBilForm(\optim{\mPoint})$ as $\sumTotA$ grows, i.e., 
\begin{equation}
 \revA{\frac{1 + \epsilon^-}{2}\frac{\sqrt[\DimInd]{\det(\Hess_{\optim{\mPoint}} \dummyFunctionA)}}{2(\DimInd+2)} \leq \frac{1}{2}\frac{\dummyFunctionA(\mPoint^\sumTotA_{(\dummyIntegerB_\sumTotA+1)})-\frac{1}{\dummyIntegerB_\sumTotA} \sum_{\sumIndB=1}^{\dummyIntegerB_\sumTotA} \dummyFunctionA(\mPoint^\sumTotA_{(\sumIndB)})}{\radiusA_\sumTotA^{2}} 
\leq \frac{3 - \epsilon^+}{2}\frac{\sqrt[\DimInd]{\det(\Hess_{\optim{\mPoint}} \dummyFunctionA)}}{2(\DimInd+2)},}
\label{eq:last-prop-Jm+1-bounds}
\end{equation}

Multiplying \revA{\cref{eq:last-prop-Jm+1-bounds}} by $\dummyIntegerB_0^{\frac{\DimInd+2}{\DimInd}}\bigl(\frac{\operatorname{vol}(\manifold)}{\ballVol_\DimInd} \bigr)^{\frac{2}{\DimInd}}$ and substituting the definition of $\radiusA_\sumTotA$ now gives the desired result.

        
\end{proof}

\section{Application of ESL to single particle Cryo-EM}
\label{sec:application-to-cryo-EM}
Remember from \cref{sec:introduction} that the aim in single particle Cryo-EM is to reconstruct the 3D structure of a biomolecule from single particle data. 
Assuming asymmetric homogeneous particles, one can formalise this as the task of recovering a 3D map $\vol \in \volSpace$ from 2D projection images $\img_1, \ldots, \img_\numImgs \in \imgSpace$. 
Here, $\volSpace$ and $\imgSpace$ are suitable function spaces that will be defined below and each $\img_{\imgInd} \in \imgSpace$ is a single sample of the $\imgSpace$–valued random variable $\stimg_{\imgInd}$ where
\begin{equation}  \label{eq:CryoEMInvProb-2}
    \stimg_\imgInd = \forward(\rot_\imgInd . \vol) + \stnoise_\imgInd 
    \quad\text{with $\rot_\imgInd\in \mathrm{SO}(3)$ and $\stnoise_{\imgInd} \sim \PDDataNoise^{\imgInd}$ for $\imgInd=1,\ldots, \numImgs$.}
\end{equation}
The operator $\forward \colon \volSpace \to \imgSpace$ (\emph{forward operator}) in \cref{eq:CryoEMInvProb-2} is assumed to be known and it associates a 3D map $\vol \colon \Real^3 \to \Real$ to a 2D image patch $\img \colon \Real^2 \to \Real$.
It is a mathematical model for image formation in a high resolution TEM in the sense that it simulates the 2D image data obtained in absence of noise from a particle with a given 3D map. 
See \cref{sec:maths-of-cryo-em} for an explicit expression for this operator.
Next, $\PDDataNoise^{\imgInd}$ denotes the statistical distribution for the $\imgSpace$-valued random variable $\stnoise_\imgInd$ that represents the observation noise.
Furthermore, each rotation $\rot_\imgInd\in \mathrm{SO}(3)$ acts on a 3D map according to the group action $\rot.\vol := \vol \circ \rot^{-1}$ for $\rot\in \mathrm{SO}(3)$.

In this section, we combine the ESL approach (\cref{alg:l2-scaled-lifting}) in an alternating updating scheme for solving the Cryo-EM inverse problem \cref{eq:CryoEMInvProb-2}, that can be generalized straightforwardly into a proper refinement method, i.e., also taking into account symmetric molecules, in-plane shifts and varying TEM optics CTF across 2D image patches.

\subsection{Measure-based lifting for rotation estimation in alternating updating}
\label{sec:new-joint-refinement}
Let $\dom\subset \Real^3$ be a bounded domain and assume our 3D maps reside in the function space $\volLtwoSpace$. 
If one chooses a regularisation functional $\regA \colon \volLtwoSpace\to\Real$ as
\begin{equation}\label{eq:MapRegFunc}
\regA(\vol) = \frac{1}{2\volRegParam_1} \| \vol\|_2^2 + \frac{1}{2\volRegParam_2} \|\rampKernel* \vol\|_2^2
\quad\text{for some $\volRegParam_1, \volRegParam_2 >0$,}
\end{equation}
where $\rampKernel \colon \volLtwoSpace\to\volLtwoSpace$ is a low-pass filter defined by its Fourier transform as $\fourier{\rampKernel}(\xi) := |\xi|$, then \revA{formally} the joint 3D map reconstruction and rotation estimation problem reads as 
\begin{align}
\label{eq:unlifted-problem-cryo-vol-rot}
    \inf_{\substack{\vol \in \volLtwoSpace \\\rot_1, \ldots, \rot_\numImgs \in \mathrm{SO}(3)}} \biggl\{ \sum_{\imgInd=1}^\numImgs \Bigl[ \frac{1}{2\noiseLevelParam} \bigl\|\forward (\rot_{\imgInd}.\vol)-\img_{\imgInd}\bigr\|_2^{2} \Bigr] + \frac{1}{2\volRegParam_1} \| \vol\|_2^2 + \frac{1}{2\volRegParam_2} \|\rampKernel* \vol\|_2^2 \biggr\}
\end{align}
for some noise level $\noiseLevelParam >0$.
We propose to solve \cref{eq:unlifted-problem-cryo-vol-rot} using the following alternating iterative scheme: 
\begin{align}
    \label{eq:-alternating-rots-update}
    \rot_{\imgInd}^{\iterInd+1} & := \argmin_{\rot \in \mathrm{SO}(3)} \biggl\{ \frac{1}{2\noiseLevelParam} \bigl\|\forward (\rot.\vol^{\iterInd})-\img_{\imgInd}\bigr\|_2^{2}  \biggr\} \quad \text{for $\imgInd= 1,\ldots,\numImgs$,} \\
    \vol^{\iterInd+1} &:= \argmin_{\vol \in \volLtwoSpace} \Bigl\{  \sum_{\imgInd=1}^\numImgs \frac{1}{2\sigma} \left\|\forward (\rot_{\imgInd}^{\iterInd+1} . \vol)-\img_{\imgInd}\right\|_2^{2} + \frac{1}{2\volRegParam_1} \| \vol\|_2^2 + \frac{1}{2\volRegParam_2} \|\rampKernel* \vol\|_2^2  \Bigr\}.
    \label{eq:-alternating-vol-update}
\end{align}
Solving \cref{eq:-alternating-vol-update} is typically straightforward. 
In the Cryo-EM setting, one can even express the solution in closed form -- as will be discussed in the following. 
Hence, the challenging part lies in \emph{jointly regularising and solving} \cref{eq:-alternating-rots-update} for the particle-specific rotations under high noise conditions. This will be addressed using the developed machinery from the previous section. 
The resulting joint recovery algorithm is summarized in \cref{alg:joint-refinement-lifting}. Below we elaborate on several steps in the method and introduce the notation.




\paragraph{Initialisation}
First note that because of the non-convexity of \cref{eq:unlifted-problem-cryo-vol-rot}, proper initialisation of either the 3D map or the rotations will be key. 
\textit{Ab initio} modelling has been developed to provide 3D maps $\vol^0$ \cite{bandeira2017estimation,bandeira2020non,bhamre2015orthogonal,levin20183d,sharon2020method,shkolnisky2012viewing,singer2011three,wang2013orientation} and although in several instances initial estimates for the rotations are also provided, this is not always the case. Hence, we will assume there is no \textit{a priori} information on rotations, i.e., we start the algorithm by estimating them.

\paragraph{Updating particle specific rotations}
Key components for using ESL (\cref{alg:l2-scaled-lifting}) for estimating particle specific rotations by \revA{jointly regularising and} solving \cref{eq:-alternating-rots-update} is the choice of sampling set $\rotsSampling$ -- or rather choice of initial set $\rotsSampling_0$ and a refinement strategy -- and the choices of the regularisation parameter $\regParam >0$ and the scaling factor $\exponentA >0$ in the regularised formulation \cref{eq:lifted-function-discr} of \cref{eq:-alternating-rots-update}.

In this work, $\rotsSampling_0$ is obtained from a symmetric uniform mesh over $\Sphere^3$ of size $3\,642$ whose uniformity has been verified numerically \cite[Table 10]{womersley2018efficient}. 
Realising that $\Sphere^3$ is a double cover of $\mathrm{SO}(3)$ and that the two can be identified \cite[Lemma 1]{graf2012unified}, yields a sampling set $\rotsSampling_0 \subset \mathrm{SO}(3)$ with $|\rotsSampling_0| = 1\,821$ and better uniformity -- in terms of integration error -- than mathematically derived sampling sets that are typically much smaller \cite{graf2009sampling}. 
For generating larger sampling sets, a triangulation of $\rotsSampling_0$ is constructed and midpoint refinement is used to construct 7 smaller simplices from each simplex in the triangulation \cite[Fig. 2]{todorov2013optimal}. The procedure can be repeated until the desired mesh resolution is obtained. 

When choosing $\regParam$ and $\exponentA$ one should consider that the effect of the former can be compensated for by the latter and vice versa for a fixed sampling set -- which is the case in a real-word scenario. Therefore, even though one could be tempted to pick $\exponentA$ yielding the lowest bound in \cref{thm:distance-bound-post-lifting}, one should rather consider \cref{thm:new-num-non-zero-coeffs} for picking $\exponentA$, i.e., prioritize on controlling the amount of non-zero coefficients in the sampling set. Nevertheless, here one also sees that choosing $\exponentA\approx2/3$ will give most control, i.e., least variability between the upper and lower bound, and hence will also be the preferred choice from a practical point of view.

For $\regParam$, the approximation in \cref{prop:choice-regulariser-lifting} is then used due to the unavailability of the Riemannian Hessian of the objective functions \cref{eq:-alternating-rots-update} -- note the bracket notation for the labeling of the rotations as in the proposition is also used in  \cref{alg:joint-refinement-lifting}. Finally, \cref{sec:numerics} investigates if one can actually choose $\exponentA$ this high and if so, how many points with non-zero coefficients one needs for proper performance, i.e., how to pick $\dummyIntegerB_0$ in \cref{prop:choice-regulariser-lifting}.

\paragraph{Updating the 3D map}
For the 3D map update it is key to have an efficient algorithm for solving \cref{eq:-alternating-vol-update} and properly chosen parameters $\noiseLevelParam,\volRegParam_1$ and $ \volRegParam_2$.

Our 3D map update builds on the approach taken in RELION \cite{scheres2012relion}.
In particular, we extend the computationally efficient closed-form solution of the RELION 3D map update scheme with an additional regulariser of the form $\vol \mapsto \frac{1}{2\volRegParam_2} \|\rampKernel* \vol\|_2^2$ term: 
\begin{align}
\vol^{\iterInd+1} 
    &:= \argmin_{\vol \in \volLtwoSpace} \Bigl\{  \sum_{\imgInd=1}^\numImgs \frac{1}{2\sigma} \bigl\|\forward (\rot_{\imgInd}^{\iterInd+1}\!. \vol)-\img_{\imgInd}\bigr\|_2^{2} + \frac{1}{2\volRegParam_1} \| \vol\|_2^2 + \frac{1}{2\volRegParam_2} \|\rampKernel * \vol\|_2^2  \Bigr\}
    \\
    & = a^{\iterInd+1}
    \ast 
    \biggl[ \sum_{\imgInd= 1}^\numImgs \frac{\volRegParam_1}{\noiseLevelParam} (\rot_\imgInd^{\iterInd+1})^{-1}\!.\,\forward^{\ast}(\img_\imgInd)
    \biggr],
\end{align}
where $a^{\iterInd+1}$ is defined by its Fourier transform as
\begin{equation}
\label{eq:volume-filter-stage-2}
    \fourier{a}^{\iterInd+1} := \frac{1}{\displaystyle{\sum_{\imgInd= 1}^\numImgs}  \frac{\volRegParam_1}{\noiseLevelParam} \fourier{\bigl((\rot_\imgInd^{\iterInd+1})^{-1}\!.\, \forward^{\ast}(\RerCTF)\bigr)} + \frac{\volRegParam_1}{\volRegParam_2} \fourier{\rampKernel}^2 + 1}
    \quad\text{with $\RerCTF$ as in \cref{eq:ctf}.}
\end{equation}
As shown above, the addition of the second regulariser does not complicate the numerics substantially, as it comes down to an additional term in the denominator of \cref{eq:volume-filter-stage-2}, and the calculations for the reconstruction will happen in the Fourier domain anyways.
Finally, $\noiseLevelParam$ can typically be determined from information on the noise level and the regularisation parameters $\volRegParam_1$ and $\volRegParam_2$ can be chosen from the norm of the (initial) 3D map:
\begin{equation}
\label{eq:vol-param-choice}
    \noiseLevelParam:= \frac{1}{\numImgs}\sum_{\imgInd=1}^\numImgs \operatorname{Var}(\img_\imgInd) \quad \text{and}\quad \volRegParam_1:=\volRegParam_2:= \|\vol^0\|_2^2,
\end{equation}
where $\operatorname{Var}(\img_\imgInd)$ is the variance of the pixel-values of the $\imgInd$th image.

\begin{remark}
We use lifting to re-phrase the non-convex optimisation problem for estimating particle-specific rotations in Cryo-EM.
This leads to \cref{alg:joint-refinement-lifting} that is then tested on a realistic single particle Cryo-EM setting to evaluate its potential for incorporation in similar joint refinement schemes. In other words, in this work we do not intend to propose a state-of-the-art joint estimation scheme with respect to the current 3D map update routine. As an example, joint schemes with more carefully chosen 3D map updates steps -- like recent ones based on more sophisticated priors \cite{gilles2022molecular} -- will most likely yield a better 3D map reconstruction. 
However, this is beyond the scope of this article.
\end{remark}




\begin{algorithm}[h!]
\caption{ESL-based joint 3D map reconstruction and rotation estimation }
\label{alg:joint-refinement-lifting}
\begin{algorithmic}
\STATE{\textit{Initialisation}: $\rotsSampling\subset\mathrm{SO}(3)$, $\vol^0\in \volLtwoSpace$, $\noiseLevelParam,\volRegParam_1, \volRegParam_2 >0$, $\dummyIntegerB_0 \geq 1$, $\exponentA\in(\frac{1}{4},\frac{2}{3})$}
\STATE{$\iterInd:=0$}
\STATE{$\dummyIntegerB :=  
\Bigl\lfloor \dummyIntegerB_0 \numRots^{\frac{2- 3\exponentA}{5}} \ \Bigr\rfloor$}
\par\medskip
\WHILE{not converged}
\par\smallskip
\STATE{$\displaystyle{\regParam_\imgInd^{\iterInd+1} := \frac{1}{2}\dummyIntegerB_0\numRots^{\frac{2\revA{+ 2\exponentA}}{\revA{5}}} \Bigl( \frac{1}{2\noiseLevelParam} \bigl\|\forward (\rot^\iterInd_{(\dummyIntegerB+1)} . \vol^\iterInd)-\img_{\imgInd} \bigr\|_2^{2} -  \frac{1}{\dummyIntegerB}\sum_{\sumIndB=1}^{\dummyIntegerB} \frac{1}{2\noiseLevelParam} \bigl\|\forward (\rot^\iterInd_{(\sumIndB)} . \vol^\iterInd)-\img_{\imgInd} \bigr\|_2^{2} \Bigr)}$}
\par\medskip
\STATE{$\displaystyle{\rotsDensCoef^{\iterInd+1}_\imgInd :=  \prod{}_{\simplex^{\numRots}}\Bigl(- \frac{\numRots^{\exponentA}}{\regParam_\imgInd^{\iterInd+1}}\Bigl( \frac{1}{2\noiseLevelParam} \bigl\|\forward (\rot . \vol^\iterInd)-\img_{\imgInd} \bigr\|_2^{2}\Bigr)_{\rot \in \rotsSampling} \Bigr) \quad\text{for $\imgInd = 1,\ldots,\numImgs$}}$}
\par\medskip
\STATE{$\iterIndB:=0$}
\par\medskip
\STATE{$\displaystyle{\rot_\imgInd^{\iterInd+1,0}:=\displaystyle{\argmax_{\rot\in \rotsSampling}}\, (\rotsDensCoef^{\iterInd+1}_\imgInd)_{\rot} \quad\text{for $\imgInd = 1,\ldots,\numImgs$}}$}
\par\medskip
\WHILE{not converged}
\par\smallskip
\STATE{$\displaystyle{\rot_\imgInd^{\iterInd+1,\iterIndB+1}:=\exp _{\rot_\imgInd^{\iterInd+1,\iterIndB}}\Bigl(\sum_{\rot \in \mathcal{X}} (\rotsDensCoef^{\iterInd+1}_\imgInd)_{\rot} \log_{\rot_\imgInd^{\iterInd+1,\iterIndB}}(\rot)\Bigr)  \quad \text{for }\imgInd = 1,\ldots,\numImgs}$}
\STATE{$\iterIndB:=\iterIndB+1$}
\ENDWHILE
\par\medskip
\STATE{$\rot_\imgInd^{\iterInd+1}:= \rot_\imgInd^{\iterInd+1,l}$}
\par\medskip
\par\smallskip
\STATE{$\displaystyle{\vol^{\iterInd+1} : =  a^{\iterInd+1} * \sum_{\imgInd= 1}^\numImgs \frac{\volRegParam_1}{\noiseLevelParam} (\rot_\imgInd^{\iterInd+1})^{-1}. \forward^{\ast}(\img_\imgInd) \quad\text{for $a^{\iterInd+1}$ from \cref{eq:volume-filter-stage-2}}}$}
\STATE{$\iterInd:= \iterInd+1$}
\ENDWHILE
\end{algorithmic}
\end{algorithm}

\section{Numerical experiments}
\label{sec:numerics}

\revA{As the theoretical results (\cref{thm:new-num-non-zero-coeffs,thm:distance-bound-post-lifting,prop:choice-regulariser-lifting}) are asymptotic, i.e., they will only be observed for a sufficiently large sampling sets, we aim to test whether these asymptotics kick in for Cryo-EM data. In particular, we test the performance of ESL (\cref{alg:l2-scaled-lifting}) as a \revA{regularising} global optimisation method and its potential when included in a joint recovery method (\cref{alg:joint-refinement-lifting}).}
\revA{All of the experiments are implemented in Python 3.6 and run on a 2 GHz Quad-Core Intel Core i5 with 16GB RAM.}





\begin{figure}[h!]
    \centering
    \begin{subfigure}{\dimexpr0.30\textwidth+20pt\relax}
    \includegraphics[width=\dimexpr\linewidth-20pt\relax]{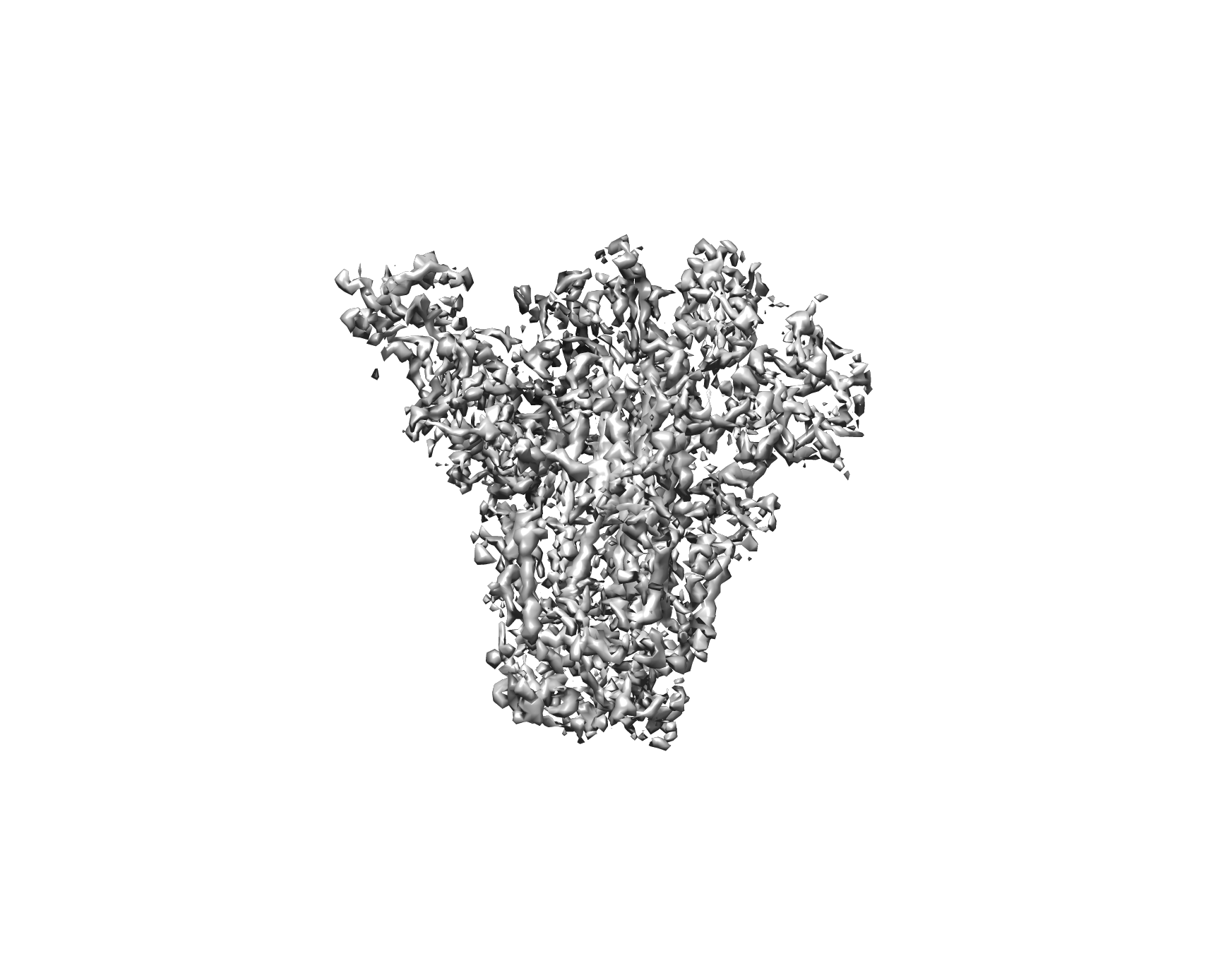}
    \caption{\revA{Iso-surface of} $\vol^{\mathrm{GT}}$}
    \label{fig:vol-gt}
    \end{subfigure}
    \hfill
    \begin{subfigure}{0.30\textwidth}
    \includegraphics[width=\linewidth]{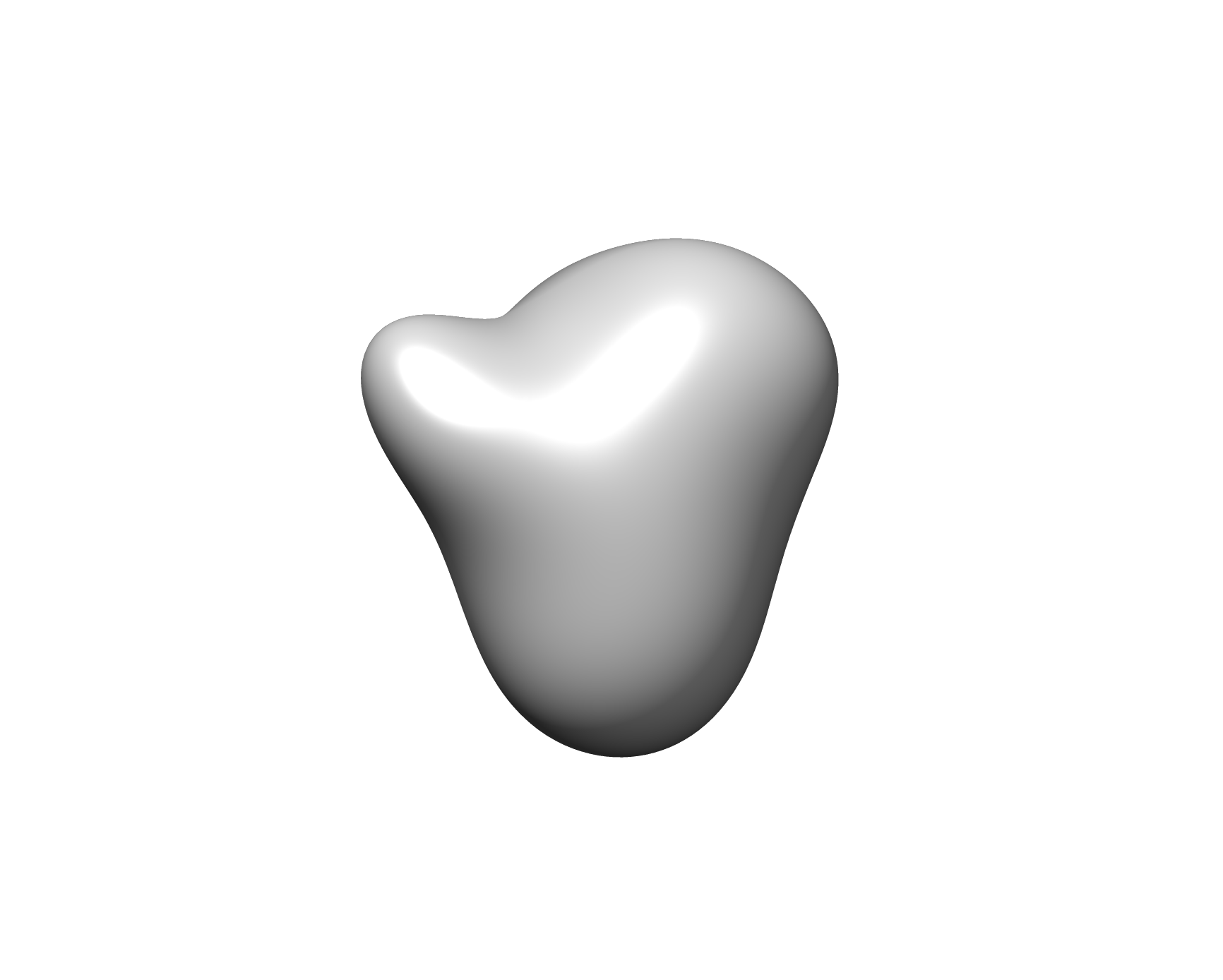}
    \caption{\revA{Iso-surface of} $\vol^{\mathrm{INIT}}$}
    \label{fig:vol-init}
    \end{subfigure}
    \hfill
    \begin{subfigure}{0.30\textwidth}
    \centering
    \includegraphics[width=0.81\linewidth]{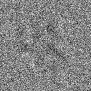}
    \caption{Example 2D image patch}
    \label{fig:projection}
    \end{subfigure}
    \caption{Synthetic dataset generated from the delta variant of the SARS-CoV-2 spike protein.}
    \label{fig:data}
\end{figure}

\paragraph{Dataset}
The algorithms will be tested on a synthetic single particle Cryo-EM dataset of the spike protein of the SARS-CoV-2 (Delta variant). Using Chimera \cite{pettersen2004ucsf}, a ground-truth 3D map $\vol^{\mathrm{GT}}$ (\cref{fig:vol-gt}) is generated on a $91\times91\times91$ grid with a voxel size of $0.21667\,\text{nm}$. This 3D map is also used to generate $\numImgs~=~2\,048$ 2D TEM image patches (projection images) $(\img_\imgInd)_{\imgInd=1}^\numImgs$ of size $91\times91$ (\cref{fig:projection}). These image patches are generated by first randomly rotating the 3D map using rotations  $(\rot^{\mathrm{GT}}_{\imgInd})_{\imgInd=1}^\numImgs$ that are uniformly distributed over $\mathrm{SO}(3)$, then applying the TEM forward operator \cref{eq:cryo-operator} as implemented in ASPIRE\footnote{\href{http://spr.math.princeton.edu/}{http://spr.math.princeton.edu/}}. The forward operator is parametrized by the electron wavenumber $\waveNumber = 0.25~\text{nm}^{-1}$, which corresponds to the (relativistically corrected) wavelength of a $200$~kV electron. The TEM optics model has defocus $\defocus = 1.5~\mu\text{m}$ and spherical aberration  $\spAberration = 2~\text{mm}$. Furthermore, the amplitude contrast is set to $\ampContrast = 0.1$. 
These parameters are chosen to mimic a realistic setting (cf. \cite[Table 2]{scheres2012relion}). Additionally, the low-pass filter $A(\xi)$ in the CTF is omitted. For each TEM projection image, noise is added corresponding to a noise level where signal-to-noise-ratio (SNR) is 1/16.
Finally, The ground true 3D map is also blurred by convolving against a normalized Gaussian kernel with standard deviation $10$ to get a realistic initial density $\vol^{\mathrm{INIT}}$ (\cref{fig:vol-init}). 


\paragraph{Overview of Experiments}
\revA{The asymptotics in} \Cref{thm:new-num-non-zero-coeffs,thm:distance-bound-post-lifting,alg:l2-scaled-lifting} are tested through rotation estimation of the 3D map $\vol^{\mathrm{GT}}$ for all  $\numImgs$ 2D TEM image patches  (\cref{sec:scaling-errors}) using a fixed regularisation parameter $\regParam$, but with varying sizes of sampling sets $\rotsSampling$ and varying values for the scaling parameter $\exponentA$. Because of the high resolution of $\vol^{\mathrm{GT}}$, the optimisation landscape over $\mathrm{SO}(3)$ will be rougher \revA{-- in the sense of the landscape having larger gradients --} than when starting from a lower-resolution (initial) volume. \revA{The main reason here is to test whether \cref{alg:l2-scaled-lifting} as a regularised rotation finding algorithm is able to perform well in a late stage of a joint recovery algorithm, e.g., \cref{alg:joint-refinement-lifting}. That is, better at finding the actual ground truth rotation rather than the global minimiser of the rotation loss.}


Once some best practices have been gathered, it is time to check the performance of \cref{alg:l2-scaled-lifting} in an iterative fashion through \cref{alg:joint-refinement-lifting} and investigate whether on can control the amount of points in the sampling set through \cref{prop:choice-regulariser-lifting} in a realistic setting (\cref{sec:join-volume-rots-experiment}).







\paragraph{Cryo-EM-related considerations}
For evaluating rotations one should consider that the joint 3D map reconstruction and rotation estimation problem can be solved up to a global rotation or reflection. Hence, any obtained rotation will be aligned with the ground truth rotations through finding a minimal orthogonal transformation as done in \cite[(5.1)]{singer2011three}. 
Then, after alignment, we use the distance in the $\mathrm{SO}(3)$ metric and related distance measures, e.g., the Wasserstein-2 distance to compare rotations $(\rot_{\imgInd})_{\imgInd=1}^\numImgs$ against the ground truth, i.e., $\distance_{\mathrm{SO}(3)}(\rot_{\imgInd},\rot^{\mathrm{GT}}_{\imgInd})$. \revA{Finally, remember from \cref{sec:lifting-for-cryo-rotations} that, in general, the ground truth rotation does not coincide with any minimiser of the rotation estimation problem -- although one expects that these coincide fairly well --, which was the main motivation of ESL to regularise the optimisation problem. Hence, we have to be careful when making claims about our theoretical results.}




\subsection{Testing asymptotic behaviour on Cryo-EM data}
\label{sec:scaling-errors}

The performance of ESL (\cref{alg:l2-scaled-lifting}) and the \revA{asymptotics} in the theoretical results \cref{thm:new-num-non-zero-coeffs,thm:distance-bound-post-lifting} are tested for a fixed $\regParam$ and varying $\rotsSampling$ and $\exponentA$ through running one rotation update step of \cref{alg:joint-refinement-lifting}.
More precisely, for fixed $\regParam$ -- possibly parametrized by $\dummyIntegerB_0$ as suggested in \cref{eq:thm-new-num-non-zero-coeffs-gamma} -- the sampling set $\rotsSampling$ and the scaling parameter $\exponentA$ are the relevant factors in the theoretical bounds. 

We know from \cref{sec:new-joint-refinement} that fixing $\regParam$ is by no means restrictive as the effect of $\regParam$ can be compensated for by $\exponentA$ and vice versa. Indeed, either could be varied while the other is fixed to establish any weight for the Tikhonov regulariser in \cref{eq:lifted-function-discr}. In other words, the effect of larger $\regParam$ on the amount of points with non-zero coefficients -- or equivalently the effect of $\dummyIntegerB_0$ -- can also be studied by varying $\exponentA$. 

Two main questions follow naturally:
\begin{enumerate}
    \item Is computing the Riemannian barycentre from the obtained weights actually beneficial in practice compared to just picking the rotation with minimal cost?
    \item Do the \revA{asymptotic} theoretical bounds on (a) the algorithm's number of non-zero coefficients and (b) the error \revA{kick in} in practice?
\end{enumerate}


Additionally, in practice there is an argument to make to pick $\exponentA$ as large as possible.
\Cref{thm:new-num-non-zero-coeffs} suggests that for $\exponentA\approx\frac{2}{3}$ and as long as the quality of the sampling set admits it, the amount of points with non-zero coefficients can be controlled more precisely as the upper and lower bound hardly vary for larger $\rotsSampling$. Again, because the effect of $\exponentA$ can be compensated by $\regParam$ -- as mentioned above -- it then remains to pick the latter properly. Two additional questions follow:
\begin{enumerate}
\setcounter{enumi}{2}
    \item Can we pick $\exponentA\approx\frac{2}{3}$, i.e., does the sampling set generating procedure described in \cref{sec:application-to-cryo-EM} satisfy the technical condition underlying the main theoretical results?
    \item How to pick $\regParam$ -- or how to choose $\dummyIntegerB_0$?
\end{enumerate}

\paragraph{Metrics}
Besides the average of the distances $\distance_{\mathrm{SO}(3)}(\rot^1_{\imgInd},\rot^{\mathrm{GT}}_{\imgInd})$, where $\rot^1_{\imgInd}$ is the rotation obtained from the update scheme, several additional metrics will be considered for each $\rotsSampling$ and each $\exponentA$ to investigate the first two questions.
\begin{enumerate}
    \item For testing \cref{alg:l2-scaled-lifting}: the average of the distances $\distance_{\mathrm{SO}(3)}(\rot^{1,0}_{\imgInd},\rot^{\mathrm{GT}}_\imgInd)$ of the rotation with largest coefficient, i.e., the initialisation of Riemannian gradient descent $\rot^{1,0}_{\imgInd}$, will be considered.
    \item (a) For testing \cref{thm:new-num-non-zero-coeffs}: the average of the $\ell^0$-norms $\|\rotsDensCoef^1_\imgInd\|_0$. 
    \par
    (b) For testing \cref{thm:distance-bound-post-lifting}: the average of the Wasserstein-2 distances $\wasserstein_2(\measure^1_\imgInd,\dirac_{\rot^{\mathrm{GT}}_\imgInd})$ of the measure $\measure^1_\imgInd:= \sum_{\rot\in\rotsSampling} (\rotsDensCoef^1_\imgInd)_\rot \dirac_\rot$ to the Dirac measure concentrated on the ground truth, will be considered. The reason is twofold: (1) this is the main ingredient for the scaling behaviour and (2) the Wasserstein distance is expected to be less prone to suffering from the discrepancy of the ground truth rotation not being the global minimiser than the estimate $\rot_\imgInd^1$, which is hypothesized to get closer to the ground truth with larger bias, i.e., larger support of the measure $\measure^1_\imgInd$.
\end{enumerate}
As we shall see, the latter two questions will follow from these metrics.

\paragraph{Algorithm settings}
\revA{We} use $\vol^0:= \vol^\mathrm{GT}$, $\noiseLevelParam$ as in \cref{eq:vol-param-choice} and $\regParam:=30$. 
Here the latter is chosen large enough such that the lifting scheme does not collapse into a single rotation. 
For computing the Riemannian barycentre $20$ descent steps are taken, \revA{which is empirically sufficient to converge.}.



Three sampling sets and six scaling parameters will be considered. 
The three sampling sets $\rotsSampling_1$, $\rotsSampling_2$ and $\rotsSampling_3$ are generated from iterative midpoint refinement of the sampling set $\rotsSampling_0$ with $|\rotsSampling_0|=1\,821$ as mentioned in \cref{sec:application-to-cryo-EM}. 
We have $|\rotsSampling_1|= 14\,761$, $|\rotsSampling_2|= 114\,564$ and $|\rotsSampling_3|= 860\,069$. 
The scaling factors are chosen as $\eta \in\{ 0.26, 0.3, 0.4, 0.5, 0.6, 0.66 \}$, so that the whole feasible range $(\frac{1}{4}, \frac{2}{3})$ is taken into account. 




\paragraph{Riemannian barycentre vs.~minimal cost rotation}
Regarding question~1, we computed averages of the distances $\distance_{\mathrm{SO}(3)}(\rot^{1,0}_{\imgInd},\rot^{\mathrm{GT}}_\imgInd)$ and  $\distance_{\mathrm{SO}(3)}(\rot^1_{\imgInd},\rot^{\mathrm{GT}}_{\imgInd})$. 
The results in \cref{tab: distances-0,tab: distances-1} indicate that averaging using \cref{alg:l2-scaled-lifting} performs better by a factor of 2 compared to picking the rotation with maximal weight, for which one normally needs about 8 times as many points on $\mathrm{SO}(3)$ to obtain the same accuracy. 

\revA{Additionally, \cref{tab: distances-1} strongly suggests that the bias introduced by the lifting approach yields a better rotation estimation than finding the global minimiser. Indeed, the global minimiser will be approximated best by choosing $\exponentA = 0.66$. However, the bold face values and not the the right-most column -- where we also average over a minimial amount of points with non-zero coefficients, cf. \cref{tab: number-non-zero-coefs-results} -- yield the lowest distance to the ground truth rotation. In other words, these results support our hypothesis that introducing bias through regularisation in the global minimisation indeed yields better rotation estimation.}



\begin{table}[h!]
    \centering
    \caption{Distance $\distance_{\mathrm{SO}(3)}(\rot^{1,0}_{\imgInd},\rot^{\mathrm{GT}}_\imgInd)$ of the rotation with minimal cost to the ground truth for different sizes of sampling set $\rotsSampling$ (baseline). Numbers are reported in degrees and averaged over all rotations.}
        \label{tab: distances-0}
    \scriptsize
        \begin{tabular}{r c}
                \toprule
                $\numRots$ & $\distance_{\mathrm{SO}(3)}(\rot^{1,0}_{\rotsSampling},\optim{\rot}) \; (^\circ)$ \\ 
                \midrule
                14\,761 & $5.026 \pm 1.727$ \\ 
                114\,564 & $2.758 \pm 0.948$\\ 
                860\,069 & $1.541 \pm 0.605$ \\ 
                \bottomrule
        \end{tabular}
\end{table}

\begin{table}[h!]
        \centering
        \scriptsize
        \caption{
        Distance $\distance_{\mathrm{SO}(3)}(\rot^1_{\imgInd},\rot^{\mathrm{GT}}_{\imgInd})$ of the (Riemannian barycentre) rotation obtained by \cref{alg:l2-scaled-lifting}
        to the ground truth for different sizes of sampling set $\rotsSampling$ and different choices of $\exponentA$ (ours). Varying $\exponentA$ while keeping $\regParam=30$ fixed directly influences the number of non-zero coefficients in the solution. Numbers are again reported in degrees and averaged over all rotations. For properly chosen $\exponentA$, the method outperforms the baseline approach \cref{tab: distances-0} by approximately a factor of two (bold).}        
        \label{tab: distances-1}
        \begin{tabular}{r cccccc}
                \toprule
                $\numRots$ & $\exponentA=0.26$  & $\exponentA=0.3$ & $\exponentA=0.4$ & $\exponentA=0.5$  & $\exponentA=0.6$ & $\exponentA=0.66$ \\ 
                \midrule
                14\,761 & $37.628 \pm 16.024$ & $32.005 \pm 14.674$ & $18.226 \pm 10.407$ & $7.116 \pm 5.548$ & $2.54 \pm 1.821$ & $\mathbf{2.026 \pm 1.178}$ \\ 
                114\,564 & $4.959 \pm 4.154$ & $2.953 \pm 2.36$ & $1.613 \pm 0.845$ & $1.229 \pm 0.583$ & $\mathbf{1.14 \pm 0.534}$ & $1.2 \pm 0.575$ \\ 
                860\,069 & $1.236 \pm 0.588$ & $1.108 \pm 0.524$ & $0.925 \pm 0.428$ & $\mathbf{0.858 \pm 0.389}$ & $0.87 \pm 0.406$ & $0.938 \pm 0.446$ \\ 
                \bottomrule
        \end{tabular}
\end{table}

\paragraph{Bounding the number of non-zero coefficients}
For part (a) of question 2, the average number of points with non-zero coefficients, i.e., $\|\rotsDensCoef^1_\imgInd\|_0$, used to compute the Riemannian barycentres is provided in \cref{tab: number-non-zero-coefs-results}. From the results in \cref{tab: number-non-zero-coefs-results} it is not directly clear whether the scaling behaviour predicted by \cref{thm:new-num-non-zero-coeffs} \revA{is observed}. As $\|\rotsDensCoef^1_\imgInd\|_0$ does not scale as predicted, but rather decreases with $\numRots$ for all instances of $\exponentA$, one would say that either there are not enough points for the scaling to kick in or the technical condition has not been satisfied. 

However, we would argue that \revA{we do observe the asymptotic behaviour predicted by theory}. That is, the predicted range is fairly broad so one should be careful immediately comparing between different sampling sets. A better indication of the validity of the scaling is to compare the ratio of upper bounds between different $\exponentA$-values with the ratio of the observed number of points for a given sampling set. The ratios relative to $\exponentA=0.66$ are provided in \cref{tab: number-non-zero-coefs-results-theory} and show that for all values of $\exponentA$ the empirical ratio gets closer to the predicted values for larger sampling sets, which is in line with \cref{thm:new-num-non-zero-coeffs}.





\begin{table}[h!]
        \centering
        \caption{The average number of points with non-zero coefficients, i.e., $\|\rotsDensCoef^1_\imgInd\|_0$, for different sizes of sampling set $\rotsSampling$ and different numbers of non-zero coefficients. Varying the latter is realized through a varying $\exponentA$, but fixed $\regParam=30$. Bold values correspond to the value or $\exponentA$ that gave minimal average error for each sampling set, cf. \cref{tab: distances-1}. Without normalization -- as done in \cref{tab: number-non-zero-coefs-results-theory} -- it is not immediately clear that the results are in line with theory for larger sampling sets.}
        \label{tab: number-non-zero-coefs-results}
        \scriptsize
        \begin{tabular}{r cccccc}
                \toprule
                $\numRots$ & $\exponentA=0.26$  & $\exponentA=0.3$ & $\exponentA=0.4$ & $\exponentA=0.5$  & $\exponentA=0.6$ & $\exponentA=0.66$ \\ 
                \midrule
                14\,761 & $760.77 \pm 92.25$ & $517.52 \pm 68.78$ & $173.46 \pm 32.22$ & $47.068 \pm 10.762$ & $16.167 \pm 2.953$ & $\mathbf{9.889 \pm 1.863}$ \\ 
                114\,564 & $253.32 \pm 51.7$ & $151.67 \pm 25.39$ & $55.684 \pm 7.449$ & $22.732 \pm 3.356$ & $\mathbf{9.655 \pm 1.916}$ & $5.948 \pm 1.413$ \\ 
                860\,069 & $204.95 \pm 25.76$ & $135.53 \pm 17.07$ & $49.646 \pm 6.585$ & $\mathbf{19.032 \pm 3.031}$ & $7.847 \pm 1.771$ & $4.721 \pm 1.339$  \\ 
                \bottomrule
        \end{tabular}
\end{table}

\begin{table}[h!]
        \centering
        \caption{A comparison of theoretical and empirical ratios of average number of points with non-zero coefficients from the results in \cref{tab: number-non-zero-coefs-results}. For larger sampling sets the discrepancy between theory and practice closes in.}
        \label{tab: number-non-zero-coefs-results-theory}
        \scriptsize
        \begin{tabular}{r r ccccc}
                \toprule
                $\numRots$ & $\frac{\|\optim{\rotsDensCoef}_\imgInd\|_0\mid_{\exponentA=(\cdot)}}{\|\optim{\rotsDensCoef}\|_0 \mid_{\exponentA=0.66}}$ & $\exponentA=0.26$  & $\exponentA=0.3$ & $\exponentA=0.4$ & $\exponentA=0.5$  & $\exponentA=0.6$ \\ 
                \midrule
                14\,761 &  predicted & 10.014 &  7.953 & 4.471 & 2.513 &  1.413 \\ 
                &  observed & 76.931 & 52.333 & 17.541 &  4.76  & 1.635 \\ 
                114\,564 & predicted & 16.375 & 12.381 & 6.155 & 3.06  & 1.521 \\ 
                & observed & 42.589 & 25.499 & 9.362 & 3.822 & 1.623 \\
                860\,069 & predicted & 26.564 & 19.136 & 8.429 & 3.713 & 1.635\\ 
                & observed & 43.412 & 28.708 & 10.516 & 4.031 & 1.662 \\
                \bottomrule
        \end{tabular}
\end{table}


\paragraph{Bounding the error}
For part (b) of question 2, getting back to \cref{tab: distances-1}, again it is not directly clear whether or not the scaling behaviour of the error bound in \cref{thm:distance-bound-post-lifting} \revA{is observed here}. The scaling seemingly fails as values for $\exponentA$ yielding more points with non-zero coefficients obtain a lower average error for larger sampling sets (bold numbers in \cref{tab: distances-1}). 

However, this effect is arguably due to \revA{the fact that the error is computed with respect to the ground truth and not the global minimiser. Subsequently, as discussed before, the  Wasserstein-2 distance is a more robust way of checking the scaling behaviour of the error in the absence of a global minimiser. These values }are provided in \cref{tab: distances-W2}. Similarly to the number of non-zero coefficients it is not obvious whether or not the scaling behaviour is kicking in. Comparing ratios -- similarly as before -- we observe in \cref{tab: distances-W2-ratios} again that theory and empirical results agree better for larger sampling sets. Overall, this experiment is in line with \cref{thm:distance-bound-post-lifting}.




\begin{table}[h!]
        \centering
        \caption{The average of the distances $\wasserstein_2(\sum_{\rot\in\rotsSampling} (\rotsDensCoef^1_\imgInd)_\rot \dirac_\rot,\dirac_{\rot^{\mathrm{GT}}_\imgInd})$ in degrees $(^\circ)$ for different sizes of sampling set $\rotsSampling$ and different numbers of non-zero coefficients, where the latter is realized through a varying $\exponentA$, but fixed $\regParam=30$.  Bold values correspond to the value or $\exponentA$ that gave minimal average error for each sampling set, cf. \cref{tab: distances-1}. Without normalization -- as done in \cref{tab: distances-W2-ratios} -- it is not immediately clear that the results are in line with theory for larger sampling sets.}
        \label{tab: distances-W2}
        \scriptsize
        \begin{tabular}{rcccccc}
                \toprule
                $\numRots$ & $\exponentA=0.26$  & $\exponentA=0.3$ & $\exponentA=0.4$ & $\exponentA=0.5$  & $\exponentA=0.6$ & $\exponentA=0.66$ \\ 
                \midrule
                14\,761 & $108.68 \pm 7.196$ & $101.495 \pm 8.781$ & $74.179 \pm 13.95$ & $34.403 \pm 14.731$ & $12.283 \pm 4.68$ & $\mathbf{9.481 \pm 2.439}$ \\ 
                114\,564 & $24.769 \pm 12.06$ & $14.76 \pm 6.628$ & $8.44 \pm 1.524$ & $6.175 \pm 0.559$ & $\mathbf{4.66 \pm 0.508}$ & $3.984 \pm 0.538$ \\ 
                860\,069 & $6.58 \pm 0.589$ & $5.731 \pm 0.489$ & $4.129 \pm 0.378$ & $\mathbf{3.049 \pm 0.33}$ & $2.318 \pm 0.339$ & $1.998 \pm 0.387$ \\ 
                \bottomrule
        \end{tabular}
\end{table}
\begin{table}[h!]
        \centering
        \caption{A comparison of theoretical and empirical ratios of Wasserstein-2 distances from the results in \cref{tab: distances-W2}. For larger sampling sets the discrepancy between theory and practice closes in.}
        \label{tab: distances-W2-ratios}
        \scriptsize
        \begin{tabular}{rrccccc}
                \toprule
                $\numRots$ &  $\frac{\wasserstein_2(\optim{\measure}_\rotsSampling,\dirac_{\rot^*})\mid_{\exponentA=(\cdot)} }{ \wasserstein_2(\optim{\measure}_\rotsSampling,\dirac_{\rot^*})\mid_{\exponentA=0.66}}$ & $\exponentA=0.26$  & $\exponentA=0.3$ & $\exponentA=0.4$ & $\exponentA=0.5$  & $\exponentA=0.6$ \\ 
                \midrule
                14\,761 &  predicted & 2.155 & 1.996 & 1.647 & 1.36 & 1.122 \\ 
                &  observed & 11.463 & 10.705 & 7.824 & 3.629 & 1.296 \\ 
                114\,564 & predicted & 2.539 & 2.313 & 1.833 & 1.452 & 1.15  \\ 
                & observed & 6.217 & 3.705 & 2.118 & 1.55 & 1.17  \\
                860\,069 & predicted & 2.984 & 2.675 & 2.035 & 1.548 & 1.178 \\ 
                & observed & 3.293 & 2.868 & 2.067 & 1.526 & 1.16 \\
                \bottomrule
        \end{tabular}
\end{table}

\paragraph{Choosing the parameter $\exponentA$}
For question~3, it is currently unknown whether our choice of sampling sets fulfills the technical condition of local low-discrepancy with $\exponentA\approx 2/3$. However, empirical evidence of the scaling behaviour (\cref{tab: number-non-zero-coefs-results-theory,tab: distances-W2-ratios}) suggests that this choice is justified for the two largest sampling sets. 

\paragraph{Choosing regularisation strength}
Finally, for question~4, choosing the regularisation parameter $\regParam$, \cref{tab: distances-1,tab: number-non-zero-coefs-results} suggest that there is a sweet spot for the number of non-zero coefficients and the regularisation parameter $\regParam$, i.e., the bold numbers in the tables suggest that \revA{inducing bias through regularisation yields better rotation estimation and that} the number of non-zero coefficients should not be chosen as small as possible but not too large either and this should be reflected in the choice of $\regParam$.

\paragraph{Best practices}
Overall, the results from this experiment lead to the following best practice for the current dataset:
\begin{enumerate}
    \item Use a sufficiently large sampling set, e.g., $\rotsSampling_2$,
    \item pick $\exponentA$ as large as possible $\exponentA\approx2/3$, and
    \item choose $\regParam$ yielding a proper amount of points with non-zero coefficient, e.g., having $\rotsSampling_2$ where  $\|\rotsDensCoef^1_\imgInd\|_0$ is between 7 and 20 on average. 
\end{enumerate}
 The latter will be achieved through a proper choice of $\dummyIntegerB_0$ and the use of \cref{prop:choice-regulariser-lifting} in the following section.

\subsection{Joint 3D map reconstruction and rotation estimation}
\label{sec:join-volume-rots-experiment}


In this experiment we will test the performance of ESL (\cref{alg:l2-scaled-lifting}) in the alternating minimisation scheme in \cref{alg:joint-refinement-lifting} and see to what extend \cref{prop:choice-regulariser-lifting} allows us to pick the amount of non-zero coefficients. We will be concerned with more straightforward questions than the ones in the previous section:
\begin{enumerate}
    \item \revA{Is} lifting-based rotation estimation in a joint scheme with 3D map reconstruction \revA{beneficial in practice compared to popular software}?
    \item \revA{Is the asymptotic behaviour in \cref{prop:choice-regulariser-lifting} observed throughout the iterative refinement scheme so that} the number of non-zero coefficients \revA{can be controlled (through \cref{thm:new-num-non-zero-coeffs})}?
\end{enumerate}



\paragraph{Metrics}
The following metrics will be used to investigate the questions above.
\begin{enumerate}
    \item For testing \cref{alg:joint-refinement-lifting}: the progression of the distribution of $\distance_{\mathrm{SO}(3)}(\rot^\iterInd_{\imgInd},\rot^{\mathrm{GT}}_{\imgInd})$ \revA{and the error distribution of $\distance_{\mathrm{SO}(3)}(\rot^{\text{RELION}}_{\imgInd},\rot^{\mathrm{GT}}_{\imgInd})$, where $\rot^{\text{RELION}}_{\imgInd}$ is the rotation with largest weight for image $\imgInd$ after a full run of RELION}.
    \item For testing \cref{prop:choice-regulariser-lifting}: the progression of the distribution of $\|\rotsDensCoef^\iterInd_\imgInd\|_0$.
\end{enumerate}
For the sake of completeness, a progression of the 3D maps will also be provided for visual inspection. 





\paragraph{Algorithm settings}
For \cref{alg:joint-refinement-lifting} we use $\rotsSampling:=\rotsSampling_2$ from the previous experiment, $\vol^0:= \vol^{\mathrm{INIT}}$, $\noiseLevelParam,\volRegParam_1,\volRegParam_2$ as in \cref{eq:vol-param-choice}, $\dummyIntegerB_0=15$ and $\exponentA=0.66$. In total, 10 iterations of \cref{alg:joint-refinement-lifting} will be done with 20 Riemannian gradient descent iterations in the inner loop, \revA{both of which are empirically sufficient to converge}.

\revA{For RELION we run the software in two ways: RELION \emph{with global search only} and \emph{unrestricted} RELION, where the software does not only do a global search, but also uses automatic mesh refinement and a local search stage. The first setting is chosen to stay as close to \cref{alg:joint-refinement-lifting} as possible. In particular, this solver uses a sampling set of size $|\rotsSampling_{\text{RELION}}| =294\,912 > \numRots$ rotations. The reason for a different sampling set is that RELION does not allow for a custom set, i.e., we could not use $\rotsSampling$, but had to take at least a larger one for fair comparison. The second setting is chosen to see how the proposed method compares to an actual state-of-the-art method, that is numerically optimised for over a decade. Before moving on to the results, it should be noted that it was impossible to turn off the marginalisation over in-plane shifts. Although in both cases the solver predicted these correctly, i.e., close to zero, it should be mentioned that  RELION is solving a harder problem.}



\begin{figure}[h!]
    \centering
    \begin{subfigure}{\dimexpr0.30\textwidth+20pt\relax}
    \makebox[20pt]{\raisebox{40pt}{\rotatebox[origin=c]{90}{$\distance_{\mathrm{SO}(3)}(\rot^\iterInd_{\imgInd},\rot^{\mathrm{GT}}_{\imgInd})$}}}%
    \includegraphics[width=\dimexpr\linewidth-20pt\relax]{"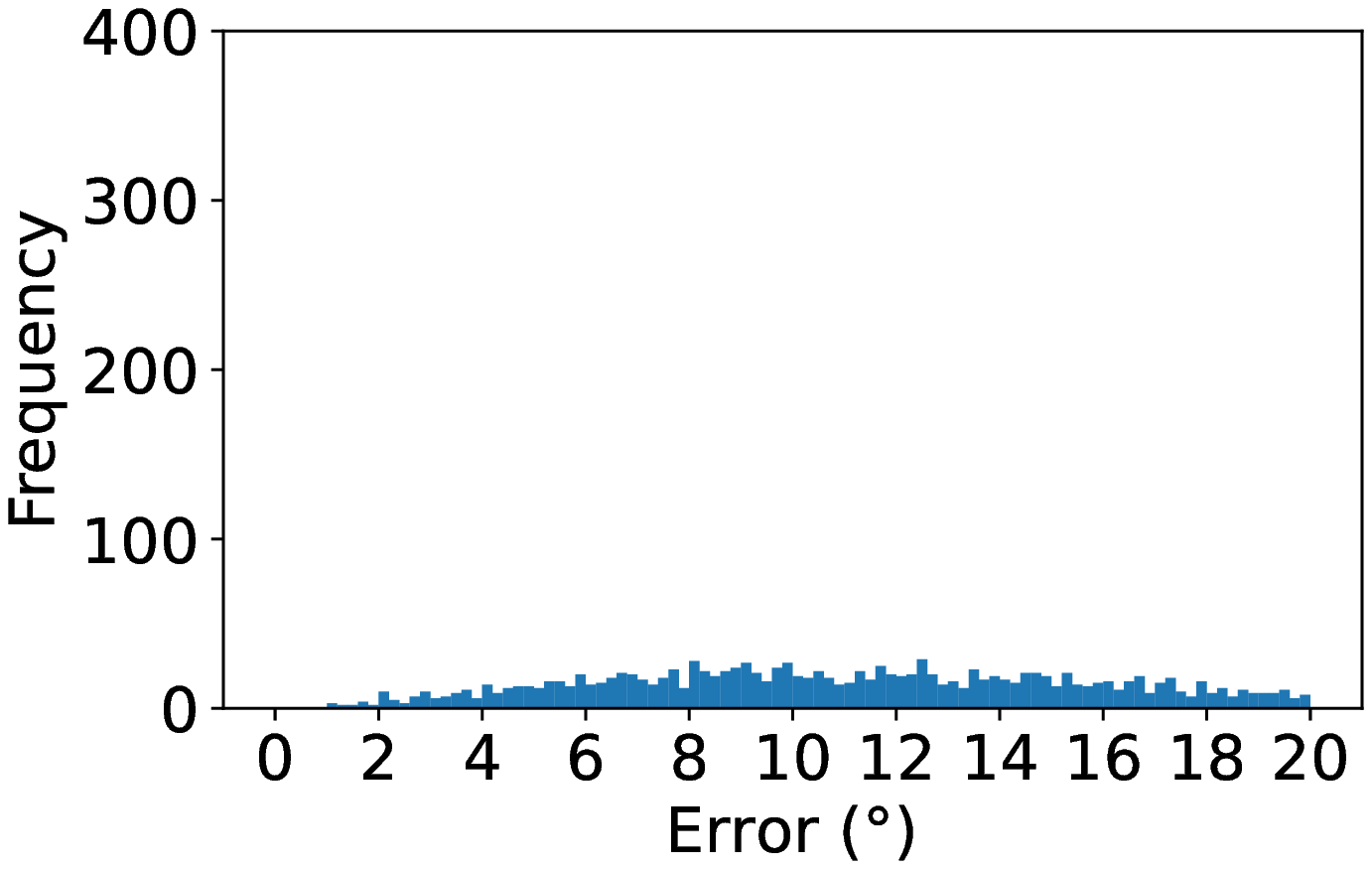"}
    \makebox[20pt]{\raisebox{40pt}{\rotatebox[origin=c]{90}{ $\|\rotsDensCoef^\iterInd_\imgInd\|_0$}}}%
    \includegraphics[width=\dimexpr\linewidth-20pt\relax]{"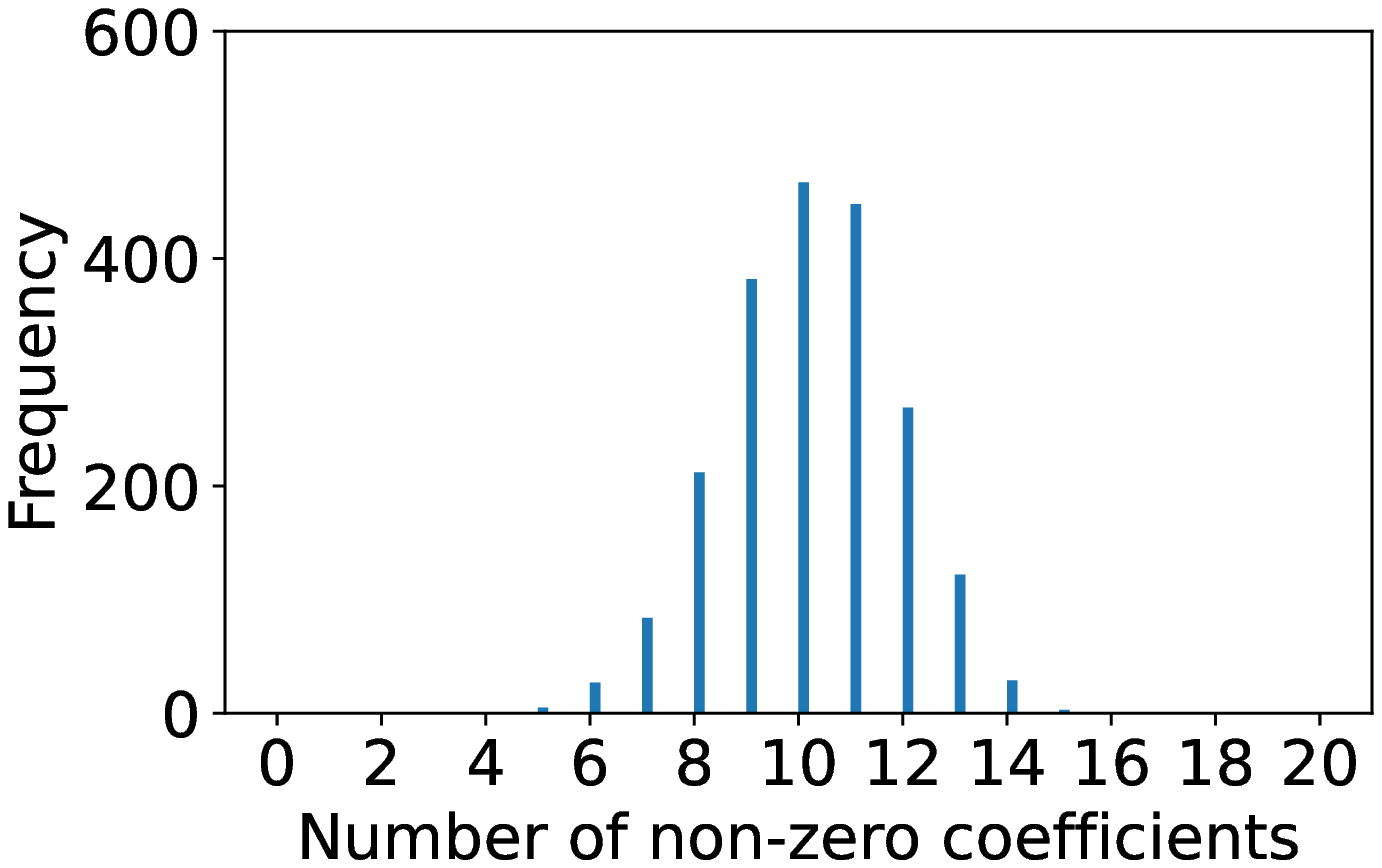"}
    \makebox[20pt]{\raisebox{40pt}{\rotatebox[origin=c]{90}{$\vol^\iterInd$}}}%
    \includegraphics[width=\dimexpr\linewidth-20pt\relax]{"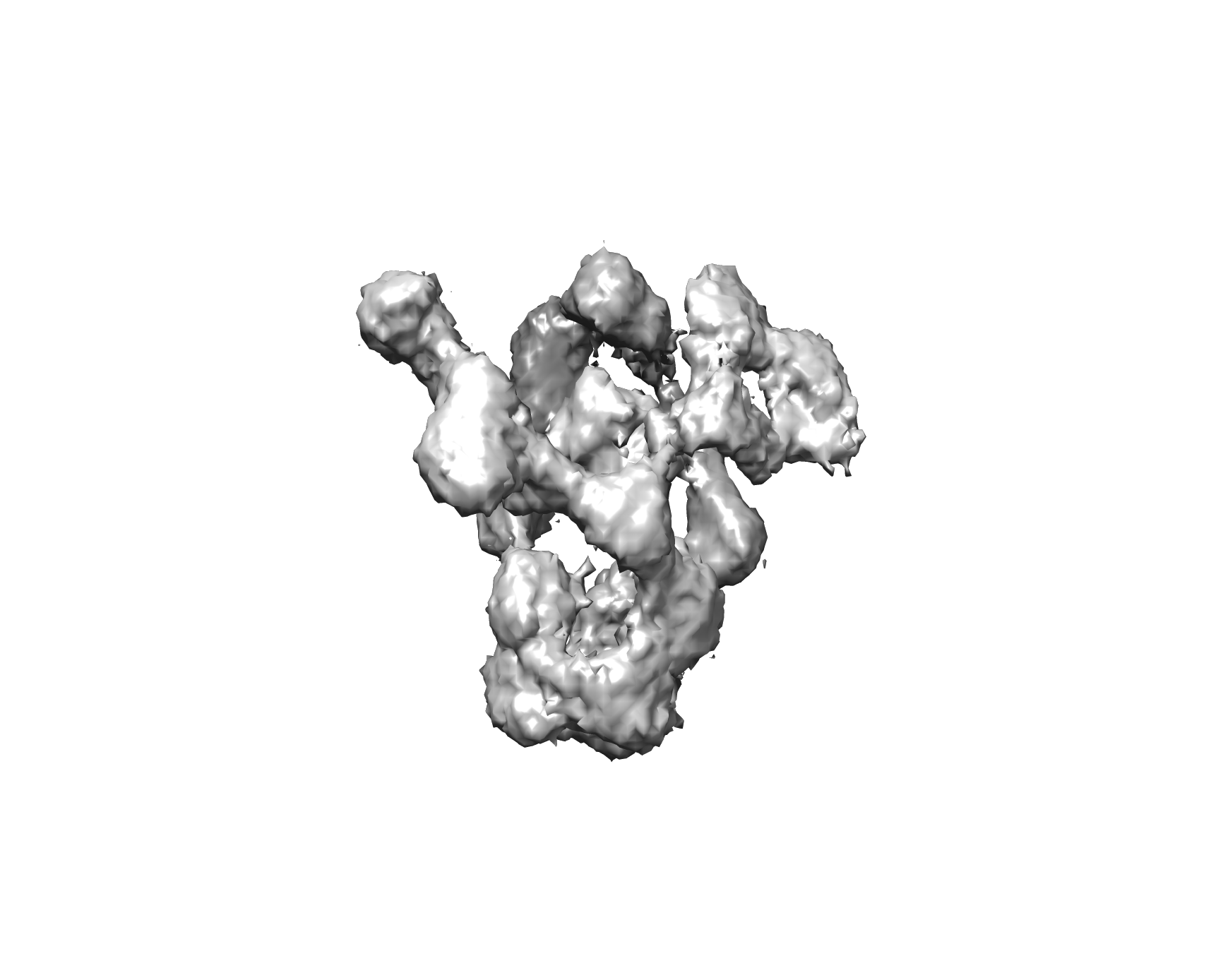"}
    \caption{$\iterInd=1$}
    \end{subfigure}
    \hfill
    \begin{subfigure}{0.30\textwidth}
    \includegraphics[width=\linewidth]{"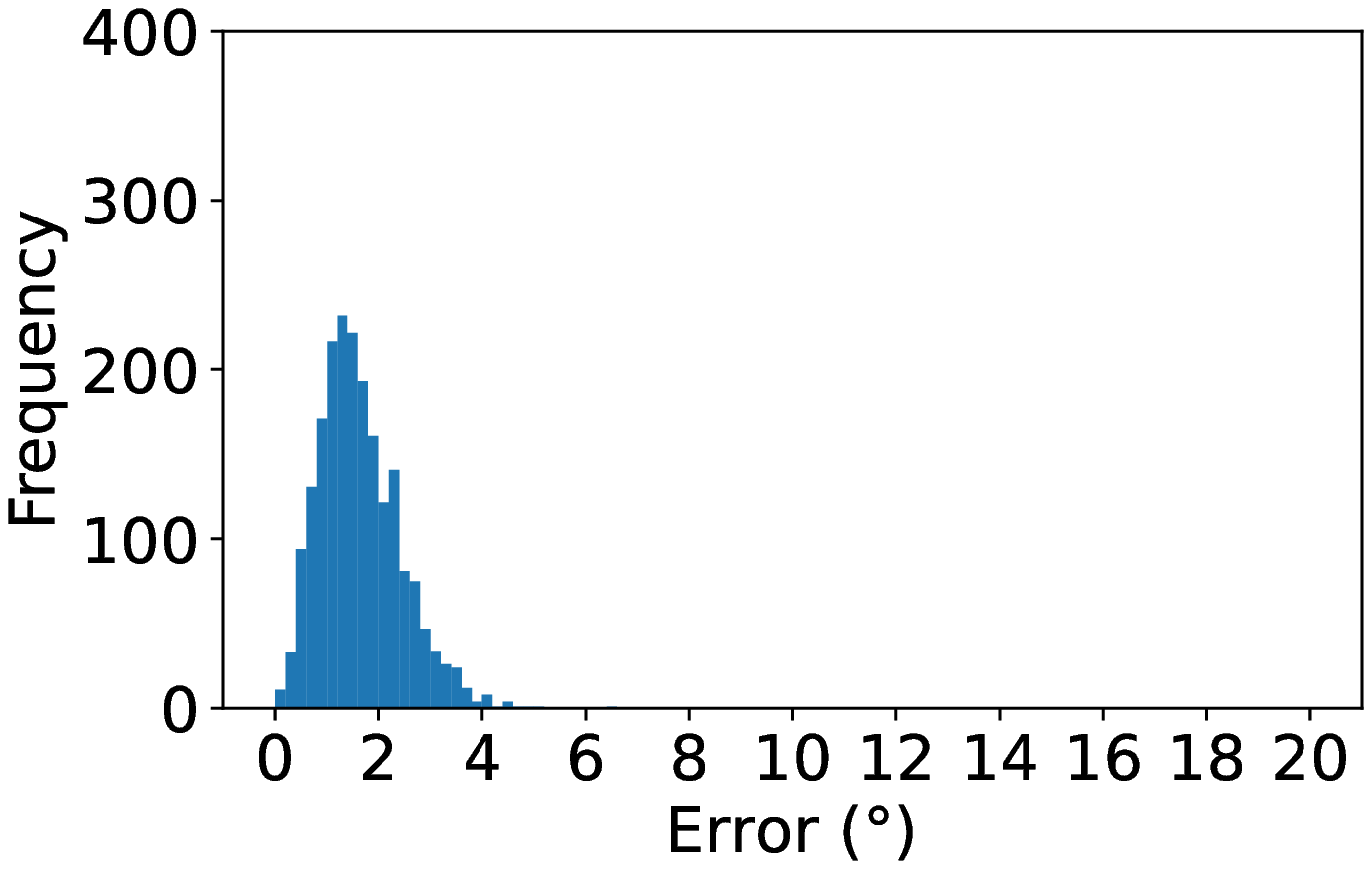"}
    \includegraphics[width=\linewidth]{"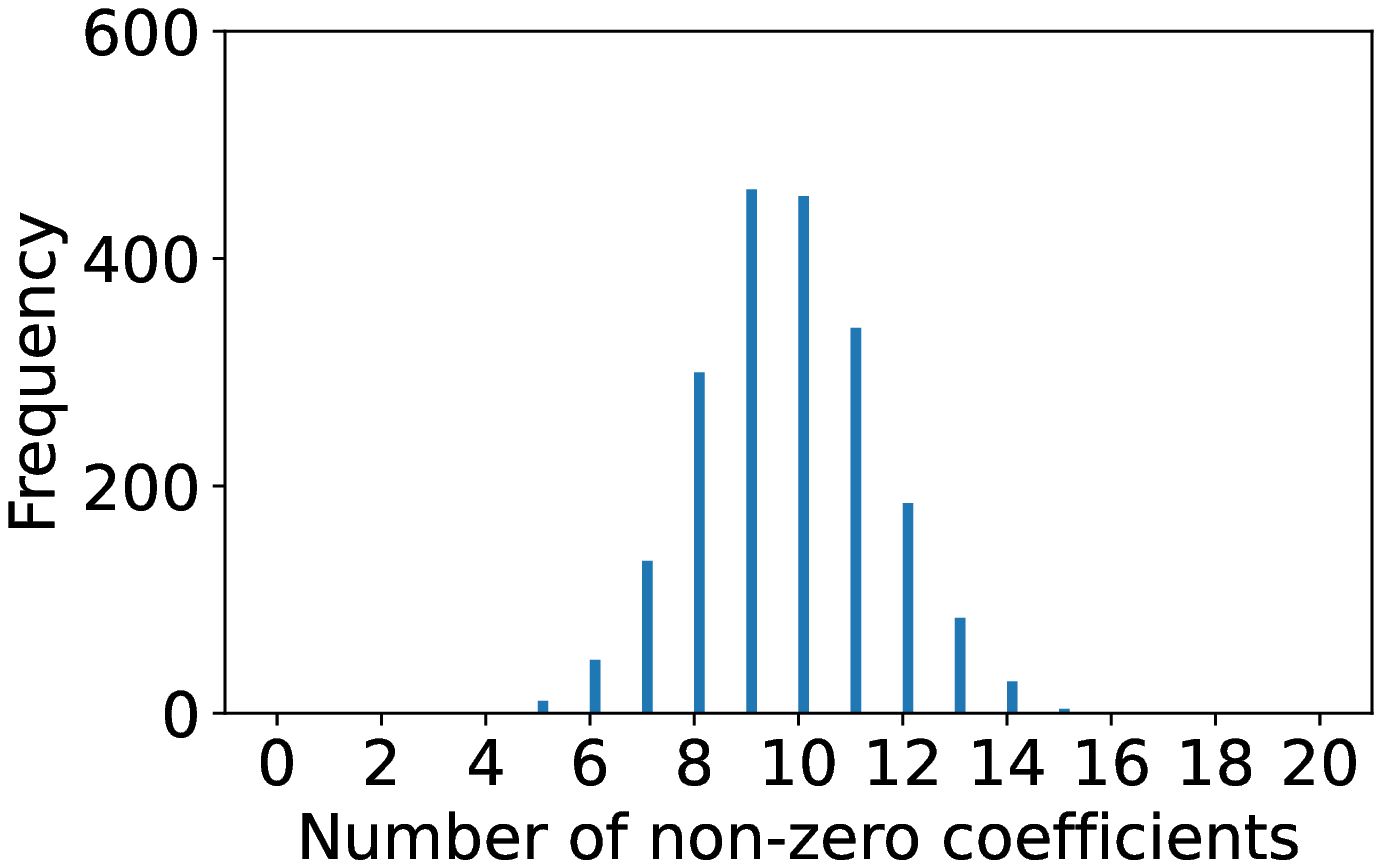"}
    \includegraphics[width=\linewidth]{"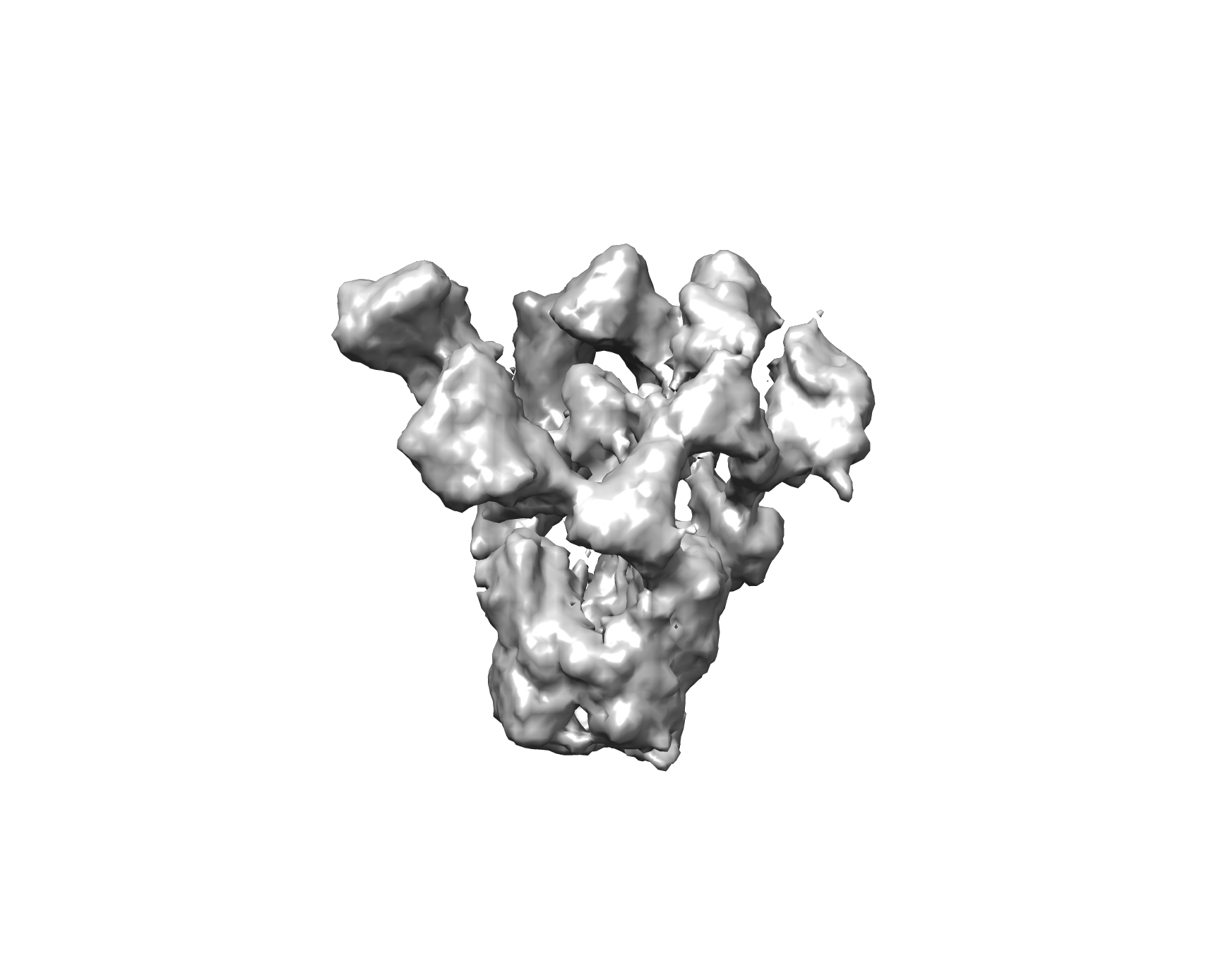"}
    \caption{$\iterInd=5$}
    \end{subfigure}
    \hfill
    \begin{subfigure}{0.30\textwidth}
    \includegraphics[width=\linewidth]{"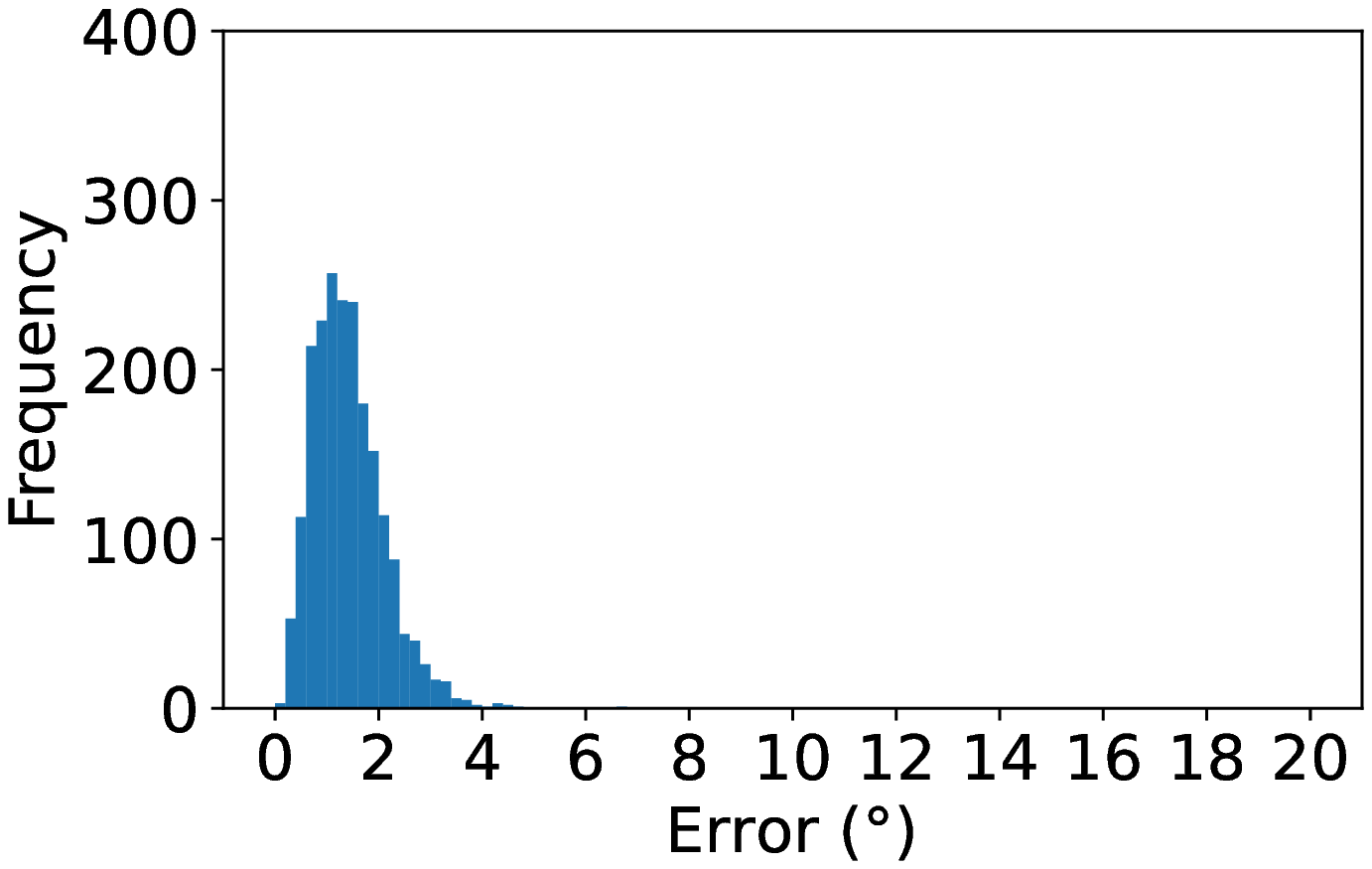"}
    \includegraphics[width=\linewidth]{"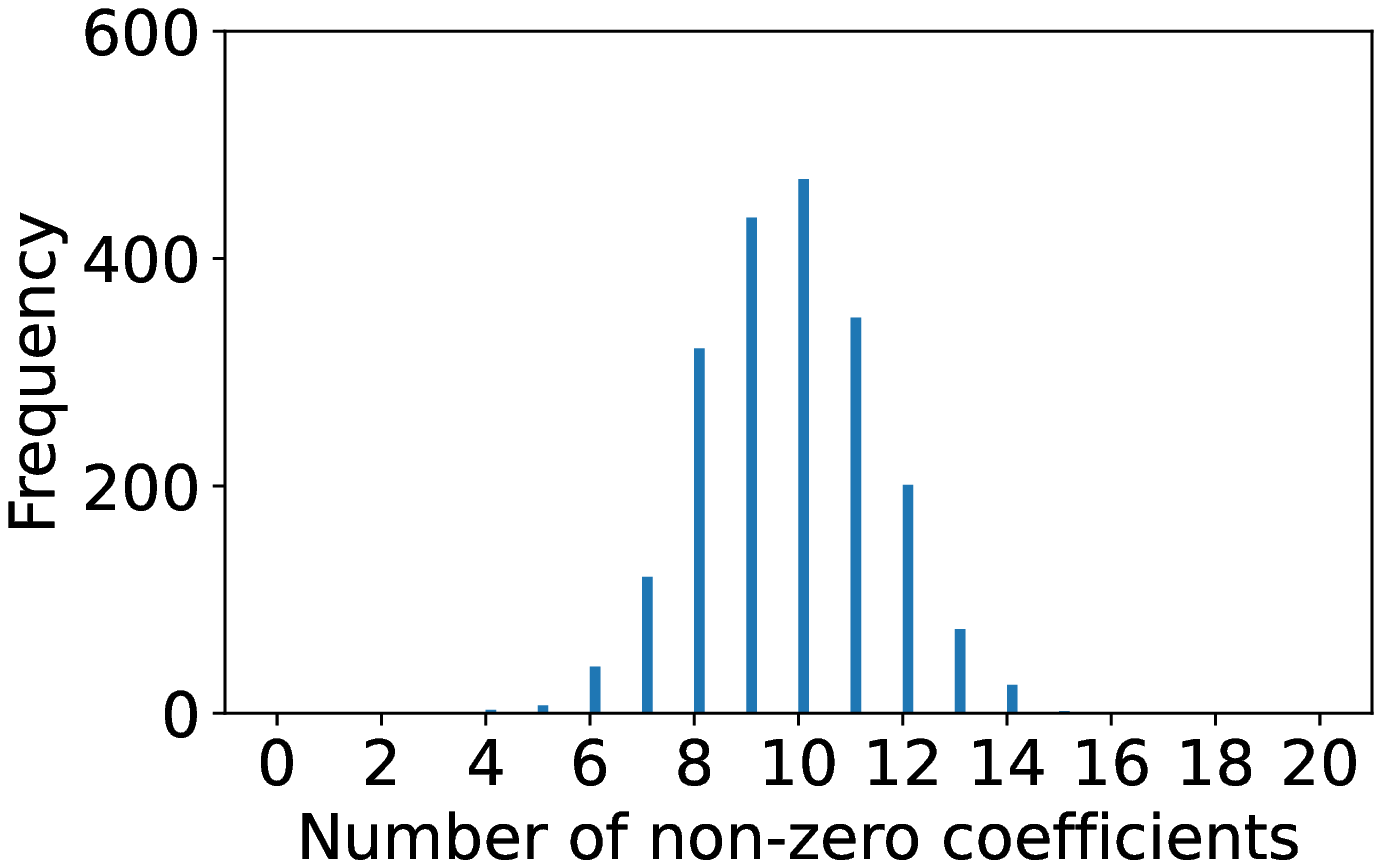"}
    \includegraphics[width=\linewidth]{"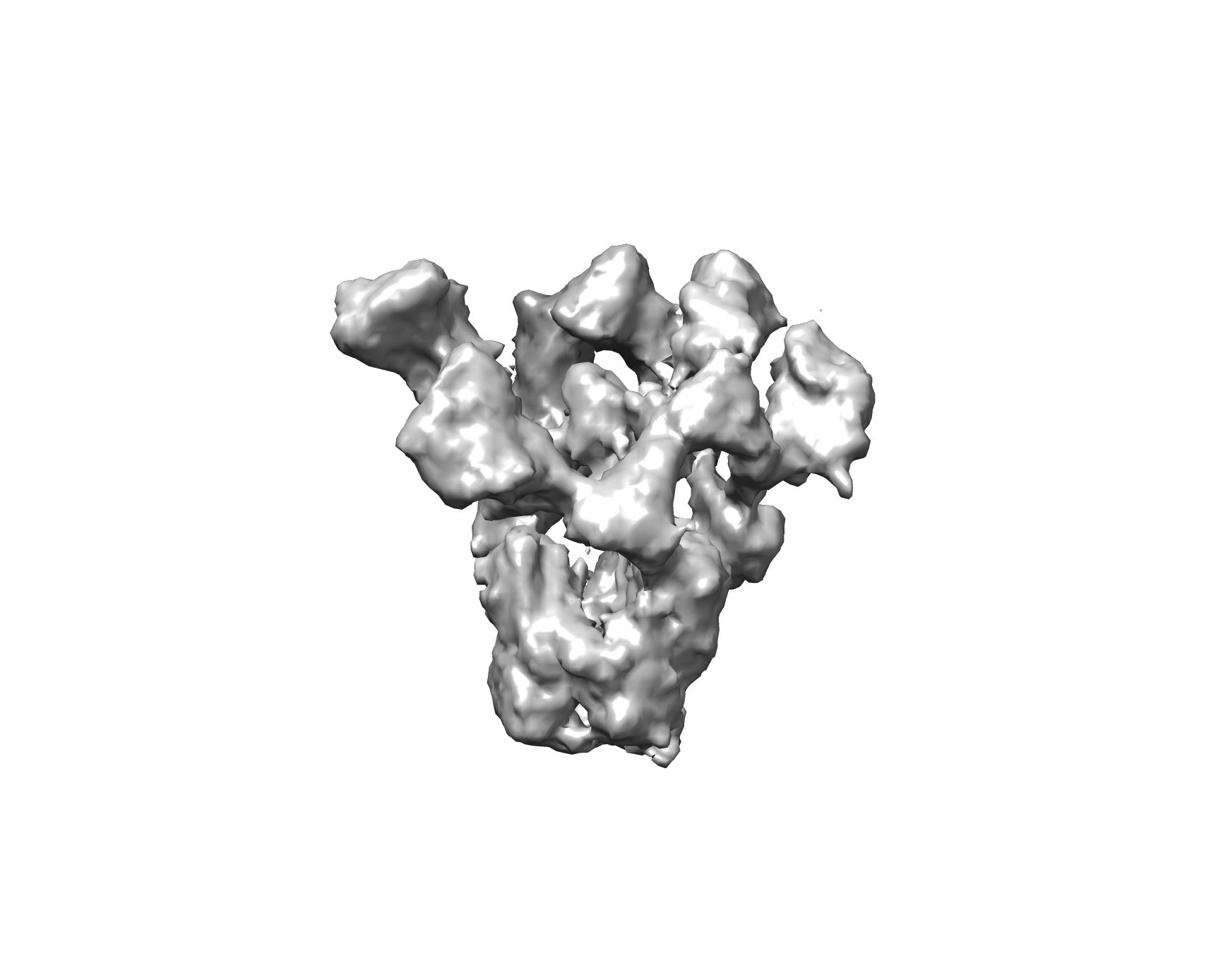"}
    \caption{$\iterInd=10$}
    \end{subfigure}
    \caption{The progression of the distribution of the rotation estimation error, the number of non-zero coefficients and the \revA{iso-surfaces of the} 3D map. The rotation estimation scheme has more or less converged after 5 iterations and the number of non-zero coefficients are in line with the theory.}
    \label{fig:full-algo-progression}
\end{figure}

\begin{figure}[h!]
    \centering
    \begin{subfigure}{\dimexpr0.30\textwidth+20pt\relax}
    \makebox[20pt]{\raisebox{40pt}{\rotatebox[origin=c]{90}{ $\distance_{\mathrm{SO}(3)}(\rot^{\text{SOL}}_{\imgInd},\rot^{\mathrm{GT}}_{\imgInd})$}}}%
    \includegraphics[width=\dimexpr\linewidth-20pt\relax]{"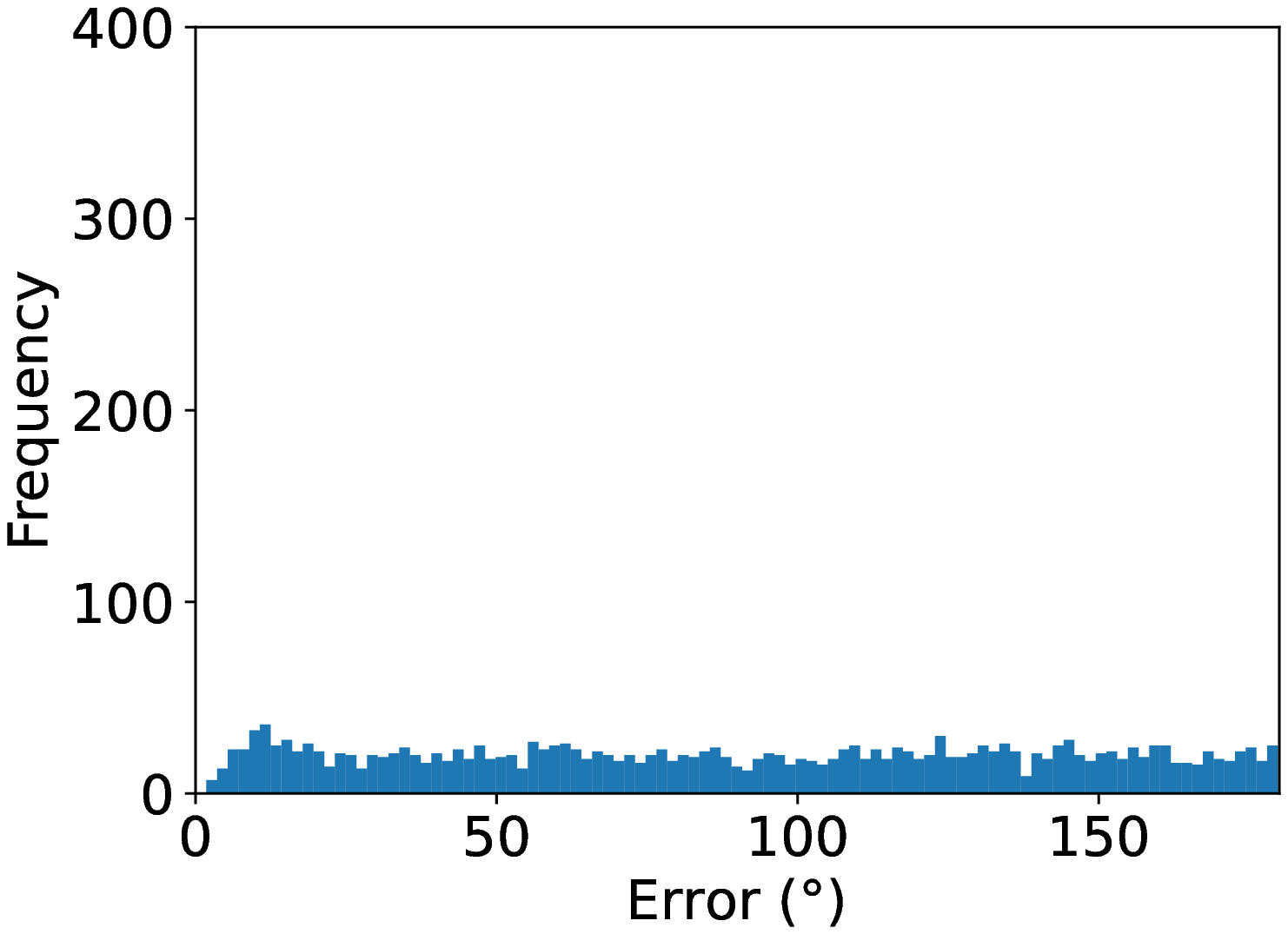"}
    \makebox[20pt]{\raisebox{40pt}{\rotatebox[origin=c]{90}{$\vol^{\text{SOL}}$}}}%
    \includegraphics[width=\dimexpr\linewidth-20pt\relax]{"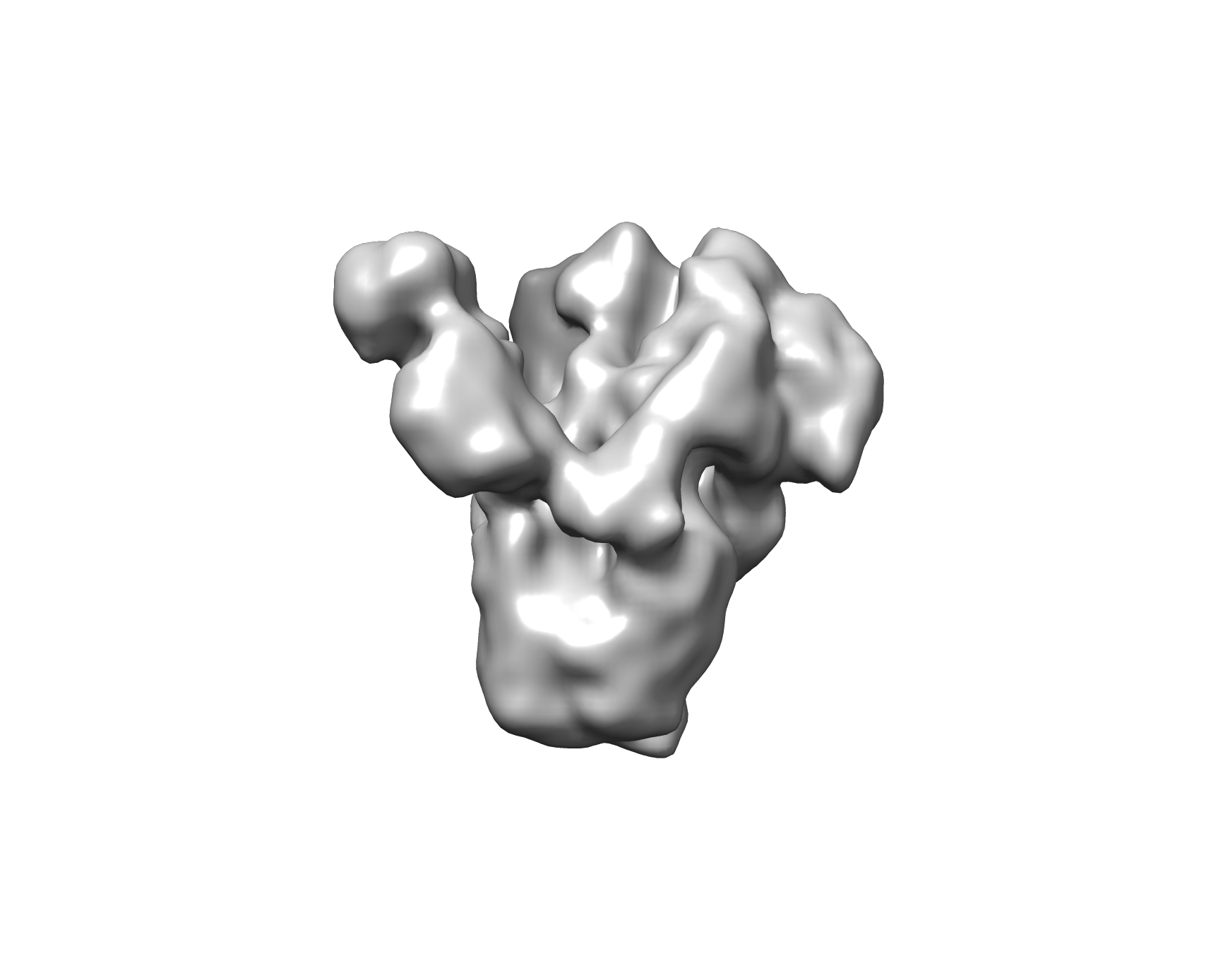"}
    \caption{RELION\\ \hspace{0.2cm} (global search only)}
    \end{subfigure}
    \hfill
    \begin{subfigure}{0.30\textwidth}
    \includegraphics[width=\linewidth]{"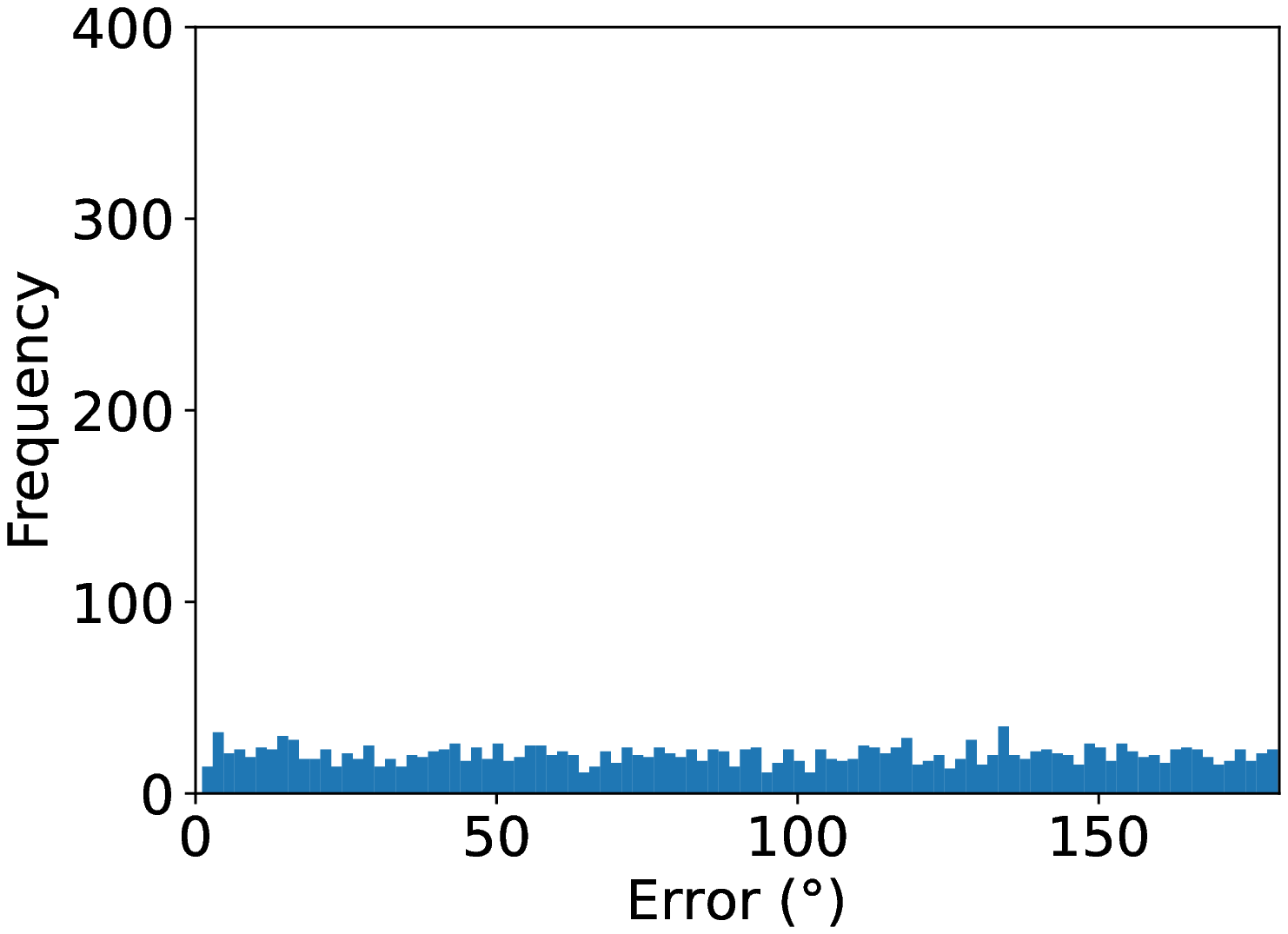"}
    \includegraphics[width=\linewidth]{"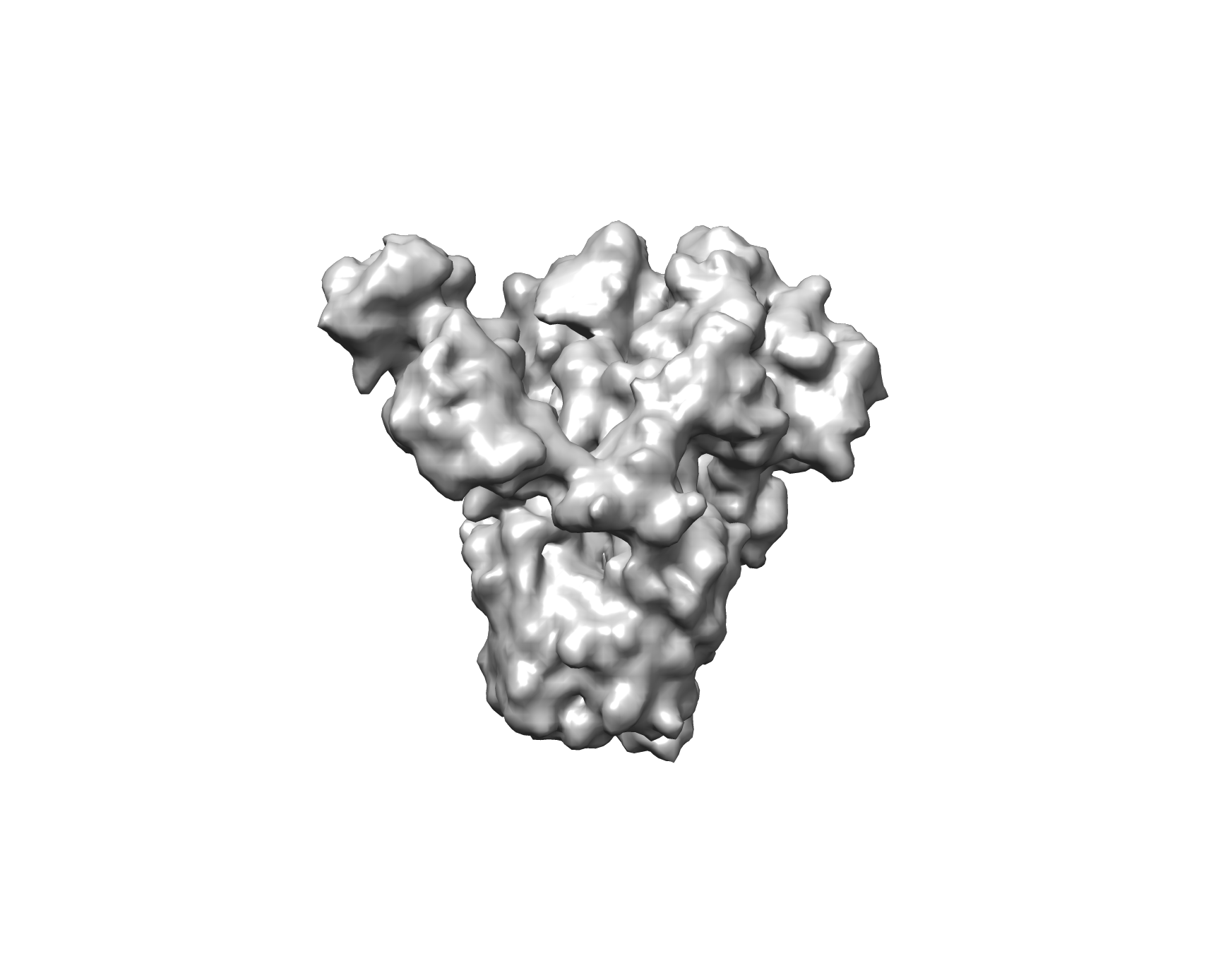"}
    \caption{RELION\\ (unrestricted)}
    \end{subfigure}
    \hfill
    \begin{subfigure}{0.30\textwidth}
    \includegraphics[width=\linewidth]{"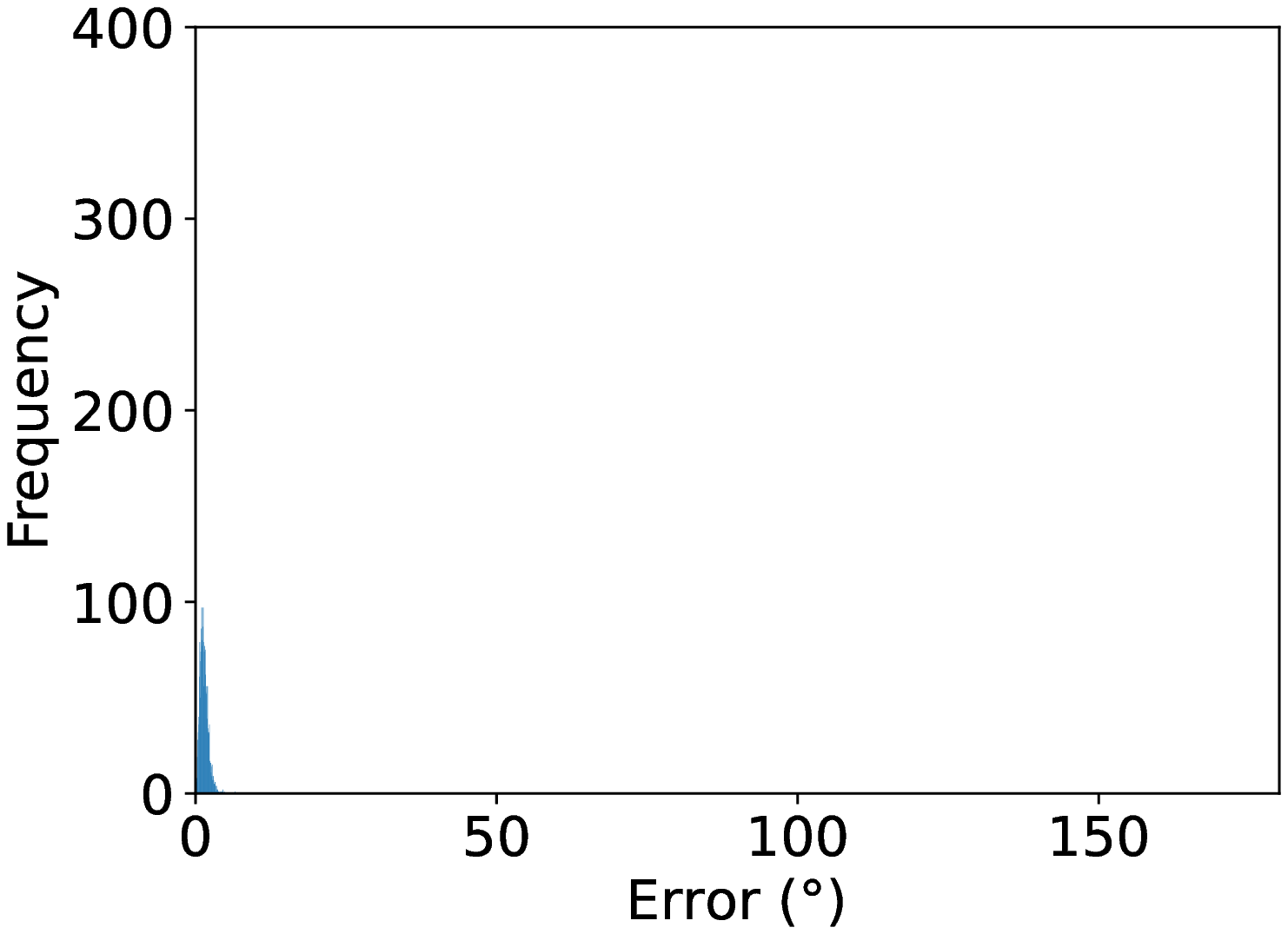"}
    \includegraphics[width=\linewidth]{"results/joint-reconstruction/expC3_22-04-06_23-59-48_SNR16_N2048_J15_r2_s10_i10/iter10.png"}
    \caption{ESL\\ (Ours)}
    \end{subfigure}
    \caption{\revA{The distributions of the rotation estimation error and the 3D maps as the result of different solvers. Both the restricted and unrestricted RELION runs do not estimate the rotations too well when compared to the proposed method. However, other factors such as RELION's Fourier-shell-correlation-regularisation and the fact that these methods marginalise over a large set of rotations -- rather than hard assignment of the above rotations --, are the reason why the 3D maps still look good.}}
    \label{fig:RELION_runs}
\end{figure}

\begin{figure}[h!]
    \centering
    \begin{subfigure}{\dimexpr0.30\textwidth+20pt\relax}
    \makebox[20pt]{\raisebox{40pt}{\rotatebox[origin=c]{90}{$\phi$}}}%
    \includegraphics[width=\dimexpr\linewidth-20pt\relax]{"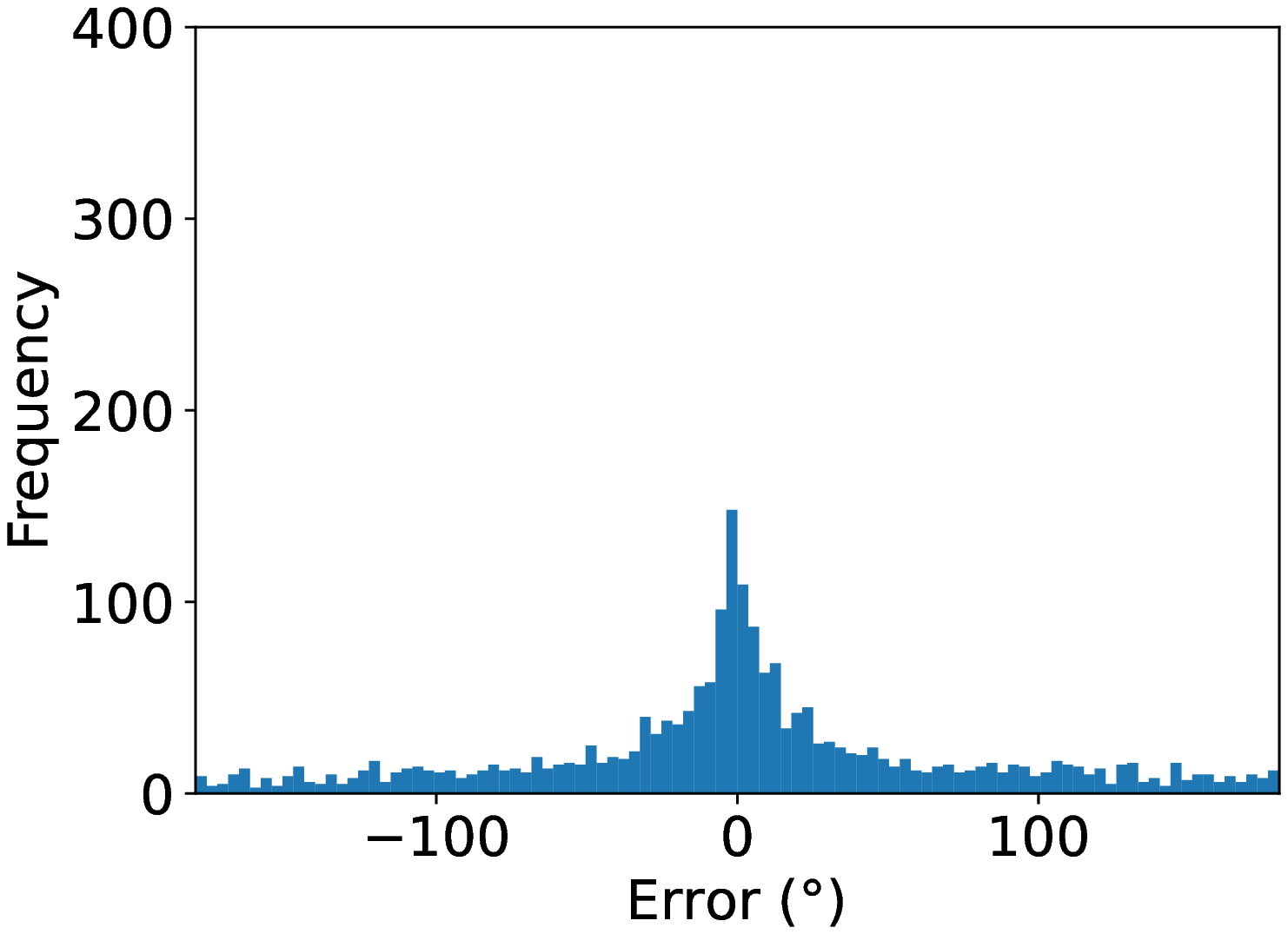"}
    \makebox[20pt]{\raisebox{40pt}{\rotatebox[origin=c]{90}{$\theta$}}}%
    \includegraphics[width=\dimexpr\linewidth-20pt\relax]{"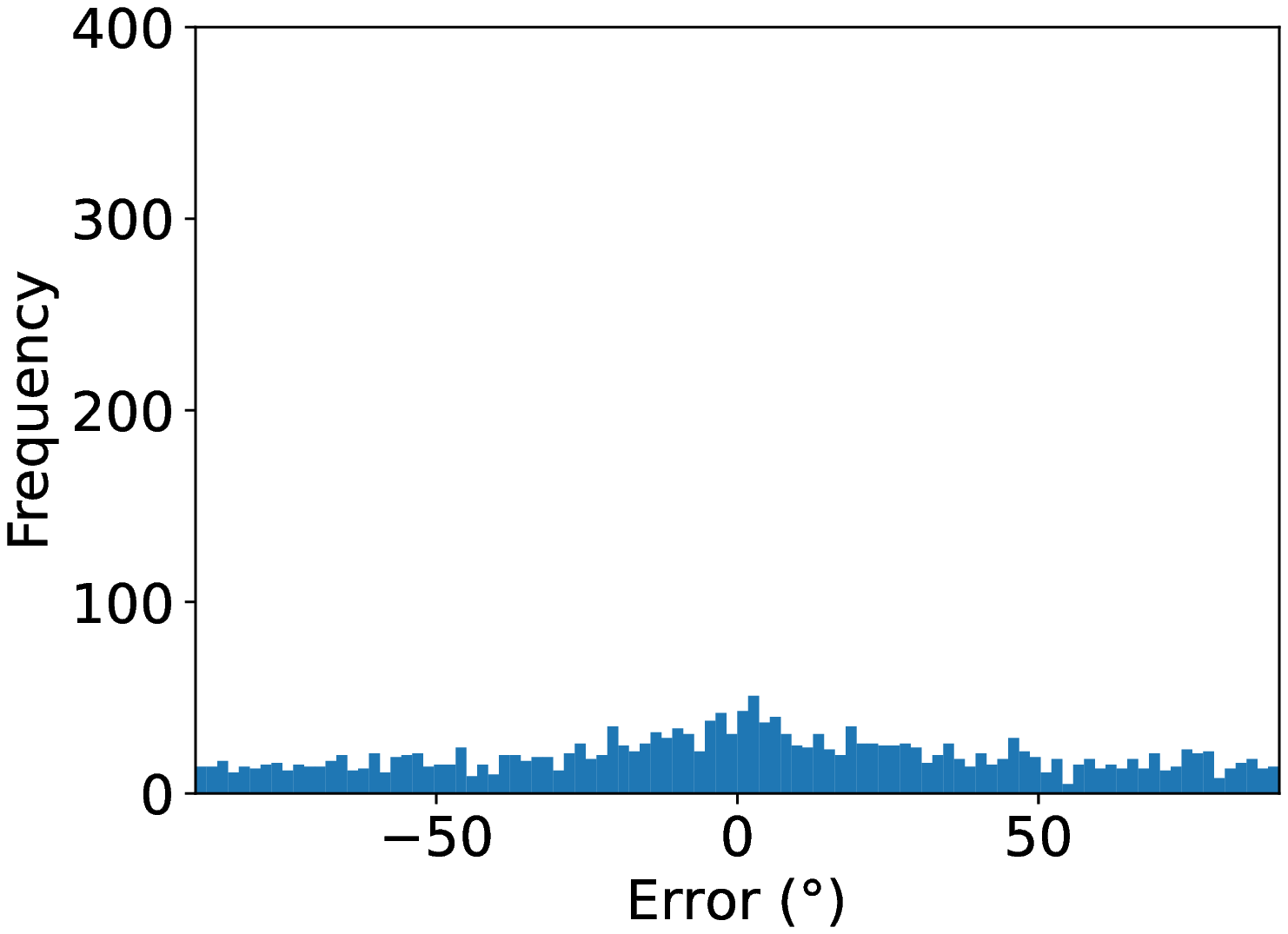"}
    \makebox[20pt]{\raisebox{40pt}{\rotatebox[origin=c]{90}{$\psi$}}}%
    \includegraphics[width=\dimexpr\linewidth-20pt\relax]{"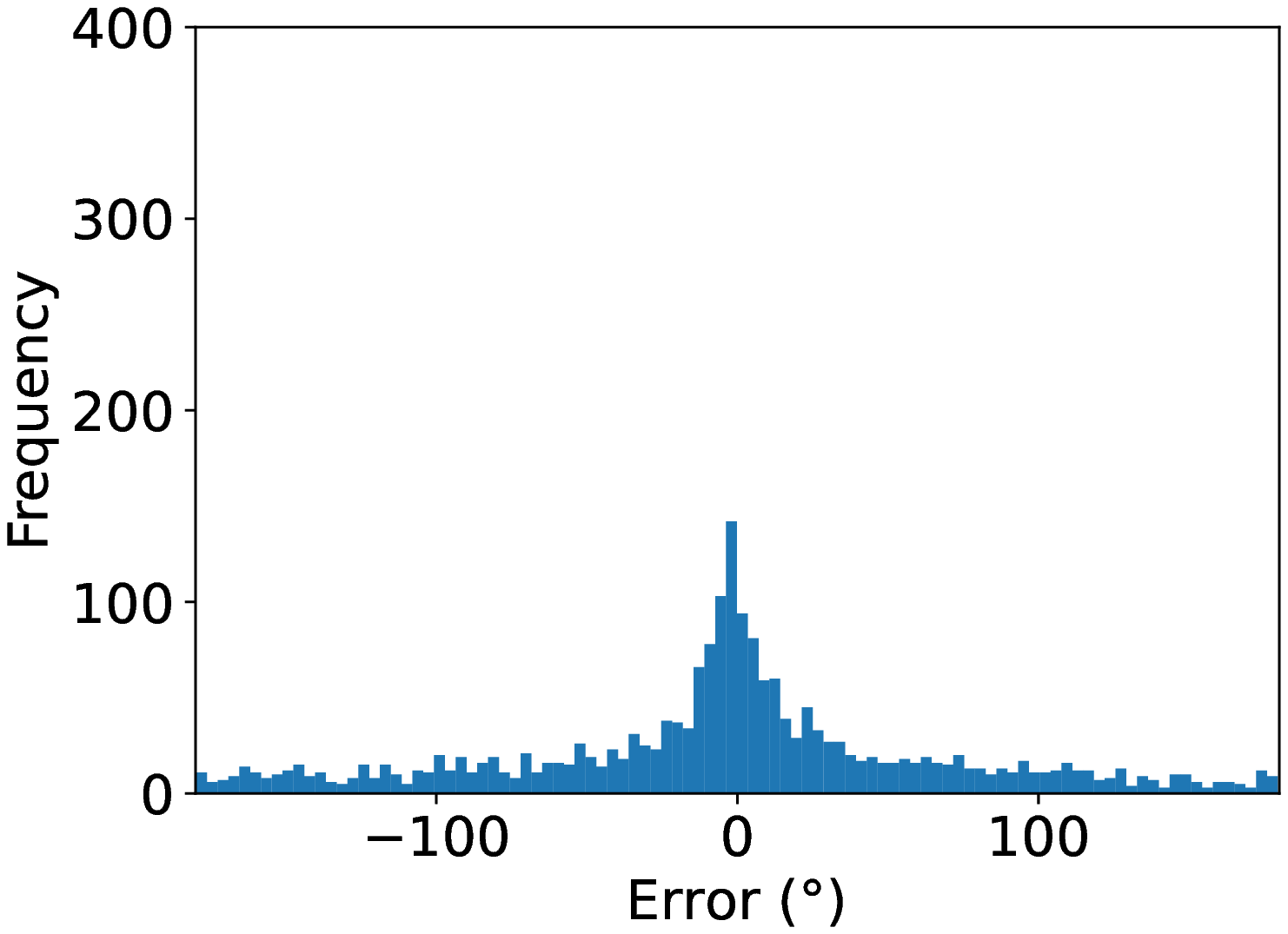"}
    \caption{RELION\\ \hspace{0.2cm} (global search only)}
    \end{subfigure}
    \hfill
    \begin{subfigure}{0.30\textwidth}
    \includegraphics[width=\linewidth]{"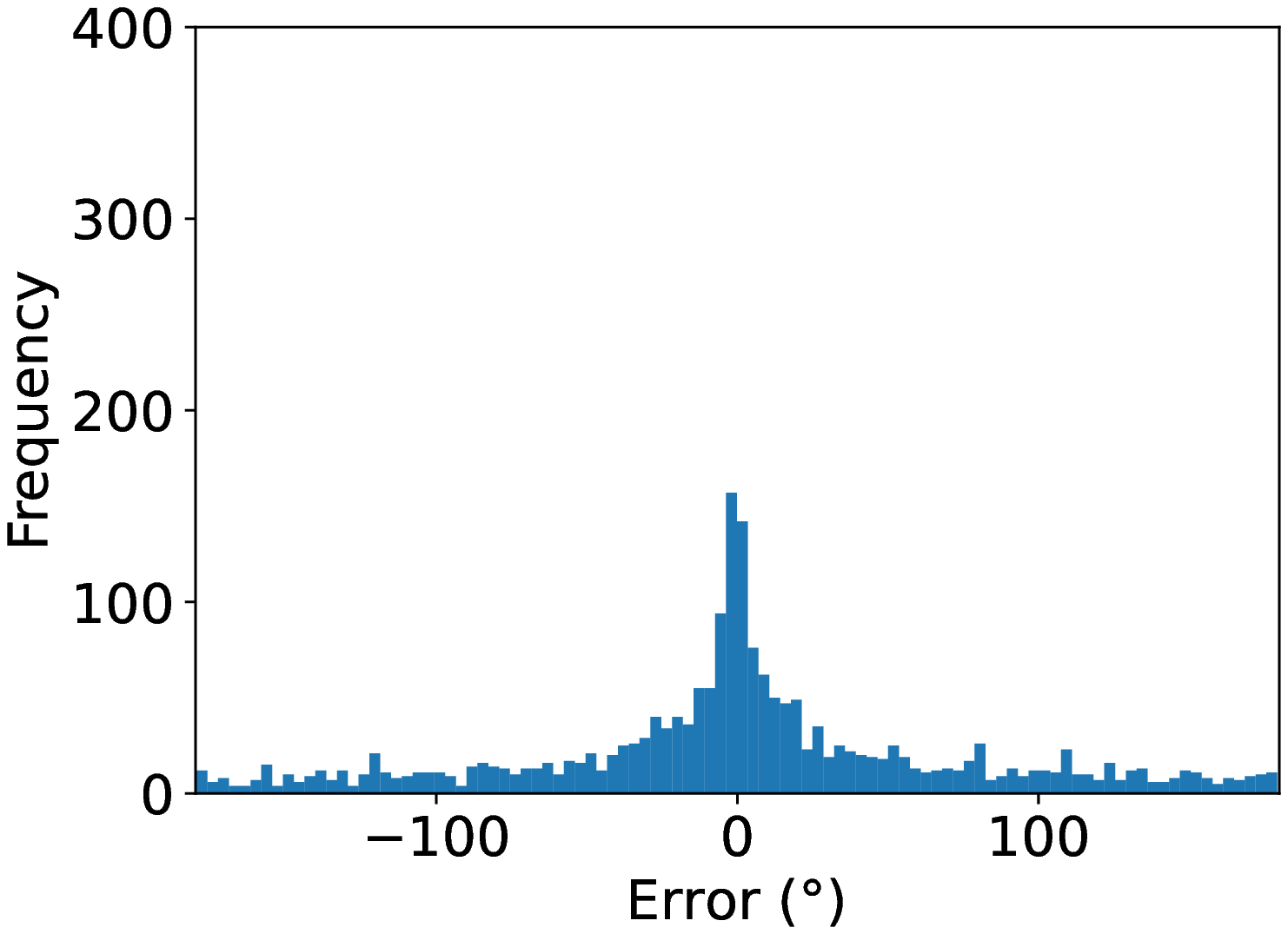"}
    \includegraphics[width=\linewidth]{"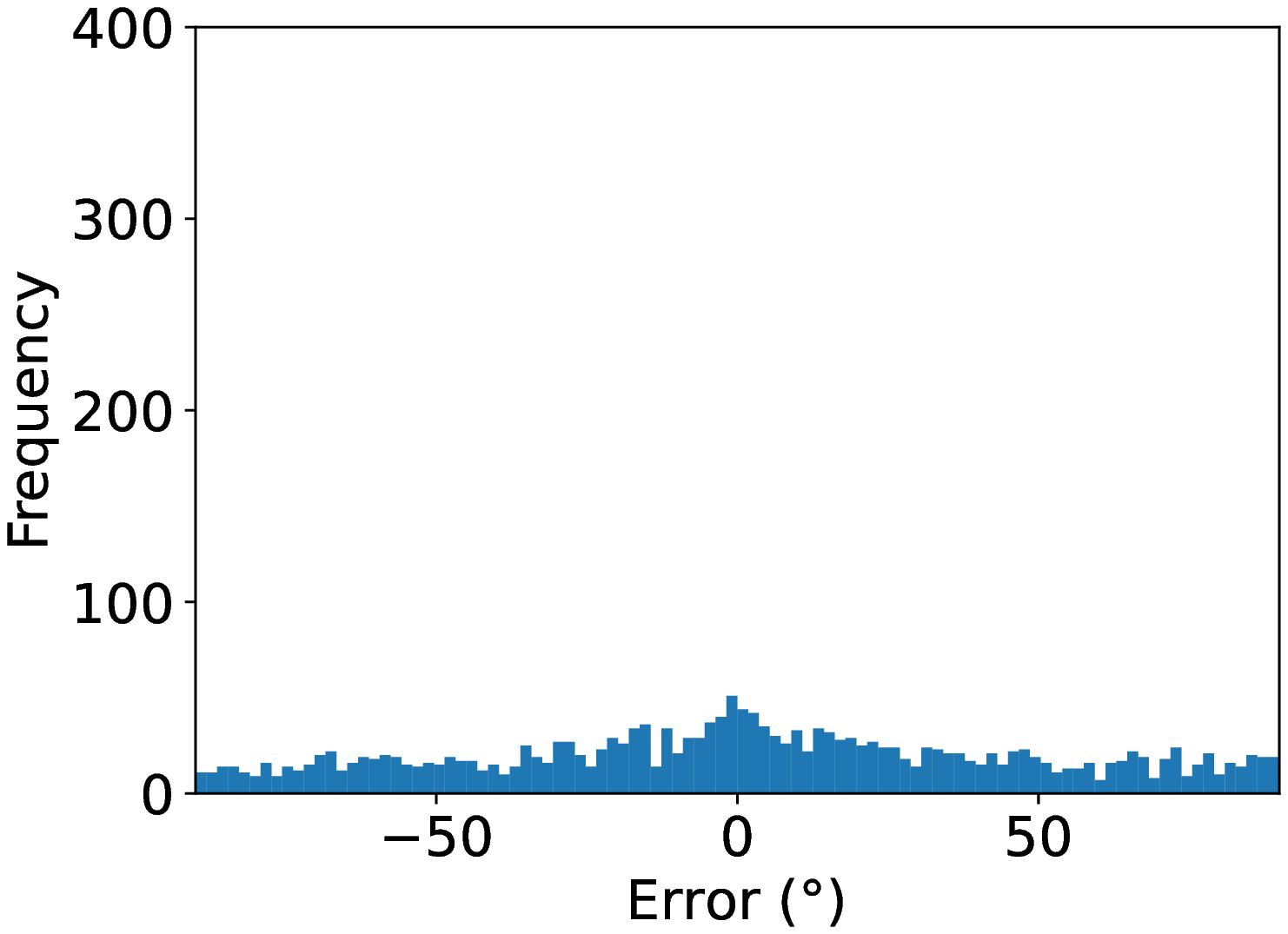"}
    \includegraphics[width=\linewidth]{"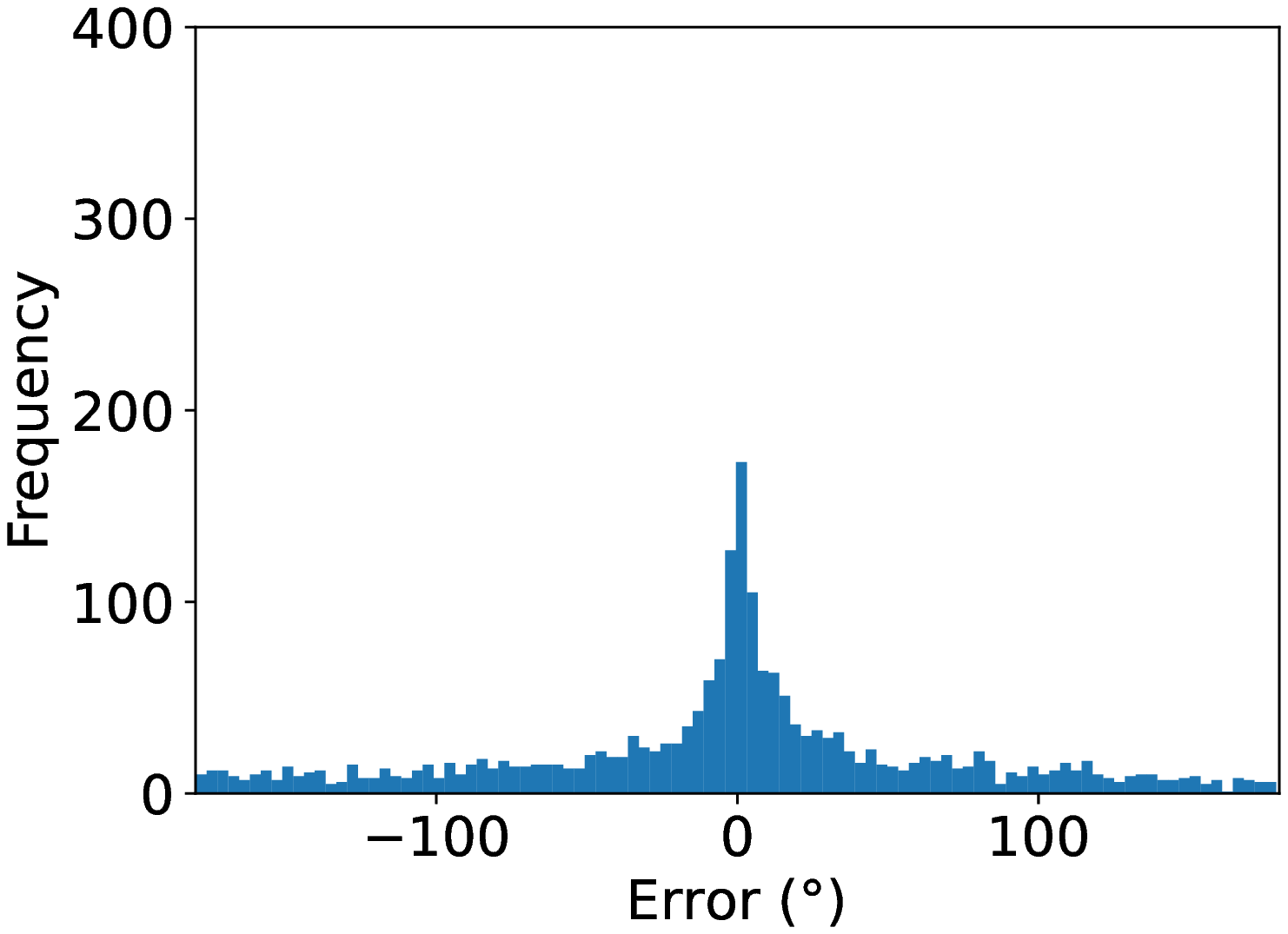"}
    \caption{RELION\\ (unrestricted)}
    \end{subfigure}
    \hfill
    \begin{subfigure}{0.30\textwidth}
    \includegraphics[width=\linewidth]{"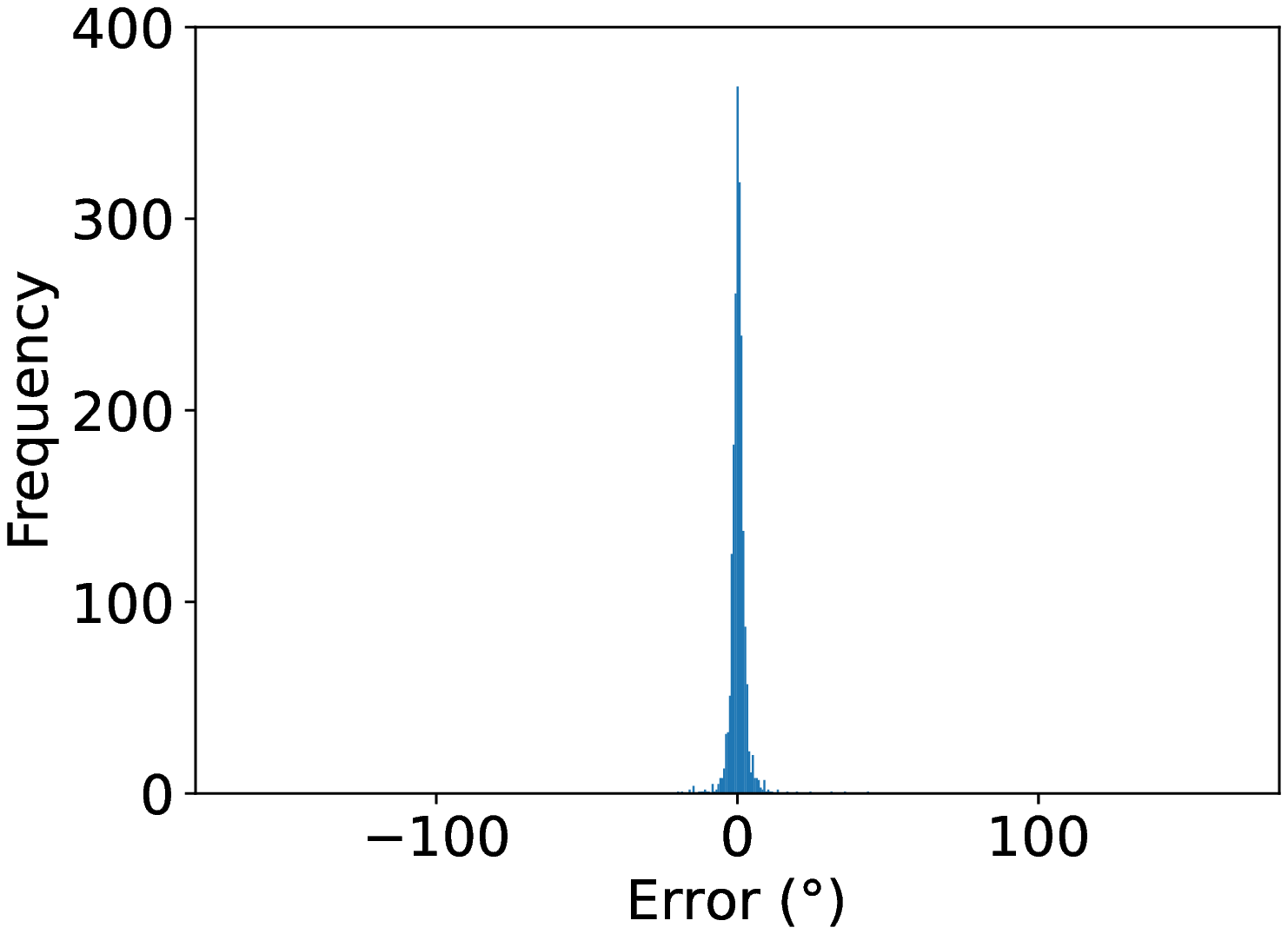"}
    \includegraphics[width=\linewidth]{"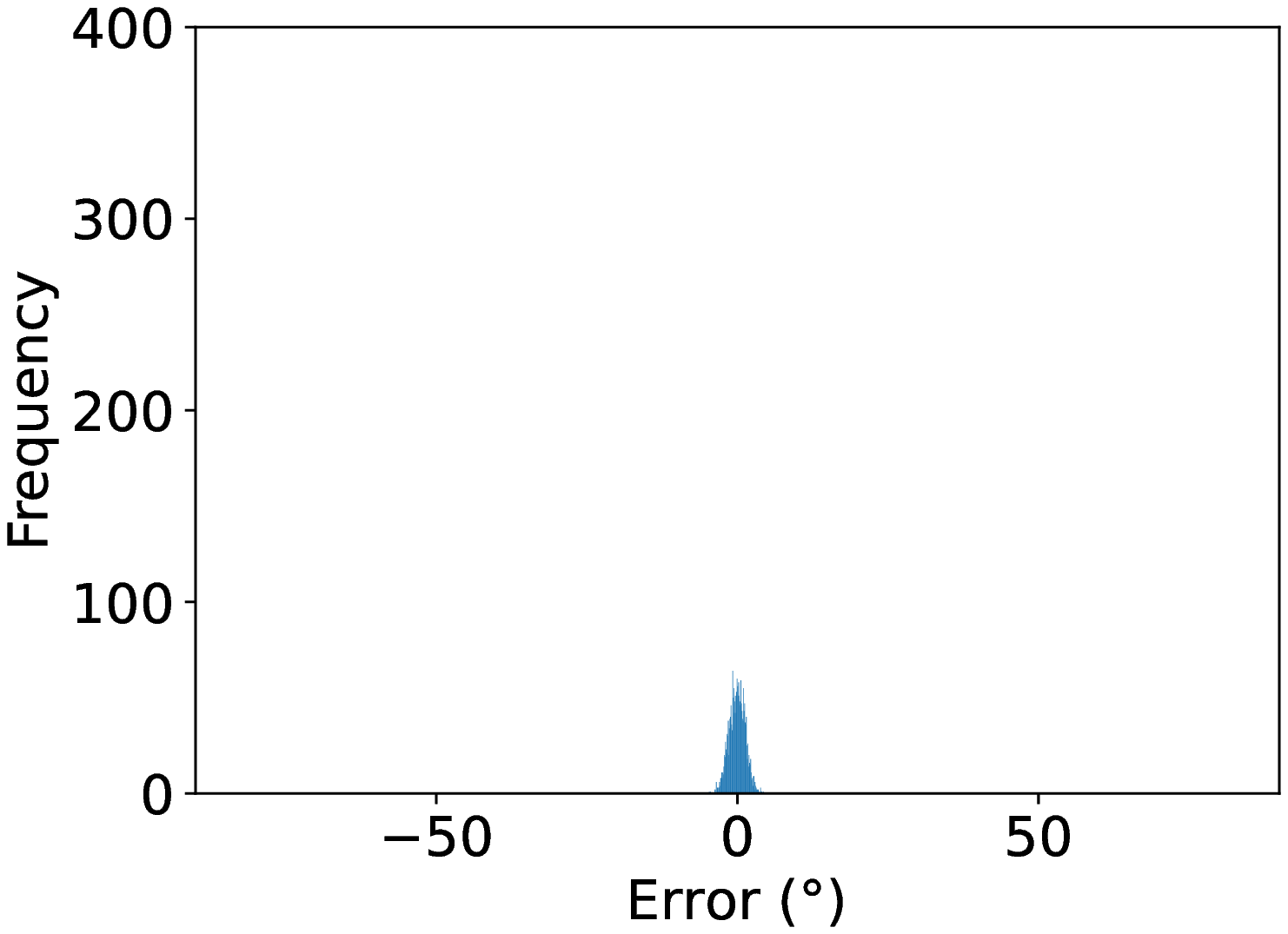"}
    \includegraphics[width=\linewidth]{"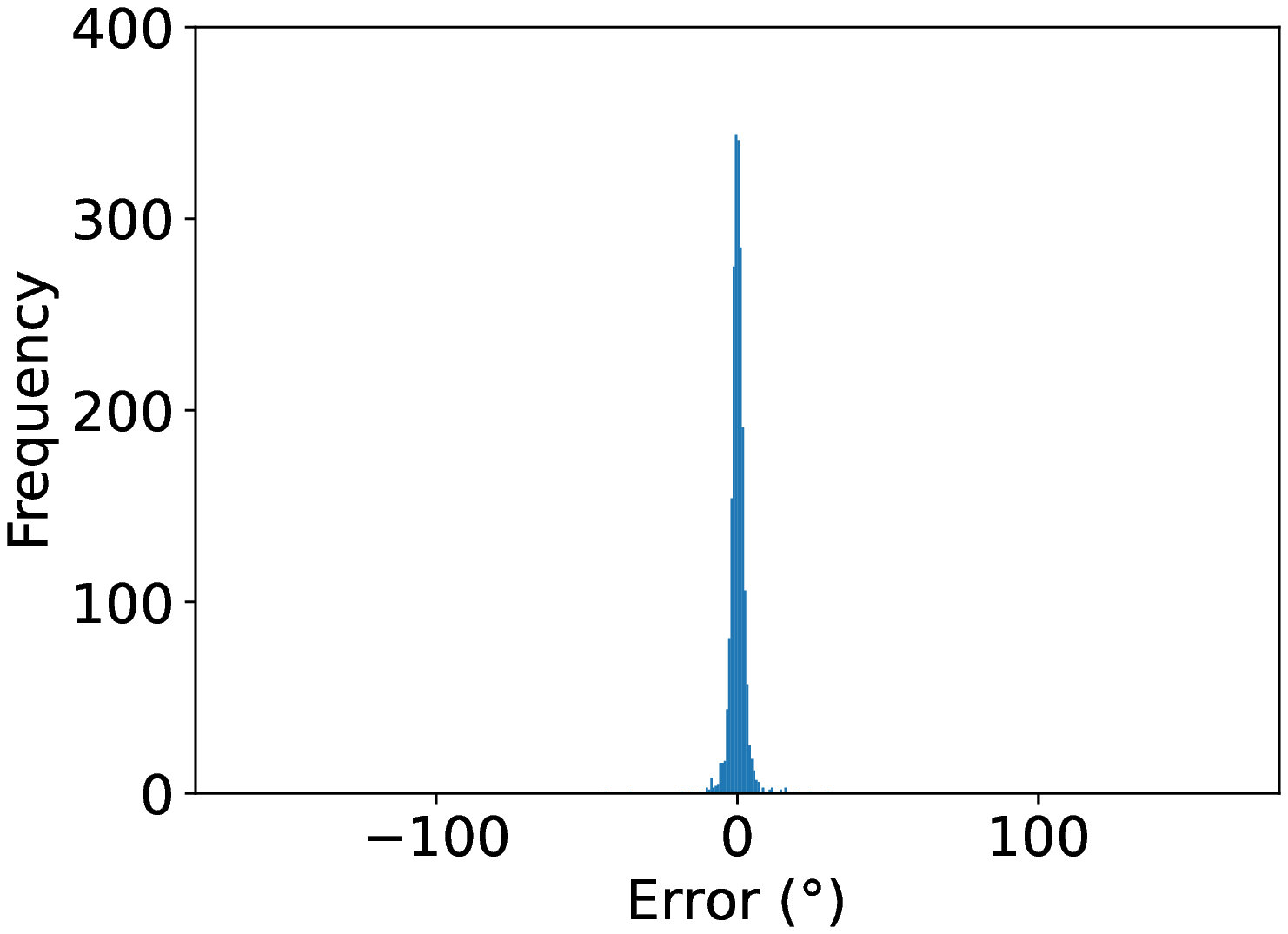"}
    \caption{ESL\\ (Ours)}
    \end{subfigure}
    \caption{\revA{The distributions of the rotation estimation error -- broken down as the rotation errors in the (ZYZ intrinsic) Euler angles $(\phi, \theta, \psi)$ -- as the result of different solvers. Both the restricted and unrestricted RELION runs do not estimate the rotations too well when compared to the proposed method. In particular, the Euler angle $\theta$ is estimated poorly.}}
    \label{fig:RELION_runs_Euler_angles}
\end{figure}

\paragraph{Results}
The progression of the distribution of errors, the amount of non-zero coefficients and the 3D map \revA{resulting from the proposed method} are shown in \cref{fig:full-algo-progression} for iterations 1, 5, and 10. \revA{The distribution of errors and the 3D maps of the final iteration of all three solvers are shown in \cref{fig:RELION_runs}. Additionally, \cref{fig:RELION_runs_Euler_angles} breaks the rotation errors down in the (ZYZ intrinsic) Euler angles $(\phi, \theta, \psi)$.}

For question~1, \cref{fig:full-algo-progression} shows that for all the images the distance between the estimated rotation to the ground truth decreases. After 10 iterations the average distance is $\revA{1.419 \pm 0.686}^{\revA{\circ}}$, which is slightly worse than the results obtained in the previous experiment. This is no surprise as the 3D map $\vol^\iterInd$ has less detail than the ground truth, which makes it arguably harder to get high angular accuracy. \revA{For either of the RELION-based solvers \cref{fig:RELION_runs} shows that the rotation with largest weight does much worse than the proposed method. Upon closer inspection \cref{fig:RELION_runs_Euler_angles} suggests that the Euler angle $\theta$ is the bottleneck. However, we should be careful making hard claims from these results as the RELION-based solvers also marginalised over the in-plane shifts.}
\revA{Nevertheless,} the lifting-based rotation estimation scheme \revA{shows promising results} in a joint \revA{refinement} scheme.

For question~2, according to theory, \revA{i.e., \cref{thm:new-num-non-zero-coeffs}},  the amount of non-zero coefficients should lie between 8 and 15\revA{, if we had $\regParam$ in closed form rather than the approximation from \cref{prop:choice-regulariser-lifting}}. From the figure we observe that the upper bound is satisfied at all times, but that the lower bound is violated at all the iterates. The latter is most likely \revA{a sign that the approximation of $\regParam$ through \cref{prop:choice-regulariser-lifting} for the current} size of the sampling set \revA{suffers from a notable but ever so slight error}. Overall, the results \revA{from the approximated $\regParam$} are reasonably in line with \cref{thm:new-num-non-zero-coeffs} and from a practical point of view \cref{prop:choice-regulariser-lifting} \revA{is a suitable approximator}.

Finally, it should be mentioned that the resolution of the obtained 3D maps after 10 iterations is visibly lower than the ground truth. The reason is twofold. First, with the small data size we cannot expect very high resolution, but -- more importantly -- in real-life data, a defocus depth \emph{range} is used rather than a single defocus depth. This is an important detail as this prevents losing certain frequency bands in Fourier space on the zero-crossings of the CTF. However, we remind the reader that this work focuses on rotation estimation and that such considerations are beyond the scope of this work.






\section{Conclusions}
\label{sec:conclusions}
In this work we presented ellipsoidal support lifting (ESL) as a general lifting-based algorithm for approximating \revA{and regularising} the global minimiser of a smooth function over a Riemannian manifold.
The main theoretical result is that the algorithm is well-posed; we also provided upper bounds on the distance to the minimiser. 
Additionally, several minor results were provided that can be used to make the algorithm more practically useful.

\revA{The asymptotic behaviour predicted by the} theoretical results has been \revA{observed} in several numerical experiments related to single particle Cryo-EM \revA{-- being the main motivation for the design of the proposed method.} Not only do these experiments support the theory, they also showcase the practical usefulness of ESL. 
\revA{That is, introducing bias through the regularising effect of ESL -- rather than actual global optimisation -- yields rotations closer to the ground truth}. \revA{Moreover, when integrated in a joint refinement scheme, ESL shows promising results when compared to state-of-the-art solvers in terms of rotation estimation.}





Finally, we see further potential of the proposed theory in the context of single particle Cryo-EM as follows:



\paragraph{Symmetric molecules}
Symmetries do not pose an issue in our framework regarding non-uniqueness of minimisers, the reason being that such symmetries can be quotiented out and our results still hold on the quotient manifold. 
As an example, assume the molecule 3D map $\vol$ has a rotation symmetry  generated by some finite subgroup $\mathrm{G}\subset \mathrm{SO}(3)$.
We can then consider the manifold $\mathrm{SO}(3)/\mathrm{G}$ and adapt \cref{alg:joint-refinement-lifting} for \revA{formally} solving
\begin{align}
\label{eq:unlifted-problem-cryo-vol-rot-symmetrised}
    \inf_{\substack{\vol \in \volLtwoSpace \\ \rot_\numImgs \in \mathrm{SO}(3)/\mathrm{G}}} \biggl\{ \sum_{\imgInd=1}^\numImgs \Bigl[ \frac{1}{|\mathrm{G}|}\sum_{\rot \in \mathrm{G}} \frac{1}{2\noiseLevelParam} \bigl\|\forward ((\rot\circ\rot_{\imgInd}).\vol)-\img_{\imgInd}\bigr\|_2^{2} \Bigr] +\regA(\vol) \biggr\}
    \quad\text{for  $\noiseLevelParam>0$.}
\end{align}
Additionally, our theoretical bounds suggest that in the case of $|\mathrm{G}|$-fold symmetry we would only need $1/|\mathrm{G}|$ times the amount of rotations in the sampling set $\rotsSampling$ to get the same error bound as for a non-symmetric molecule.

\paragraph{In-plane shifts}
In experimental 3D map and orientation reconstruction, one also considers 2D in-plane shifts of the molecule due to improper alignment of the particles, i.e., the actual manifold of interest is a bounded subset of $\mathrm{SO}(3)\times \Real^2$. 
\Cref{alg:joint-refinement-lifting} does not change except that our theoretical results suggest that one should choose $\exponentA\in (\frac{1}{6},\frac{2}{5})$ since the manifold $\mathrm{SO}(3)\times \Real^2$ is now five-dimensional.

\paragraph{Local search}
The theory suggests that instead of increasing the amount of points in the sampling set it can be equally efficient to reduce the volume of the manifold, i.e., the search area. 
This is already done in RELION as a local search and in cryoSPARC through the multi-scale nature of branch-and-bound. 
Our framework would also be suitable here as bounded sub-manifolds satisfy the conditions for the theory to hold.



%
%
%

\section*{Acknowledgments}
We would like to thank Sjors Scheres, Dari Kimanius and Johannes Schwab for fruitful discussions throughout this project \revA{and for assistance with the comparison to RELION.}

\bibliography{bibliography}

\appendix
\section{Local low discrepancy on $(0,1)\subset \Real$}
\label{sec:low-discr-R}
\revA{
For fixed $0<\exponentA<2$ and $\dummyIntegerB>0$, consider the mapping $\radiusB:\Natural \to \Real$ given by 
\begin{equation}
    \radiusB(M):=\frac{\dummyIntegerB}{2} M^{-\frac{1+\exponentA}{3}}.
    \label{eq:rm-local-R}
\end{equation}
and construct the sequence $(M_m)_{m=1}^\infty\subset \Natural$, where $M_1 = 1$ and for $m>1$ the number $M_m$ is defined as the first natural number larger than $M_{m-1}$ satisfying
\begin{equation}
    2\radiusB(M_m-1) \cdot M_m \leq  \lfloor 2\radiusB(M_m) \cdot(M_m+1) \rfloor < 2\radiusB(M_m) \cdot (M_m + 1)
    \label{eq:sequence-len-construction-R}
\end{equation}
This sequence is well-defined because the mapping $M \mapsto 2\radiusB(M) \cdot (M+1)$ is monotonically increasing and unbounded for our choice of $\exponentA>0$.


\begin{proposition}
    Let $0<\exponentA<2$ and let $\dummyIntegerB_0>0$. Consider the sequence of sampling sets $(\rotsSampling_m)_{m=1}^\infty \subset (0,1)$ given by
    \begin{equation}
        \rotsSampling_m := \Bigl\{\frac{i}{M_m+1} \, \mid \, i = 1, \ldots, M_m \Bigr\},
    \end{equation}
    where $(M_m)_{i=1}^\infty\subset \Natural $ as constructed above through $\dummyIntegerB$ and $\exponentA$.
    
    Then, $(\rotsSampling_m)_{m=1}^\infty$ is a local low discrepancy sequence with $\exponentA$ and $\dummyIntegerB$ on $\manifold := (0,1)$.
\end{proposition}
\begin{proof}
    Fix $0<\exponentA<2$ and $\dummyIntegerB >0$ and define the mapping $\radiusB:\Natural \to \Real$ as in \cref{eq:rm-local-R}. Consider the radius $\radiusA_\sumTotA$ from \cref{def:eta-local-low-discrepancy-sequence}, which is given by 
    \begin{equation}
        \radiusA_\sumTotA :=  \sqrt[1]{\dummyIntegerB \frac{\operatorname{vol}((0,1))}{\ballVol_1}} |\rotsSampling_m|^{- \frac{1 + \exponentA}{1 + 2}} = \frac{\dummyIntegerB}{2} M_m^{-\frac{1+\exponentA}{3}}.
    \end{equation}
    Note that $\radiusB(M_\sumTotA)  = \radiusA_\sumTotA$. Next, choose a point $\mPoint\in (0,1)$ and a positive definite bi-linear form $\dummyBilForm: \Real\times\Real\to\Real$ given by $\dummyBilForm(x,y):= a xy$ for some $a>0$. Then, $\dummyBilForm(\log_p q,\log_p q) = a(q - p)^2$. The following proof consists of two steps: (i) show \cref{eq:eta-local-low-discrepancy-sequence-1} and (ii) show \cref{eq:eta-local-low-discrepancy-sequence-2}.

    (i) It is easy to check that 
    \begin{equation}
        |\ellipsoid^{\dummyBilForm}_{\radiusA_\sumTotA}(\mPoint) \cap \rotsSampling_m| = \lfloor 2\radiusB(M_m)(M_m+1) \rfloor.
    \end{equation}
    Subsequently, we have
    \begin{equation}
        \frac{\operatorname{vol}(\ellipsoid^{\dummyBilForm}_{\radiusA_\sumTotA}(\mPoint))}{\operatorname{vol}((0,1))} = 2\radiusA_\sumTotA = 2\radiusB(M_m).
    \end{equation}
    Then, for showing \cref{eq:eta-local-low-discrepancy-sequence-1}, i.e., 
    \begin{equation}
        \biggl| \frac{|\ellipsoid_{\radiusA_\sumTotA}^\dummyBilForm(\mPoint) \cap \rotsSampling_\sumTotA|}{|\rotsSampling_\sumTotA|} - \frac{\operatorname{vol}\bigl(\ellipsoid_{\radiusA_\sumTotA}^\dummyBilForm(\mPoint)\bigr)}{\operatorname{vol}(\manifold)} \biggr|
    \in o\Bigl(|\rotsSampling_\sumTotA|^{-\frac{\DimInd(1+\exponentA)}{\DimInd+2}}\Bigr), 
    \end{equation}
    it suffices to show that the right-hand inequality in the following equation holds for large $m$:
    \begin{equation}
        |\frac{\lfloor 2\radiusB(M_m)(M_m+1) \rfloor}{M_m} - 2\radiusB(M_m)| = \frac{1}{M_m}|\lfloor 2\radiusB(M_m)(M_m+1) \rfloor - 2\radiusB(M_m)M_m| \leq \frac{4\radiusB(M_m)}{M_m}
        \label{eq:local-1-property-R}
    \end{equation}
    as $\frac{\radiusB(M_m)}{M_m}\in \mathcal{O}(M_m^{-\frac{2}{3}(1+\exponentA)}) \overset{M_m = |\rotsSampling_m|}{=} \mathcal{O}(|\rotsSampling_m|^{-\frac{2}{3}(1+\exponentA)})$ holds due to the assumption $\exponentA<2$ and $\mathcal{O}(|\rotsSampling_m|^{-\frac{2}{3}(1+\exponentA)}) \subset o(|\rotsSampling_m|^{-\frac{1}{3}(1+\exponentA)})$ due to the assumption $\exponentA >0$. 
    
    Therefore, it remains to show the inequality in \cref{eq:local-1-property-R}. Consider that
    \begin{multline}
        |\lfloor 2\radiusB(M_m)(M_m+1) \rfloor - 2\radiusB(M_m)M_m| \\
        \leq \max (| 2\radiusB(M_m)(M_m+1) - 2\radiusB(M_m)M_m|,|2\radiusB(M_m-1)M_m  - 2\radiusB(M_m)M_m|)
    \end{multline}
    by \cref{eq:sequence-len-construction-R}. In other words, if both
    \begin{equation}
        |2\radiusB(M_m)(M_m+1)  - 2\radiusB(M_m)M_m| \leq 4\radiusB(M_m)
        \label{eq:local-1-property-R-upper}
    \end{equation}
    and 
    \begin{equation}
        |2\radiusB(M_m-1)M_m  - 2\radiusB(M_m)M_m| \leq 4\radiusB(M_m)
        \label{eq:local-1-property-R-lower}
    \end{equation}
    hold, \cref{eq:local-1-property-R} holds and we are done.
    
    We start by showing \cref{eq:local-1-property-R-upper}. This bound actually follows straight away. Indeed,
    \begin{equation}
        |2\radiusB(M_m)(M_m+1)  - 2\radiusB(M_m)M_m|  =  2\radiusB(M_m) \leq 4\radiusB(M_m).
    \end{equation}
    So it remains to show \cref{eq:local-1-property-R-lower}. For that we will rewrite the condition several times. Dividing by $2M_m$ on both sides and using that $\radiusB(M_m-1) \geq \radiusB(M_m)$ for all $m$ gives the equivalence
    \begin{equation}
        |2\radiusB(M_m-1)M_m  - 2\radiusB(M_m)M_m| \leq 4\radiusB(M_m) \quad \Leftrightarrow \quad \radiusB(M_m-1)  - \radiusB(M_m) \leq \frac{2\radiusB(M_m)}{M_m}.
        \label{eq:lower-first-eqi-R}
    \end{equation}
    Then,
    \begin{multline}
        \radiusB(M_m-1)  - \radiusB(M_m) 
        \leq \frac{2\radiusB(M_m)}{M_m} \quad 
        \Leftrightarrow \quad  \radiusB(M_m-1) \leq \Bigl(1 + \frac{2}{M_m}\Bigr)  \radiusB(M_m)\\
        \overset{\cref{eq:rm-local-R}}{\Leftrightarrow} \quad (M_m - 1)^{-\frac{2}{3}(1 + \exponentA)} \leq \Bigl(1 + \frac{2}{M_m}\Bigr)  M_m^{-\frac{2}{3}(1 + \exponentA)} 
        \quad
        \Leftrightarrow \quad \Bigl(\frac{M_m}{M_m - 1}\Bigr)^{\frac{2}{3}(1 + \exponentA)} \leq 1 + \frac{2}{M_m}\\
        \Leftrightarrow \quad  1 + \frac{2}{M_m} - \Bigl(\frac{M_m}{M_m - 1}\Bigr)^{\frac{2}{3}(1 + \exponentA)} \geq 0.
        \label{eq:lower-R-equiv}
    \end{multline}
    Notice that
    \begin{equation}
        \lim_{m\to \infty}  1 + \frac{2}{M_m} - \Bigl(\frac{M_m}{M_m - 1}\Bigr)^{\frac{2}{3}(1 + \exponentA)} = 0.
    \end{equation}
    So if we can show that the sequence decreases to zero for all $m\geq m'$, the equivalences in \cref{eq:lower-R-equiv} yield the lower bound \cref{eq:local-1-property-R-lower} through \cref{eq:lower-first-eqi-R}.
    
    Let $f:(1,\infty) \to \Real$ given by 
    \begin{equation}
        f(x):=1 + \frac{2}{x} - \Bigl(\frac{x}{x-1}\Bigr)^{\frac{2}{3}(1 + \exponentA)}
    \end{equation}
    Showing that $f$ is monotonically decreasing (to $0$) for large enough $x$ proofs \cref{eq:lower-R-equiv} and with that our claim. And indeed,
    \begin{align}
    \begin{split}
         \frac{\mathrm{d} f}{\mathrm{d}x} \leq 0 \!\!\quad \!\!\Leftrightarrow \!\! \quad \!\!
         - \frac{2}{x^2} + \frac{2}{3}(1 + \exponentA) \frac{\bigl(\frac{x}{x-1} 
         \bigr)^{\frac{2}{3}(1 + \exponentA)-1}}{(x-1)^2} \leq 0 \!\! \quad \!\!
        \Leftrightarrow \!\! \quad  \!\! \Bigl( \frac{x}{x-1}\Bigr)^{\frac{2}{3}(1 + \exponentA)+1} \leq \frac{2}{\frac{2}{3}(1 + \exponentA)} = \frac{3}{1 + \exponentA}.
        \label{eq:-lower-deriv-R}
    \end{split}
    \end{align}
    Note that the right-most term in the last line is larger than 1 because $\exponentA <2$ and the left-hand-side of the inequality is decreasing to 1, which ensures the existence oa an $x'$ so that \cref{eq:-lower-deriv-R} holds for all $x\geq x'$. Consequently, there exists an $m'\in \Natural$ such that for all $m\geq m'$ we have $1 + \frac{2}{M_m} - \bigl(\frac{M_m}{M_m - 1}\bigr)^{\frac{2}{3}(1 + \exponentA)} \geq 0$. We conclude that \cref{eq:eta-local-low-discrepancy-sequence-1} holds for this special choice of $\rotsSampling_m$.
    
    \vspace{0.2cm}
    (ii) For the integration of quadratic forms we aim to prove \cref{eq:eta-local-low-discrepancy-sequence-2}, i.e., bounding the term
    
    \begin{multline}
        \biggl|\sum_{\mPointB\in \ellipsoid_{\radiusA_\sumTotA}^\dummyBilForm(\mPoint) \cap \rotsSampling_\sumTotA}\!\!\!\!
        \frac{\dummyBilForm\bigl(\log_\mPoint(\mPointB),\log_\mPoint(\mPointB)\bigr)}{|\rotsSampling_\sumTotA|} - \frac{\int_{\ellipsoid_{\radiusA_\sumTotA}^\dummyBilForm(\mPoint)}\dummyBilForm\bigl(\log_\mPoint(\mPointB),\log_\mPoint(\mPointB)\bigr)\; \mathrm{d}\mPointB}{\operatorname{vol}(\manifold)} \biggr| \\
        = \biggl|\sum_{\mPointB\in \ellipsoid_{\radiusA_\sumTotA}^\dummyBilForm(\mPoint) \cap \rotsSampling_\sumTotA}\!\!\!\!
        \frac{a(\mPointB - \mPoint)^2}{|\rotsSampling_\sumTotA|} - \int_{\mPoint - \radiusA_\sumTotA}^{\mPoint + \radiusA_\sumTotA} a(\mPointB - \mPoint)^2\; \mathrm{d}\mPointB \biggr|. 
        \label{eq:quadratic-approx-R}
    \end{multline}
    Note that for large enough $\sumTotA\in\Natural$ there exist $\radiusA_\sumTotA^-, \radiusA_\sumTotA^+ \in (\radiusA_\sumTotA - \frac{1}{M_\sumTotA +1}, \radiusA_\sumTotA + \frac{1}{M_\sumTotA +1})$ such that the discrete sum is a mid-point integration approximation of $\int_{\mPoint - \radiusA_\sumTotA^-}^{\mPoint + \radiusA_\sumTotA^+} a(\mPointB - \mPoint)^2\; \mathrm{d}\mPointB$, which allows to use its well-known error bound 
    \begin{equation}
        \biggl|\sum_{\mPointB\in \ellipsoid_{\radiusA_\sumTotA}^\dummyBilForm(\mPoint) \cap \rotsSampling_\sumTotA}\!\!\!\!
        \frac{a(\mPointB - \mPoint)^2}{|\rotsSampling_\sumTotA|} - \int_{\mPoint - \radiusA_\sumTotA^-}^{\mPoint + \radiusA_\sumTotA^+} a(\mPointB - \mPoint)^2\; \mathrm{d}\mPointB \biggr| \leq  \frac{2a (\radiusA_\sumTotA^- + \radiusA_\sumTotA^+)^3}{24 |\ellipsoid_{\radiusA_\sumTotA}^\dummyBilForm(\mPoint) \cap \rotsSampling_\sumTotA|} .
        \label{eq:midpoint-int-r1}
    \end{equation}
    In order to also bound the effect of this slight modification of the integral bounds, we use the triangle inequality to \cref{eq:quadratic-approx-R},
    \begin{multline}
        \biggl|\sum_{\mPointB\in \ellipsoid_{\radiusA_\sumTotA}^\dummyBilForm(\mPoint) \cap \rotsSampling_\sumTotA}\!\!\!\!
        \frac{a(\mPointB - \mPoint)^2}{|\rotsSampling_\sumTotA|} - \int_{\mPoint - \radiusA_\sumTotA}^{\mPoint + \radiusA_\sumTotA} a(\mPointB - \mPoint)^2\; \mathrm{d}\mPointB \biggr| \\
        \leq \biggl| \int_{\mPoint - \radiusA_\sumTotA}^{\mPoint - \radiusA_\sumTotA^-} a(\mPointB - \mPoint)^2\; \mathrm{d}\mPointB\biggr| + \biggl|\sum_{\mPointB\in \ellipsoid_{\radiusA_\sumTotA}^\dummyBilForm(\mPoint) \cap \rotsSampling_\sumTotA}\!\!\!\!
        \frac{a(\mPointB - \mPoint)^2}{|\rotsSampling_\sumTotA|} - \int_{\mPoint - \radiusA_\sumTotA^-}^{\mPoint + \radiusA_\sumTotA^+} a(\mPointB - \mPoint)^2\; \mathrm{d}\mPointB \biggr| + \biggl| \int_{\mPoint + \radiusA_\sumTotA^+}^{\mPoint + \radiusA_\sumTotA} a(\mPointB - \mPoint)^2\; \mathrm{d}\mPointB \biggr|
        \label{eq:three-term-split-local-low-R}
    \end{multline}
    For showing that \cref{eq:eta-local-low-discrepancy-sequence-2} holds, we will show that the upper bound in \cref{eq:three-term-split-local-low-R} is $o(|\rotsSampling_\sumTotA|^{-(1+\exponentA)})$.
    
    Using the fact that by definition $|\radiusA^{\pm}_m - \radiusA_m| < \frac{1}{M_m +1}$, and that $\frac{1}{M_\sumTotA +1} \leq \radiusA_\sumTotA$ for $\exponentA<2$, the two outer-most terms can be bounded by
    \begin{equation}
        \biggl| \int_{\mPoint - \radiusA_\sumTotA}^{\mPoint - \radiusA_\sumTotA^-} a(\mPointB - \mPoint)^2\; \mathrm{d}\mPointB\biggr|  \leq \frac{a(2\radiusB(M_\sumTotA))^2}{M_\sumTotA +1}\leq \frac{4a \radiusB(M_\sumTotA)^2}{M_\sumTotA}, 
        \label{eq:outer-bound-rewrite}
    \end{equation}
    and
    \begin{equation}
        \biggl| \int_{\mPoint + \radiusA_\sumTotA^+}^{\mPoint + \radiusA_\sumTotA} a(\mPointB - \mPoint)^2\; \mathrm{d}\mPointB \biggr|  \leq \frac{a(2\radiusB(M_\sumTotA))^2}{M_\sumTotA +1}\leq \frac{4a \radiusB(M_\sumTotA)^2}{M_\sumTotA}. 
        \label{eq:outer-bound-rewrite-u}
    \end{equation}
    
    For the middle term in \cref{eq:three-term-split-local-low-R} we have already shown \cref{eq:midpoint-int-r1}, which can be further bounded by
    \begin{multline}
        \frac{2a (\radiusA_\sumTotA^- + \radiusA_\sumTotA^+)^3}{24 |\ellipsoid_{\radiusA_\sumTotA}^\dummyBilForm(\mPoint) \cap \rotsSampling_\sumTotA|} \overset{\frac{1}{M_\sumTotA +1} \leq \radiusA_\sumTotA}{\leq} \frac{2a (4\radiusA_\sumTotA)^3}{24 |\ellipsoid_{\radiusA_\sumTotA}^\dummyBilForm(\mPoint) \cap \rotsSampling_\sumTotA|} \overset{\cref{eq:sequence-len-construction-R}}{\leq} \frac{2a (4\radiusA_\sumTotA)^3}{24 (2\radiusB(M_m-1) M_m)} \\
        \overset{\substack{M\mapsto \radiusB(M)\\ \text{is decreasing}}}{\leq}  \frac{2a (4\radiusA_\sumTotA)^3}{24 (2\radiusB(M_m) M_m)} = \frac{2a (4\radiusA_\sumTotA)^3}{24 (2\radiusA_\sumTotA M_m)} = \frac{8a \radiusA_\sumTotA^2}{3 M_m} = \frac{8a \radiusB(M_\sumTotA)^2}{3 M_m}
        \label{eq:mid-bound-rewrite}
    \end{multline}
    
    Finally, because $\frac{\radiusB(M_\sumTotA)^2}{M_m} \in o(M_\sumTotA^{-(1+\exponentA)}) = o(|\rotsSampling_\sumTotA|^{-(1+\exponentA)})$
    for $\exponentA<2$, we conclude that \cref{eq:eta-local-low-discrepancy-sequence-2} holds by inserting \cref{eq:outer-bound-rewrite,eq:outer-bound-rewrite-u,eq:mid-bound-rewrite} into \cref{eq:three-term-split-local-low-R}, which was the second and final part of the claim.
\end{proof}

}
\section{The Cryo-EM forward operator}
\label{sec:maths-of-cryo-em}
When imaging thin unstained cryofixated biological specimens in a high resolution TEM, one can assume that the weak-phase object approximation holds \cite{vulovic2014use,oktem2015mathematics}. 
This approximation allows us to express $\forward$ as a parallel beam ray transform on lines parallel to the TEM optical axis (tomographic projection) followed by a 2D convolution in the detector plane to account for TEM optics and detector response \cite{oktem2015mathematics}. 

Phrasing the above mathematically, we get that $\forward \colon \volSpace \to \imgSpace $ in \cref{eq:CryoEMInvProb} is of the form
\begin{equation}
    \forward(\vol) := \RerCTF \ast \bigl(\projection(\vol)\bigr)
    \quad\text{for $\vol \in \volSpace.$}
    \label{eq:cryo-operator}
\end{equation}
In the above, $\projection \colon \volSpace \to \imgSpace$ is the parallel beam ray transform (projection mapping) restricted to lines parallel to the TEM optical axis.
Choosing a coordinate system where the $\eC$-axis is parallel to the TEM optical axis allows us to express it as
\begin{equation}
    \projection \vol (\eA,\eB) := \int_{-\infty}^\infty \vol(\eA,\eB,\eC) \mathrm{d} \eC.
\end{equation}
Next, the convolution in \cref{eq:cryo-operator} is a 2D convolution in the detector plane and the corresponding convolution kernel $\RerCTF \colon \Real^2 \to \Real$ is given analytically by its Fourier transform:
\begin{equation}
    \ReCTF(\xi) := -  A(\xi) \biggl(\sqrt{1 - \ampContrast^2} \sin\Bigl( W\bigl( |\xi|\bigr)\Bigr) 
    + \alpha \cos\Bigl( W\bigl( |\xi|\bigr)\Bigr)
    \biggr)
    \label{eq:ctf}
\end{equation}
with $0 < \ampContrast < 1$ denoting the \emph{amplitude contrast ratio}, $A \colon \Real^2 \to \Real$ is the aperture function that is commonly an indicator function on a disc with a radii given by the TEM objective aperture, and
\[ W(s) := \dfrac{\defocus}{2 k} s^{2}-\dfrac{\spAberration}{4 k^{3}} s^{4} \] models imperfections in the TEM optics. Here, $\defocus$ is the \emph{defocus}, $\spAberration$ is the \emph{spherical aberration}, and $\waveNumber$ is the (relativistically corrected) \emph{wave number} of the imaging electron. The function $\ReCTF$ in \cref{eq:ctf} is often referred to as the \emph{contrast transfer function} (CTF).

\end{document}